\documentclass[11pt, reqno]{amsart}

\usepackage{geometry}
\usepackage[utf8]{inputenc}
\usepackage{amstext,latexsym,amsbsy,amsmath,amssymb,amsthm,mathtools,relsize,geometry,enumerate}
\usepackage{hyperref}
\hypersetup{
	colorlinks   = true, 
	urlcolor     = blue, 
	linkcolor    = red, 
	citecolor   = red 
}
\usepackage{mathtools,enumerate,enumitem}

\usepackage{mathrsfs}

\usepackage{comment}
\usepackage{graphicx}
\usepackage{epstopdf}
\setkeys{Gin}{width=\linewidth,totalheight=\textheight,keepaspectratio}
\graphicspath{{./images/}}

\usepackage[normalem]{ulem}

\usepackage{multicol}
\usepackage{booktabs}
\usepackage{mathrsfs}
\usepackage{color}

\numberwithin{equation}{section}

\usepackage{bbm}

\newcommand{\E}{\mathcal{E}}

\newcommand{\mP}{\mathbb{P}}

\newcommand{\mE}{\mathbb{E}}

\newcommand{\Ome}{\Omega}
\newcommand{\p}{\partial}

\def\E{\mathbb{E}}
\def\P{\mathbb{P}}

\newcommand{\f}{\varphi}
\newcommand{\Zbb}{\mathbb{Z}}

\numberwithin{equation}{section}

\theoremstyle{plain}
\newtheorem{theorem}{Theorem}[section]

\newtheorem{remark}{Remark}[section]
\newtheorem{lemma}[theorem]{Lemma}

\allowdisplaybreaks[4]

\begin{document}
	
	\title{Fully discrete finite element methods for the stochastic Kuramoto-Sivashinsky equation with multiplicative noise}
	
	\author{
		Hung D. Nguyen$^1$ and 
		Liet Vo$^2$ 
	}
	
	\thanks{$^1$ Department of Mathematics, The University of Tennessee, Knoxville, TN 37996, U.S.A. ({\tt hnguye53@utk.edu}).} 
	
	\thanks{$^2$ School of Mathematical and Statistical Sciences, The University of Texas Rio Grande Valley, Edinburg, TX 78539, U.S.A.  ({\tt liet.vo@utrgv.edu}). This author was partially supported by the NSF grant DMS-2530211.}

		\begin{abstract} We investigate a fully discrete finite element approximation for the stochastic Kuramoto–Sivashinsky equation, combining the standard finite element methods in spatial discretization with the implicit Euler–Maruyama scheme in time. Rigorous error estimates are established for two distinct noise regimes. In the case of bounded multiplicative noise, we prove optimal strong convergence rates in full expectation. The analysis relies crucially on a stochastic Gronwall inequality and an exponential stability estimate for the PDE solution, which together control the interplay between the nonlinear drift and the multiplicative stochastic forcing. For general multiplicative noise, where boundedness no longer holds, we derive sub-optimal convergence rates in probability by introducing a localization technique based on carefully constructed subsets of the sample space. This dual framework demonstrates that the proposed fully discrete scheme achieves strong convergence under bounded noise and probabilistic convergence under general multiplicative noise, thus providing the first comprehensive error analysis for numerical approximations of the stochastic Kuramoto–Sivashinsky equation. {Numerical experiments are also provided to demonstrate the efficiency of the numerical method and validate the theoretical results.}
	\end{abstract}

	\maketitle

	
	\medskip
	

	\section{Introduction}\label{sec-1}
We consider the stochastic Kuramoto--Sivashinsky (SKS) equation with It\^o multiplicative noise
\begin{subequations}\label{eq1.1}
	\begin{align} \label{eq1.1a}
		du  + \big(\nu \partial_x^4 u + \partial_x^2 u + u\partial_x u\big)\,dt &=  B(u)\,dW(t), 
		&&\mbox{in }(0,T)\times D,\,a.s.\\ 
		u(0) &= u_0, &&\mbox{on } D,\, a.s.,
	\end{align}
\end{subequations}
where $D=[0,L]\subset \mathbb{R}$ with $L>0$, and the solution is taken to be periodic in space with period $L$. 
Here, $\nu>0$ is a constant representing the diffusion effect, $\{W(t): t\geq 0\}$ denotes a real-valued Wiener process, and the operator $B$ depends on the solution $u$, leading to a multiplicative noise structure. 
A precise definition of $B$ and the assumptions imposed on it are given in Section \ref{sec2}.  

The deterministic Kuramoto--Sivashinsky (KS) equation has long been regarded as a canonical model in the study of dissipative nonlinear systems that exhibit spatio-temporal chaos, turbulence, and pattern formation. Originally introduced in the context of laminar flame-front propagation \cite{kuramoto1978diffusion, michelson1977nonlinear,ashinsky1988nonlinear}, the KS equation has also appeared in thin-film dynamics and plasma physics (see \cite{adams2025well, cohen1976non,laquey1975nonlinear} and references therein). 
It has become a standard testbed for both analytical and numerical methods in nonlinear PDEs. 
When stochastic forcing is introduced, the SKS equation serves as a model for complex physical systems subject to random perturbations, such as thin-film growth under noise or noisy reaction--diffusion processes. 

From a mathematical perspective, the stochastic Kuramoto--Sivashinsky (SKS) equation poses significant challenges due to the interplay of several competing mechanisms. The nonlinear convective term $u\partial_x u$ is non-Lipschitz,  preventing the direct application of standard SPDE techniques that assume Lipschitz drift globally. 
This term interacts strongly with the multiplicative stochastic noise, where fluctuations depend on the solution itself and can be amplified in regions of large amplitude, potentially causing rapid growth or instability. 
At the same time, the underlying PDE operator $\nu\partial_x^4 +\partial_x^2$ is not guaranteed to be sign-definite for any value of $\nu$, which complicates stability analysis and undermines standard coercivity arguments. These three effects do not act in isolation, but interact in subtle and mutually reinforcing ways. Controlling this interaction is the central difficulty in both the analytical study and the design of convergent numerical schemes for the SKS equation.

The well-posedness of the SKS equation has been studied in the past decades. Duan and Ervin \cite{duan2001stochastic} proved global existence and uniqueness under additive noise, while Wu, Cui, and Duan \cite{wu2018global} extended these results to multiplicative noise. The work of Ferrario \cite{ferrario2008invariant} investigated the long-time behaviors of the SKS through the analysis of statistically steady states.
More recently, Gao and Nguyen \cite{gao2025exponential} showed that the SKS equation with additive noise exhibits an exponential stability estimate for the analytical solution. These results provide a solid analytical foundation that motivates the numerical analysis developed in this work. 

Over the last three decades, intensive numerical approaches have been developed for the deterministic KS equation, including finite different method \cite{mackenzie2000holistic}, local discontinuous galerkin method \cite{xu2006local}, spectral collocation method \cite{khater2008numerical, manickam1998second}, Implicit–Explicit BDF method \cite{akrivis2004implicit}, pseudospectral method \cite{lopez1994numerical}, meshless method \cite{haq2010meshless} and finite element method \cite{akrivis1994finite, anders2012higher, doss2019fourth}.  
For the stochastic KS equation, much of the literature has focused on simulation and qualitative behavior rather than on rigorous numerical analysis. 
For example, spectral-based schemes have been employed in \cite{drotar1999numerical,ueno2005renormalization} to study the scaling behavior of surface growth and the large-scale and long-time behavior of the SKS equation, but without providing sharp convergence guarantees. 
To the best of our knowledge, no comprehensive finite element error analysis for the SKS equation with multiplicative noise is currently available.  

By contrast, numerical analysis of other stochastic PDEs is more advanced. For semilinear parabolic SPDEs, the convergence of finite element methods and implicit Euler--Maruyama schemes has been studied in works such as \cite{printems2001discretization, vo2023higher,feng2021optimally, carelli2012rates,breit2021convergence}. These results, however, do not directly apply to the SKS equation because of its higher-order derivatives, the delicate balance of stability and instability, and the nontrivial interaction of multiplicative noise with the nonlinear drift. The purpose of this paper is therefore to fill this gap by developing a rigorous convergence theory for a fully discrete finite element approximation of the SKS equation. 

Our approach combines a standard continuous finite element method for spatial discretization with the implicit Euler--Maruyama scheme for temporal discretization. For bounded multiplicative noise, we establish optimal strong convergence rates in full expectation, relying on three key analytical tools: a discrete stochastic Gronwall inequality \cite{kruse2018discrete}, an exponential stability estimate for the exact solution, and a higher-moment bootstrap technique developed in \cite{vo2025high, vo2024high}. Together, these allow us to control the nonlinear drift--noise interaction and recover sharp error bounds. It should be noted that the idea has been successfully applied to establish the full moment error estimates of a time-discrete scheme for the stochastic Navier-Stokes equations \cite{feng2025NS}. For the more general case of unbounded multiplicative noise, where exponential stability estimates are no longer available, we employ a localization strategy based on subsets of the sample space, following the approach of \cite{carelli2012rates}. This yields sub-optimal error estimates, with expectations taken over the localized subsets, thereby extending the applicability of our analysis to a broader class of noise structures. Beyond the SKS equation itself, the techniques developed here, in particular, the combination of exponential stability estimates with stochastic Gronwall inequalities and bootstrapping arguments, are expected to be useful in the analysis of other nonlinear SPDEs with non-Lipschitz drift and multiplicative noise. Lastly, we remark that the results obtained in this work for the SKS are only valid in dimension one. In fact, for arbitrary $\nu>0$, the global well-posedness of strong solutions to the deterministic KS equation in higher dimensions is not yet completely resolved \cite{larios2024remarks}. 

The remainder of the paper is organized as follows. 
In Section~\ref{sec2}, we introduce the functional framework, assumptions on the noise operator, and the regularity of the PDE solution, including presenting an exponential stability estimate, high moment H\"older continuity estimates, and high moment stability estimates in $H^m$ norm for $m \geq 0$.
Section~\ref{sec3} describes the fully discrete finite element method, and its stability estimates. In particular, Section~\ref{sub-sec3.2} is devoted to the error analysis under bounded multiplicative noise, while Section~\ref{sub-sec3.3} extends the results to general multiplicative noise using localization. Numerical experiments are presented in Section~\ref{sec-4} to confirm the efficiency of the proposed method and to validate the theoretical results. Concluding remarks and possible extensions are given in Section~\ref{sec-5}. Finally, we present and prove some useful auxiliary results in Appendix \ref{auxiliary}.

\section{Preliminaries}\label{sec2}
\subsection{Notations and assumptions}
Standard function and space notation will be adopted in this paper. 
We denote ${L}^p(D)$ and ${ H}^{k}(D)$ as the Lebesgue and Sobolev spaces of the functions that are periodic with period $L$ and have vanishing mean. In particular, $L^2(D)=H^0(D)$.
$C$ denotes a generic constant that is independent of the mesh parameters $h$ and $k$.

In this paper, we will make the following assumptions on $B: L^2(D) \rightarrow L^2(D)$.
\begin{enumerate}[label=(\Alph*)]
	\item There exists a constant $L_0>0$ such that
	\begin{align}\label{Assump_Lineargrowth}
		\|B(u)\|_{L^2} \leq L_0.
	\end{align}
	\item For all $m\ge 0$, suppose that $B: H^{m}(D) \rightarrow H^m(D)$. Moreover, there exists a constant $C_{B}>0$, {depending on $m$}, such that
	\begin{align}\label{Assump_Lipschitz}
		\|B(u) - B(v)\|_{H^m} \leq C_{B}\|u-v\|_{H^m}.
	\end{align}
	It should be noted that \eqref{Assump_Lipschitz} implies
	\begin{align} \label{cond:|B(u)|_Hm<|u|_Hm}
		\|B(u)\|_{H^m} \leq c_B (1 + \|u\|_{H^m}),
	\end{align}
	where $c_B = \max\{C_B, \|B(0)\|_{H^m}\}$.
\end{enumerate}
In addition, we recall the Sobolev embedding inequality in one dimension:
\begin{align}\label{sobolev_ineq}
	\|u\|_{L^{\infty}} \leq C_e\|u\|_{H^1},
\end{align}
where $C_e>0$ is a constant, {which depends only on the domain $D$}.

\begin{remark}
	1. We note that while the Lipschitz condition \eqref{Assump_Lipschitz} will be assumed throughout, the boundedness condition \eqref{Assump_Lineargrowth} is only needed to derive an exponential moment estimate for the solutions of \eqref{eq1.1}; see Lemma \ref{lemma_expo_moment_H}. This result is then used in Section \ref{sub-sec3.2} to establish convergence in $L^p$ for the finite element scheme.
	
	{ 2. In this paper, we consider $W$ to be a standard one-dimensional Brownian motion and $B$ to be an operator on $L^2(D)$ for notational simplicity. The analysis can be extended to the case of a cylindrical Wiener process or a $Q$-Wiener process with Hilbert--Schmidt operators $B$ satisfying assumptions analogous to \eqref{Assump_Lipschitz} and \eqref{Assump_Lineargrowth}. In this more general setting, the main additional effort lies in handling the noise terms in the error estimates, in particular controlling the corresponding Hilbert--Schmidt norms. However, these difficulties are mild and do not affect the overall structure of the analysis. For instance, a $Q$-Wiener process can be represented via its Karhunen--Loève expansion, which allows the arguments to be carried out similarly. 
		
		Therefore, we restrict ourselves to the one-dimensional noise case primarily to avoid cumbersome notation.}
\end{remark}

\subsection{Solution concepts and PDE results}
First, we state the following theorem, {which is a consequence of} \cite[Theorem 1.1]{wu2018global} about the existence and uniqueness of {a variational solution} to \eqref{eq1.1}.

\begin{theorem} \label{thm:well-posed} Let $(\Omega, \mathcal{F}, \mathbb{P}, \{\mathcal{F}_t\}_{t \geq 0} )$ be a probability space. Let $W$ be an ${\mathbb R}$-valued Wiener process. 
	Suppose ${u}_0\in L^2(\Omega; L^2(D))$ {and that $B$ satisfies condition \eqref{Assump_Lipschitz}}. Then, there exists a unique weak solution $u \in L^2(\Omega; L^2[0,T]; L^2(D))$ that satisfies $\mP$-a.s.
	\begin{align}\label{Weak_formulation}
		&(u(t), \phi) + \int_{0}^{t}\nu\bigl(\p_x^2 u(s) , \p_x^2\phi\bigr)\, ds - \int_0^t \bigl(\p_x u(s), \p_x\phi\bigr)\, ds + \int_{0}^t \bigl(u(s)\p_xu(s), \phi\bigr)\, ds
		\\\nonumber&\qquad= \bigl(u_0,\phi\bigr) + \left(\int_0^t (B(u)\, dW(s), \phi\right)\qquad \forall \phi \in H^2(D).
	\end{align}

\end{theorem}

{
	\begin{remark}
		1. We note that although \cite[Theorem~1.1]{wu2018global} establishes only the well-posedness of a mild solution $u$ to \eqref{eq1.1}, the argument in \cite[Theorem~6.5]{da2014stochastic}, together with the assumption $B(v)\in L^2(D)$ for any $v\in L^2(D)$, ensures that $u$ is also a weak solution to \eqref{eq1.1} and satisfies \eqref{Weak_formulation}.
		
		2. We further remark that, although the global solutions in \cite[Theorem~1.1]{wu2018global} are proved for the case $\nu = 1$, the argument therein can be readily adapted to establish Theorem~\ref{thm:well-posed} for any $\nu > 0$.
	\end{remark}
}

Having obtained the well-posedness of equation \eqref{eq1.1}, we turn to the analysis of higher regularities that will be exploited later in Section \ref{sec3} to establish the error estimates of the finite element method. We start with Lemma \ref{lemma_expo_moment_H} below, giving an exponential stability estimate for $u$ assuming that $B(u)$ is bounded in $L^2$ norm, i.e., $B(u)$ satisfies \eqref{Assump_Lineargrowth}. Its proof is a slight rework of that of \cite[Lemma 3.2]{gao2025exponential}, adapting to the multiplicative noise setting.

\medskip 

\begin{lemma} \label{lemma_expo_moment_H}
	Suppose that $u_0\in L^2(\Omega;L^2(D))$ and that $B$ satisfies \eqref{Assump_Lineargrowth}. Then, for all {$\kappa\in \left(0,\frac{1}{32L_0^2}\right)$} , there holds
	\begin{align} \label{expo_moment_H}
		\E\left[\exp\left\{ \kappa\|u(t)\|^2_{L^2} + \frac{\kappa}{2}\nu \int_0^t \|\partial_x^2u(s)\|^2_{L^2}d s \right\} \right]& \le {\sqrt{2}\,\E\left[\exp\left\{ 16\kappa\|u_0\|^2_{L^2}+\kappa c_0t+\kappa \tilde{c}_0 \right\}\right]},
	\end{align}
	for all $t \in [0,T]$, and for some positive constants $c_0$ and $\tilde{c}_0$ independent of $\kappa$, $u_0$ and $t$.
\end{lemma}
\begin{proof}
	First of all, we consider $\|u(t)-\varphi_{b(t)}\|^2_{L^2([0,2L])}$, where $\f_{b(t)}$ is the function from \eqref{form:phi_b(t)}. In particular, \eqref{form:phi_b(t)}-\eqref{eqn:b(t)} imply that
	\begin{align*}
		d\,\f_{b(t)} = \partial_x\f(\cdot+b(t)) b'(t)\,dt=\frac{1}{4L}\partial_x\f(\cdot+b(t)) \left( u,\partial_x\f(\cdot+b(t)) \right)_{L^2([0,2L])}\,dt.
	\end{align*}
	From the above identity and \eqref{eq1.1}, using the It\^o's formula gives
	\begin{align} \label{eqn:d|u-Phi_b|^2}
		\frac{1}{2}d\|u-\f_{b}\|^2_{L^2([0,2L])} & = \left( u-\f_b,-\nu \partial_x^4 u-\partial_x^2 u-u\partial_x u  -\partial_x\f_b b'\right)_{L^2([0,2L])}d t \notag \\
		&\qquad + \left( u-\f_b,B(u) d W\right)_{L^2([0,2L])}+\frac{1}{2}\|B(u)\|^2_{L^2([0,2L])}d t.
	\end{align}
	Since $u$ is $L$-periodic and $\f_b$ is $2L$-periodic, it holds that
	\begin{align*}
		\left(u^2,\partial_x u  \right)_{L^2([0,2L])}=0,\quad \left(\f_b,\partial_x \f_b  \right)_{L^2([0,2L])}=0.
	\end{align*}
	We use integration by parts to compute
	\begin{align*}
		& \left( u-\f_b,-\nu \partial_x^4 u-\partial_x^2u-u\partial_x u  -\partial_x\f_b b'\right)_{L^2([0,2L])}\\
		&= -\nu\|\partial_x^2u\|^2_{L^2([0,2L])}+\|\partial_x u \|^2_{L^2([0,2L])} -\left( u,\partial_x\f_b\right)_{L^2([0,2L])}b' \\
		&\qquad +\nu\left( \f_b,\partial_x^4  u\right)_{L^2([0,2L])}+\left( \f_b,\partial_x^2u\right)_{L^2([0,2L])}+\left( u\partial_x u ,\f_b\right)_{L^2([0,2L])}\\
		&=  \|\partial_x u \|^2_{L^2([0,2L])} -\nu\|\partial_x^2 u\|^2_{L^2([0,2L])}  -\frac{1}{2}\left( u^2,\partial_x\f_b\right)_{L^2([0,2L])}\\
		&\qquad -\left( u,\partial_x\f_b\right)_{L^2([0,2L])}b'+\left( \f_b,\partial_x^2 u\right)_{L^2([0,2L])}+\nu\left( \partial_x^2 \f_b,\partial_x^2  u\right)_{L^2([0,2L])}.
	\end{align*}
	With regard to the last two terms on the above right-hand side, we employ the Cauchy-Schwarz inequality to infer
	\begin{align*}
		\left( \f_b,\partial_x^2 u\right)_{L^2([0,2L])}+\nu\left( \partial_x^2 \f_b,\partial_x^2  u\right)_{L^2([0,2L])} 
		&\le \frac{1}{\nu}\|\f_b\|^2_{L^2([0,2L])}+\nu\|\partial_x^2 \f_b\|^2_{L^2([0,2L])}+\frac{\nu}{2}\|\partial_x^2 u\|^2_{L^2([0,2L])},
	\end{align*}
	whence
	\begin{align*}
		& \left( u-\f_b,-\nu  \partial_x^4 u-\partial_x^2 u-u\partial_x u  -\partial_x\f_b b'\right)_{L^2([0,2L])}\\
		&\le \|\partial_x u \|^2_{L^2([0,2L])} -\frac{1}{2}\nu \|\partial_x^2 u\|^2_{L^2([0,2L])}  -\frac{1}{2}\left( u^2,\partial_x\f_b\right)_{L^2([0,2L])}\\
		&\qquad -\left( u,\partial_x\f_b\right)_{L^2([0,2L])}b'+\frac{1}{\nu }\|\f_b\|^2_{L^2([0,2L])}+\nu \|\partial_x^2 \f_b\|^2_{L^2([0,2L])}\\
		&\le - \frac{1}{8}\nu \|\partial_x^2u\|^2_{L^2([0,2L])} -\frac{1}{2}\|u\|^2_{L^2([0,2L])}+\frac{1}{4L}\big|\left( u,\partial_x\f_b\right)_{L^2([0,2L])}  \big|^2\\
		&\qquad-\left( u,\partial_x\f_b\right)_{L^2([0,2L])}b'+\frac{1}{\nu }\|\f_b\|^2_{L^2([0,2L])}+\nu \|\partial_x^2\f_b\|^2_{L^2([0,2L])}.
	\end{align*}
	where the last implication follows from \eqref{ineq:phi}. At this point, we recall that $b(t)$ satisfies \eqref{eqn:b(t)}, implying
	\begin{align*}
		& \left( u-\f_b,-\nu  \partial_x^4 u-\partial_x^2 u-u\partial_x u  -\partial_x\f_b b'\right)_{L^2([0,2L])}\\
		&\le - \frac{1}{8}\nu \|\partial_x^2u\|^2_{L^2([0,2L])} -\frac{1}{2}\|u\|^2_{L^2([0,2L])}+\frac{1}{\nu }\|\f_b\|^2_{L^2([0,2L])}+\nu \|\partial_x^2\f_b\|^2_{L^2([0,2L])}.
	\end{align*}
	Also, we employ the elementary inequality $2(a^2+b^2)\ge (a+b)^2$, $a,b\in\mathbb{R}$, to see that
	\begin{align*}
		-\frac{1}{2}\|u\|^2_{L^2([0,2L])}\le -\frac{1}{4}\|u-\f_b\|^2_{L^2([0,2L])}+\frac{1}{2}\|\f_b\|^2_{L^2([0,2L])}.
	\end{align*}
	As a consequence, we obtain
	\begin{align} \label{ineq:<u-phi>}
		&\left( u-\f_b,-\partial_x^2u-\nu  \partial_x^4 u-u\partial_x u  -\partial_x\f_b b'\right)_{L^2([0,2L])} \notag \\
		&\le - \frac{1}{8}\nu \|\partial_x^2u\|^2_{L^2([0,2L])} -\frac{1}{4}\|u-\f_b\|^2_{L^2([0,2L])} \\\nonumber
		&\qquad+\left(\frac{1}{2}+\frac{1}{\nu }\right)\|\f_b\|^2_{L^2([0,2L])}+\nu \|\partial_x^2\f_b\|^2_{L^2([0,2L])}.
	\end{align}
	From \eqref{eqn:d|u-Phi_b|^2}, {we recall condition \eqref{Assump_Lineargrowth} and $L$-periodicity to see that
		\begin{align*}
			\|B(u)\|^2_{L^2([0,2L])} = 2 \|B(u)\|^2_{L^2([0,L])} \le 2L_0^2.
	\end{align*}}
	This together with \eqref{ineq:<u-phi>} implies
	\begin{align} \label{ineq:d|u-Phi_b|^2}
		d\|u-\f_{b}\|^2_{L^2([0,2L])} & \le - \frac{1}{4}\nu \|\partial_x^2u\|^2_{L^2([0,2L])}d t -\frac{1}{2}\|u-\f_b\|^2_{L^2([0,2L])}d t+d M+c_0dt,
	\end{align}
	where we have set
	\begin{align} \label{form:M(t)}
		M(t) := 2\int_0^t\left( u(s)-\f_{b(s)},B(u(s))d W(s)\right)_{L^2([0,2L])},
	\end{align}
	and {
		\begin{align*}
			c_0 := \Big(1+\frac{2}{\nu }\Big)\|\f_b\|^2_{L^2([0,2L])}+2\nu \|\partial_x^2\f_b\|^2_{L^2([0,2L])} +2L_0^2.
	\end{align*}}
	We integrate both sides of \eqref{ineq:d|u-Phi_b|^2} with respect to time $t$ and obtain
	\begin{align*}
		&\|u(t)-\f_{b(t)}\|^2_{L^2([0,2L])} + \frac{1}{4}\nu \int_0^t \|\partial_x^2u(s)\|^2_{L^2([0,2L])}d s +\frac{1}{2}\int_0^t\|u(s)-\f_{b(s)}\|^2_{L^2([0,2L])}d s\\
		& \le\|u_0-\f\|^2_{L^2([0,2L])}+M(t)+c_0t.
	\end{align*}
	It follows that
	\begin{align*}
		&\frac{1}{2}\|u(t)\|^2_{L^2([0,2L])} + \frac{1}{4}\nu \int_0^t \|\partial_x^2u(s)\|^2_{L^2([0,2L])}d s +\frac{1}{2}\int_0^t\|u(s)-\f_{b(s)}\|^2_{L^2([0,2L])}d s\\
		& \le \|\f_{b(t)}\|^2_{L^2([0,2L])}+2\|u_0\|^2_{L^2([0,2L])}+2\|\f\|^2_{L^2([0,2L])} +M(t)+c_0t.
	\end{align*}
	Since $u$ is $L$-periodic and $\f$ is $2L$-periodic, we have
	\begin{align} \label{eqn:L-periodic}
		\|u\|^2_{L^2([0,2L])} =2\|u\|^2_{L^2},\quad \|\f_{b(t)}\|^2_{L^2([0,2L])}= \|\f(\cdot +b(t))\|^2_{L^2([0,2L])} = \|\f\|^2_{L^2([0,2L])}.
	\end{align}
	So, for all $\kappa>0$, it holds that
	\begin{align} \label{ineq:|u|^+int.|u|^2<M}
		& \kappa\|u(t)\|^2_{L^2} + \frac{\kappa}{2}\nu \int_0^t \|\partial_x^2u(s)\|^2_{L^2}d s +\kappa\int_0^t\|u(s)-\f_{b(s)}\|^2_{L^2([0,2L])}d s 
		\notag\\
		& \le 8\kappa\|u_0\|^2_{L^2}+2\kappa M(t)+2\kappa c_0t+6\kappa\|\f\|^2_{L^2([0,2L])}.
	\end{align}
	
	Turning back to \eqref{expo_moment_H}, we aim to employ the exponential Martingale inequality to establish the exponential moment bound on $u(t)$. To see this, we recall expression \eqref{form:M(t)} of $M(t)$ and note that the corresponding quadratic variation process $\langle M\rangle(t)$ satisfies
	\begin{align*}
		d\langle M\rangle = 4 \big|  \left( u-\f_b,B(u)\right)_{L^2([0,2L])}\big|^2d t&\le 4 \|B(u)\|^2_{L^2([0,2L])} \|u-\f_b\|^2_{L^2([0,2L])} d t\\
		&\le 8L_0^2\|u-\f_b\|^2_{L^2([0,2L])} d t.
	\end{align*}
	In the last implication above, we invoked periodicity and the extra hypothesis that $\|B(u)\|_{L^2}\le L_0$. In view of \eqref{ineq:|u|^+int.|u|^2<M}, provided that 
	\begin{align*}
		0<\kappa < \frac{1}{32L_0^2},
	\end{align*}
	we deduce that
	\begin{align} \label{ineq:|u|^+int.|u|^2<M-<M>}
		& \kappa\|u(t)\|^2_{L^2} + \frac{\kappa}{2}\nu \int_0^t \|\partial_x^2u(s)\|^2_{L^2}d s -4\kappa\|u_0\|^2_{L^2}-2\kappa c_0t-6\kappa\|\f\|^2_{L^2([0,2L])} \notag \\
		&\le 2\kappa M(t)-\kappa\int_0^t\|u(s)-\f_{b(s)}\|^2_{L^2([0,2L])}d s \notag \\
		& \le 2 \kappa M(t) -4\kappa^2\langle M\rangle(t).
	\end{align}
	Next, we recall the exponential martingale inequality
	\begin{align*}
		\P\Big( \sup_{t\ge 0}\Big[2\kappa M(t)-\frac{1}{2}\lambda\cdot 4\kappa^2 \langle M\rangle(t)\Big]\ge R\Big) \le e^{-\lambda R},\quad \lambda>0,\, R>0.
	\end{align*}
	Based on the above right-hand side, we pick $\lambda=2$ and obtain
	\begin{align*}
		\P\Big( \sup_{t\ge 0}\Big[2\kappa M(t)-4\kappa^2\langle M\rangle(t)\Big]\ge R\Big) \le e^{-2 R},
	\end{align*}
	implying
	\begin{align*} 
		\E\left[\exp\Big\{ \sup_{t\ge 0}\Big[2\kappa M(t)-4\kappa^2\langle M\rangle(t)\Big] \Big\}\right] \le 2.
	\end{align*}
	From \eqref{ineq:|u|^+int.|u|^2<M-<M>}, we get immediately that
	\begin{align*}
		& \E\left[\exp\left\{ \sup_{t\ge 0} \left(\kappa\|u(t)\|^2_{L^2} + \frac{\kappa}{2}\nu \int_0^t \|\partial_x^2u(s)\|^2_{L^2}\, d s -8\kappa\|u_0\|^2_{L^2}-2\kappa c_0t-6\kappa\|\f\|^2_{L^2([0,2L])}\right)\right\}\right]\\
		&\le \E\left[\exp\Big\{ \sup_{t\ge 0}\Big[\kappa M(t)-4\kappa^2\langle M\rangle(t)\Big] \Big\}\right] \le 2,
	\end{align*}
	whence
	\begin{align*}
		& \E\left[\exp\left\{ \sup_{t\in[0,T]} \left[\kappa\|u(t)\|^2_{L^2} + \frac{\kappa}{2}\nu \int_0^t \|\partial_x^2u(s)\|^2_{L^2}\, d s -8\kappa\|u_0\|^2_{L^2}\right]\right\}\right]\\
		&\le 2\exp\left\{2 \kappa c_0T+6\kappa\|\f\|^2_{L^2([0,2L])}\right\}.
	\end{align*}
	Lastly, we invoke H\"older's inequality to deduce
	\begin{align*}
		& \E\left[\exp\left\{ \sup_{t\in [0,T]} \left[\kappa\|u(t)\|^2_{L^2} + \frac{\kappa}{2}\nu \int_0^t \|\partial_x^2u(s)\|^2_{L^2}\, d s\right]\right\}\right]\\
		&\le \left(\E\left[\exp\left\{ \sup_{t\in [0,T]} \left[2\kappa\|u(t)\|^2_{L^2} + \kappa\nu \int_0^t \|\partial_x^2u(s)\|^2_{L^2}d s -16\kappa\|u_0\|^2_{L^2}\right]\right\}\right]\right)^{\frac12}\\
		&\qquad\qquad
		\times\left(\E\left[\exp\left\{ 16\kappa\|u_0\|^2_{L^2} \right\}\right]\right)^{\frac12}.
	\end{align*}
	It follows that
	\begin{align*}
		& \E\left[ \exp\left\{ \sup_{t\in [0,T]} \left[\kappa\|u(t)\|^2_{L^2} + \frac{\kappa}{2}\nu \int_0^t \|\partial_x^2u(s)\|^2_{L^2}d s\right]\right\}\right]\\
		&\le\sqrt{2}\exp\left\{ 2\kappa c_0T+6\kappa\|\f\|^2_{L^2([0,2L])}\right\} \left(\E\left[\exp\left\{ 16\kappa\|u_0\|^2_{L^2} \right\}\right]\right)^{\frac12}.
	\end{align*}
	In turn, this establishes \eqref{expo_moment_H}, thereby finishing the proof.
	
\end{proof}

Having established an exponential stability estimate, we turn to higher regularity for $u$ provided the solution starts from sufficiently smooth initial conditions. More precisely, we have the following result, whose proof is drawn upon the technique in \cite{flandoli1995martingale} dealing with the same energy estimates for the Navier-Stokes equation.

\begin{lemma} \label{lemma_poly_moment_H^s}
	Let $m,q\in \mathbb{N}$ with $q\ge 2$, and $B$ satisfy \eqref{Assump_Lipschitz}. {Suppose that $u_0$ satisfies 
		\begin{align} \label{cond:moment_H^s:int_0^T}
			u_0\in \bigcap_{i=0}^m L^{2^{m-i}q}(\Omega;H^{i}(D)).
		\end{align}
		Then, the following hold
		\begin{align} \label{ineq:poly_moment_H^s:int_0^T}
			\sup_{s\in[0,T]}\E[\|u(s)\|^q_{H^m}]+ \int_0^T \mE\left[\|u(s)\|^{q-2}_{H^{m}}\|u(s)\|^2_{H^{m+2}}\right]ds \leq C_{m,q},
		\end{align}
		and
		\begin{align} \label{ineq:poly_moment_H^s:sup_[0,T]}
			\mE\left[\sup_{s\in [0,T]} \|u(s)\|^q_{H^{m}}\right]   \leq C_{m,q},
		\end{align}
		where $C_{m,q} = C_T\left(\sum_{i=0}^m \E\left[\|u_0\|^{2^{m-i}q}_{H^i}\right]+1\right)$.}

\end{lemma}

\begin{proof}
	We proceed to establish \eqref{ineq:poly_moment_H^s:int_0^T} and \eqref{ineq:poly_moment_H^s:sup_[0,T]} by induction on the regularity parameter $m$. We start with the base case $m=0$ and apply It\^o's formula to $\|u\|^q_{L^2}$ and obtain the identity
	\begin{align*}
		&d \|u\|^q_{L^2} -q\|u\|^{q-2}_{L^2}\|\partial_x u\|^2_{L^2} dt+\nu q\|u\|^{q-2}_{L^2}\|\partial_x^2 u\|^2_{L^2}\,dt {+ \frac{1}{2}q\|u\|^{q-2}_{L^2}\|B(u)\|^2_{L^2}dt}\\
		&\quad= q\|u\|^{q-2}_{L^2}\left(u,B(u)d W\right)+\frac{1}{2}q(q-2)\|u\|^{q-4}_{L^2}|(u,B(u))|^2dt.
	\end{align*}
	On the one hand, we employ integration by parts to estimate
	\begin{align*}
		-q\|u\|^{q-2}_{L^2}\|\partial_x u\|^2_{L^2}  = -q\|u\|^{q-2}_{L^2}(u,\partial_x^2 u) \le \frac{1}{100}\nu q \|u\|^{q-2}_{L^2}\|\partial_x^2 u\|^{2}_{L^2}+ C\|u\|^q_{L^2},
	\end{align*}
	for a positive constant $C=C(\nu,q)$. On the other hand, we invoke the estimate \eqref{cond:|B(u)|_Hm<|u|_Hm} to infer
	\begin{align*}
		\frac{1}{2}q\|u\|^{q-2}_{L^2}\|B(u)\|^2_{L^2}+\frac{1}{2}q(q-2)\|u\|^{q-4}_{L^2}|(u,B(u))|^2 \le C\big(1+\|u\|^q_{L^2}\big).
	\end{align*}
	Altogether, we get
	\begin{align} \label{ineq:d|u|^q_L2}
		&d \|u\|^q_{L^2}+ \frac{1}{2}\nu q\|u\|^{q-2}_{L^2}\|\partial_x^2 u\|^2_{L^2}\, dt\le  q\|u\|^{q-2}_{L^2}\left(u,B(u)d W\right)+ C\|u\|^{q}_{L^2}dt+Cdt.
	\end{align}
	We integrate the above estimate with respect to time and take expectations on both sides to deduce
	\begin{align*}
		&\E\left[\|u(t)\|^q_{L^2}\right] + \frac{1}{2}\nu q\int_0^t\E\left[\|u(s)\|^{q-2}_{L^2}\|\partial_x^2 u(s)\|^2_{L^2}\right]ds\\
		&\le \E\left[\|u_0\|_{L^2}^q\right] + C\int_0^t\E\left[\|u(s)\|^{q}_{L^2}\right]\, ds+Ct.
	\end{align*}
	It follows immediately from Gronwall's inequality that
	\begin{align*}
		\E\left[\|u(t)\|^q_{L^2}\right] \le Ce^{Ct}\E\left[\|u_0\|^q_{L^2}\right],\quad t\ge 0.
	\end{align*}
	In turn, 
	\begin{align} \label{ineq:int_0^T|u|^(q-2)_L2.|u|^2_H2}
		\frac{1}{2}\nu q\int_0^T\mE\left[\|u(s)\|^{q-2}_{L^2}\|\partial_x^2 u(s)\|^2_{L^2}\right]ds &\le \E\left[\|u_0\|_{L^2}^q\right] + C \int_0^T \mE\left[\|u(s)\|^{q}_{L^2}\right]\, ds+CT\notag \\
		&\le Ce^{CT}\big(\mE\left[\|u_0\|^q_{L^2}\right]+1\big).
	\end{align}
	In particular, this produces \eqref{ineq:poly_moment_H^s:int_0^T} for the base case $m=0$ and for all $q\ge 2$.
	
	Now, from \eqref{ineq:d|u|^q_L2}, we observe that
	\begin{align*}
		\E\left[\sup_{s\in[0,T]}\|u(s)\|^q_{L^2} \right]\le \E\left[\|u_0\|^q_{L^2}\right]+\E\left[\sup_{s\in[0,T]}|M_0(s)| \right]+C\int_0^T\E[\|u(s)\|^q_{L^2}]ds+CT,
	\end{align*}
	where $M_0$ is the Martingale process defined as
	\begin{align*}
		dM_0= q\|u\|^{q-2}_{L^2}\left(u,B(u)d W\right),
	\end{align*}
	and whose associated quadratic variation process is given by
	{\begin{align*}
			\langle M_0(t)\rangle= \int_0^t q^2\|u(s)\|^{2q-4}_{L^2}|\left(u(s),B(u(s))\right)|^2ds.
		\end{align*}
		Furthermore, {estimate \eqref{cond:|B(u)|_Hm<|u|_Hm}} together with Holder's inequality implies that
		\begin{align*}
			\langle M_0(t)\rangle \le C\int_0^t \|u(s)\|_{L^2}^{2q-2}(1+\|u(s)\|^{2}_{L^2})ds.
	\end{align*}}
	{	So, we invoke the Burkholder-Gundy-Davis inequality to infer
		\begin{align*}
			\E\left[\sup_{s\in[0,T]}|M_0(s)| \right] &\le C\E\left[   \left(\int_0^T\|u(s)\|_{L^2}^{2q-2}(1+\|u(s)\|^{2}_{L^2})ds\right)^{\frac{1}{2}}\right]\\
			& \le C \E\left[\sup_{s\in [0,T]}\|u(s)\|_{L^2}^{q-1} \left(\int_0^T(1+\|u(s)\|^{2}_{L^2})ds\right)^{\frac{1}{2}} \right].
		\end{align*}
		An application of Young's inequality produces (recalling $q\ge 2$)
		\begin{align*}
			& C \E\left[\sup_{s\in [0,T]}\|u(s)\|_{L^2}^{q-1} \left(\int_0^T(1+\|u(s)\|^{2}_{L^2})ds\right)^{\frac{1}{2}} \right] \\
			&\le \frac{1}{100} \E\left[\sup_{s\in [0,T]}\|u(s)\|_{L^2}^{q}  \right] + C \E\left[\left(\int_0^T(1+\|u(s)\|^{2}_{L^2})ds\right)^{\frac{q}{2}} \right]\\
			&\le \frac{1}{100} \E\left[\sup_{s\in [0,T]}\|u(s)\|_{L^2}^{q}  \right] + C \E\left[\int_0^T(1+\|u(s)\|^{q}_{L^2})ds \right].
		\end{align*}
		As a consequence, we obtain
		\begin{align*}
			\E\left[\sup_{s\in[0,T]}|M_0(s)| \right]&\le \frac{1}{100} \E\left[\sup_{s\in [0,T]}\|u(s)\|_{L^2}^{q}  \right] + CT+C\int_0^T\E[ \|u(s)\|^{q}_{L^2}]ds.
		\end{align*}
		It follows that
		\begin{align*}
			\frac{99}{100}\E\left[\sup_{s\in[0,T]}\|u(s)\|^q_{L^2} \right]& \le \E\left[\|u_0\|^q_{L^2}\right]+CT ++C\int_0^T\E[\|u(s)\|^q_{L^2}]ds\\
			&\le Ce^{CT} \big( \E[\|u_0\|^{q}_{L^2}]+1\big) .
	\end{align*} }
	This establishes \eqref{ineq:poly_moment_H^s:sup_[0,T]} and concludes the proof for the base case $m=0$ and for all $q\ge 2$.
	
	Next, suppose that \eqref{ineq:poly_moment_H^s:int_0^T} and \eqref{ineq:poly_moment_H^s:sup_[0,T]} hold for all $0,\dots,m-1$ (where $m\ge 1$) and for all $q\ge 2$. Let us consider the general case $m$. Similar to the base case, we apply It\^o's formula to $\|\partial_x^{m}u\|^q$ and obtain the identity
	\begin{align} \label{eqn:d|u|^q_H^n}
		&d \|\partial_x^m u\|^q_{L^2} -q\|\partial_x^m u\|^{q-2}_{L^2}\|\partial_x^{m+1} u\|^2_{L^2} dt+\nu q\|\partial_x^m u\|^{q-2}_{L^2}\|\partial_x^{m+2} u\|^2_{L^2}dt\notag \\
		&\qquad+ q\|u\|^{q-2}_{L^2}\left(\partial_x^{m-2}[u\partial_x u],\partial_x^{m+2}u\right)\,dt   \notag \\
		&\quad= q\|\partial_x^m u\|^{q-2}_{L^2}\left(\partial_x^m u, \partial_x^m B(u)d W\right)+\frac{1}{2}q(q-2)\|\partial_x^m u\|^{q-4}_{L^2}|(\partial_x^m u,\partial_x^m B(u))|^2dt\notag \\
		&\hspace{3cm} {+\frac{1}{2}q \|\partial_x^m u\|^{q-2}_{L^2}\|\partial_x^mB(u)\|^2_{L^2}dt  }.
	\end{align}
	Once again, we employ integration by parts to infer
	\begin{align*}
		q\|\partial_x^m u\|^{q-2}_{L^2}\|\partial_x^{m+1} u\|^2_{L^2} &= -q\|\partial_x^m u\|^{q-2}_{L^2}(\partial_x^{m} u, \partial_x^{m+2}u)\\
		&\le  \frac{1}{100}\nu q \|\partial_x^m u\|^{q-2}_{L^2}\|\partial_x^{m +2}u\|^2_{L^2} +C\|\partial_x^m u\|^2_{L^2}.
	\end{align*}
	Also, we invoke {estimate \eqref{cond:|B(u)|_Hm<|u|_Hm}} with Holder's inequality to estimate
	\begin{align*}
		\frac{1}{2}q \|\partial_x^m u\|^{q-2}_{L^2}\|\partial_x^mB(u)\|^2_{L^2} & +\frac{1}{2}q(q-2)\|\partial_x^m u\|^{q-4}_{L^2}|(\partial_x^m u,\partial_x^m B(u))|^2\\
		&\le C\big(1+ \|\partial_x^m u\|^{q}_{L^2}\big).
	\end{align*}
	With regard to the last term on the left-hand side of \eqref{eqn:d|u|^q_H^n}, observe that
	\begin{align*}
		\partial_x^{m-2}[u\partial_x u] = \sum_{i=0}^{m-2}a_i \partial_x^i u \partial_x^{m-1-i}u.
	\end{align*}
	We employ the embedding $H^1\subset L^\infty$  (in dimension one), cf. \eqref{sobolev_ineq}, to estimate for each $i=0,\dots,m-2$
	\begin{align*}
		\left(\partial_x^i u \partial_x^{m-1-i}u,\partial_x^{m+2}u\right)\le \|\partial_x^i u\|_{L^\infty}\|\partial_x^{m-1-i}u\|_{L^2}\|\partial_x^{m+2}u\|_{L^2}&\le C_e\|\partial_x^{m-1}u\|^2_{L^2}\|\partial_x^{m+2}u\|_{L^2}.
	\end{align*}
	{In the above, $C_e$ is the constant in \eqref{sobolev_ineq}.} It follows that
	\begin{align*}
		&q\|\partial_x^m u\|^{q-2}_{L^2}\left|\left(\partial_x^{m-2}[u\partial_x u],\partial_x^{m+2}u\right)\right| \\
		&\quad  \le Cq\|\partial_x^m u\|^{q-2}_{L^2}\|\partial_x^{m-1}u\|^2_{L^2}\|\partial_x^{m+2}u\|_{L^2}\\
		&\quad \le \frac{1}{100}\nu q\|\partial_x^m u\|^{q-2}_{L^2}\|\partial_x^{m+2}u\|_{L^2}^2 + C\|\partial_x^m u\|^{q}_{L^2} + C \|\partial_x^{m-1}u\|^{2q}_{L^2}.
	\end{align*}
	We collect all of the above estimates with \eqref{eqn:d|u|^q_H^n} to deduce the bound
	\begin{align} \label{ineq:|D^nu|^q}
		&d \|\partial_x^m u\|^q_{L^2} +\frac{1}{2}\nu q\|\partial_x^m u\|^{q-2}_{L^2}\|\partial_x^{m+2} u\|^2_{L^2}\,dt \notag \\
		&\quad\le q\|\partial_x^m u\|^{q-2}_{L^2}\left(\partial_x^m u, \partial_x^m B(u)d W\right)+C\|\partial_x^m u\|^{q}_{L^2}dt + C \|\partial_x^{m-1}u\|^{2q}_{L^2}dt+Cdt,
	\end{align}
	implying
	\begin{align*}
		&\E\left[\|\partial_x^m u(t)\|^q_{L^2}\right] +\frac{1}{2}\nu q\int_0^t\mE\left[\|\partial_x^m u(s)\|^{q-2}_{L^2}\|\partial_x^{m+2} u(s)\|^2_{L^2}\right]\,ds \\
		&\quad\le \mE\left[\|\partial_x^m u_0\|^q_{L^2}\right]+C\int_0^t\E\left[\|\partial_x^m u(s)\|^{q}_{L^2}\right]\, ds + C\int_0^t \E\left[\|\partial_x^{m-1}u(s)\|^{2q}_{L^2}\right]\,ds+Ct.
	\end{align*}
	In light of the induction hypothesis, we readily have
	\begin{align*}
		\int_0^t \E\left[\|\partial_x^{m-1}u(s)\|^{2q}_{L^2}\right]\,ds \le C\Big(\sum_{i=0}^{m-1}\mE\left[\|\partial_{x}^{i} u_0\|^{2^{m-i}q}_{L^2}\right]+1\Big),
	\end{align*}
	implying
	\begin{align*}
		&\E\left[\|\partial_x^m u(t)\|^q_{L^2}\right] +\frac{1}{2}\nu q\int_0^t\mE\left[\|\partial_x^m u(s)\|^{q-2}_{L^2}\|\partial_x^{m+2} u(s)\|^2_{L^2}\right]ds \\
		&\quad\le \mE\left[\|\partial_x^m u_0\|^q_{L^2}\right] + C\int_0^t\mE\left[\|\partial_x^m u(s)\|^{q}_{L^2}\right]\,ds + C\sum_{i=0}^{n-1}\mE\left[\|\partial_{x}^{i} u_0\|^{2^{m-i}q}_{L^2}\right]+C\\
		&\quad \le C\int_0^t\mE\left[\|\partial_x^m u(s)\|^{q}_{L^2}\right]\,ds + C\sum_{i=0}^{m}\mE\left[\|\partial_{x}^{i} u_0\|^{2^{m-i}q}_{L^2}\right]+C.
	\end{align*}
	As a consequence of Gronwall's inequality, we obtain
	\begin{align*}
		\mE\left[\|\partial_x^m u(t)\|^q_{L^2}\right] \le Ce^{Ct}\Big(\sum_{i=0}^{m}\mE\left[\|\partial_{x}^{i} u_0\|^{2^{m-i}q}_{L^2}\right]+1\Big).
	\end{align*}
	In turn, it holds that
	\begin{align*}
		&\frac{1}{2}\nu q\int_0^t\mE\left[\|\partial_x^m u(s)\|^{q-2}_{L^2}\|\partial_x^{m+2} u(s)\|^2_{L^2}\right]\,ds \\
		&\quad\le C\int_0^t\mE\left[\|\partial_x^m u(s)\|^{q}_{L^2}\right]\,ds + C\sum_{i=0}^{m}\mE\left[\|\partial_{x}^{i} u_0\|^{2^{m-i}q}_{L^2}\right] +C\\
		&\quad \le C\sum_{i=0}^{m}\mE\left[\|\partial_{x}^{i} u_0\|^{2^{m-i}q}_{L^2}\right]+C.
	\end{align*}
	This produces \eqref{ineq:poly_moment_H^s:int_0^T} for the general case $m\ge 1$ and for all $q\ge 2$.
	
	Turning to \eqref{ineq:poly_moment_H^s:sup_[0,T]}, from \eqref{ineq:|D^nu|^q}, we have
	\begin{align*}
		&\E\left[\sup_{s\in[0,T]}\|\partial_x^m u(s)\|^q_{L^2}\right] \\
		&\le \E[\|\partial_x^m u_0\|^q_{L^2}]+ \E \left[\sup_{s\in[0,T]}|M_m(s)|\right] +C\int_0^T\E[\|\partial_x^m u(s)\|^{q}_{L^2}]+ \E[\|\partial_x^{m-1}u(s)\|^{2q}_{L^2}]ds+C,
	\end{align*}
	where $M_m$ is the Martingale defined as
	\begin{align*}
		dM_m= q\|\partial_x^m u\|^{q-2}_{L^2}\left(\partial_x^m u, \partial_x^m B(u)d W\right),
	\end{align*}
	and whose quadratic variation process satisfies the bound (recalling inequality \eqref{cond:|B(u)|_Hm<|u|_Hm})
	\begin{align*}
		d\langle M_m\rangle\le C\|\partial_x^m u\|^{2q-2}_{L^2}(1+\|\partial_x^m u\|^{2}_{L^2})dt.
	\end{align*}
	Similar to the estimate for $M_0$ in the base case, we also obtain
	\begin{align*}
		\E \left[\sup_{s\in[0,T]}|M_m(s)|\right] \le \frac{1}{100} \E\left[\sup_{s\in [0,T]}\|\partial_x^m u(s)\|_{L^2}^{q}  \right] + CT+C\int_0^T\E[ \|\partial_x^m u(s)\|^{q}_{L^2}]ds,
	\end{align*}   
	implying
	\begin{align*}
		\frac{99}{100}\E\left[\sup_{s\in[0,T]}\|\partial_x^m u(s)\|^q_{L^2}\right] \le \E[\|\partial_x^m u_0\|^q_{L^2}]+CT+C\int_0^T \E[\|\partial_x^{m}u(s)\|^{q}_{L^2}]ds.
	\end{align*}
	In view of \eqref{ineq:poly_moment_H^s:int_0^T}, it holds that
	\begin{align*}
		\int_0^T \E[\|\partial_x^{m}u(s)\|^{q}_{L^2}]+\E[\|\partial_x^{m-1}u(s)\|^{2q}_{L^2}]ds\le C\sum_{i=0}^{m}\mE\left[\|\partial_{x}^{i} u_0\|^{2^{m-i}q}_{L^2}\right]+C.
	\end{align*}
	As a consequence, we get
	\begin{align*}
		&\E\left[\sup_{s\in[0,T]}\|\partial_x^m u(s)\|^q_{L^2}\right]\le C\left(\sum_{i=0}^{m}\mE\left[\|\partial_{x}^{i} u_0\|^{2^{m-i}q}_{L^2}\right]+1\right).
	\end{align*}
	This establishes \eqref{ineq:poly_moment_H^s:sup_[0,T]} for the general case $m\ge 1$ and for all $q\ge 2$. The proof is thus complete.
\end{proof}

As a consequence of Lemma \ref{lemma_poly_moment_H^s}, we state and prove the following H\"older continuity estimates for the solution of \eqref{Weak_formulation}.

\begin{lemma} \label{lemma_Holder}
	Given integers $m\ge 0$ and $q\ge 1$, suppose that $B$ satisfies \eqref{Assump_Lipschitz} and that
	
	{   1. when $q=1$, 
		\begin{align*}
			u_0\in \bigcap_{i=0}^{m+2} L^{2^{m+3-i}}(\Omega;H^i(D)); 
	\end{align*}  }
	
	{ 2. when $q\ge2$
		\begin{align*}
			u_0\in \bigcap_{i=0}^{(m+1)q} L^{2^{(m+1)q-i}2q}(\Omega;H^i(D)).
	\end{align*}}
	Then, for all $0\le s\le t\le T$, the following holds
	\begin{align} \label{ineq:|u(t)-u(s)|^(2q)_Hn}
		\E\left[\|\partial_x^m (u(t) -u(s)) \|_{L^2}^{2q}\right] \le K_{m,q} \,(t-s)^{q}.
	\end{align}
	for a positive constant $K_{m,q}=K(m,q,T,u_0)$.
\end{lemma}

\bigskip

Since the linear operator $\nu\partial_x^4+\partial_x^2$ does not have a sign-definite, we are not able to follow the semigroup technique in \cite{breit2021convergence,carelli2012rates} to establish Lemma \ref{lemma_Holder}. To circumvent the challenge, we will draw upon the approach of \cite{ferrario2008invariant} dealing with the same issue for the well-posedness of \eqref{eq1.1} in the additive noise setting. The argument essentially consists of two main steps as follows:

\textit{Step 1:} Fixing $a=a(\nu)>0$ such that 
\begin{align}\label{cond:a}
	{\nu \Big(\frac{2\pi}{L}\Big)^4 \ell^4-\Big(\frac{2\pi}{L}\Big)^2 \ell^2+a\ge \frac{1}{2}\nu \Big(\frac{2\pi}{L}\Big)^4\ell^4,\quad \ell\in\Zbb\setminus\{0\},}
\end{align}
let $z(t)$ be the process solving
\begin{align} \label{eqn:z}
	d z+ \nu \partial_x^4 zdt +\partial_x^2 zd t +azdt = B(u)dW,\quad z(0)=0.
\end{align}
It is not difficult to see under the choice of $a$ in \eqref{cond:a}, $\nu \partial_x^4 +\partial_x^2 +a$ becomes strictly positive. We then exploit this fact to prove that $z(t)$ satisfies H\"older estimates provided $B(u)$ satisfies certain regularity properties. This result is captured in Lemma \ref{lemma_Holder_z}, whose detailed proof is presented in Appendix \ref{auxiliary}.

\textit{Step 2:} Setting $v=u-z$, observe that $v$ satisfies the random PDE
\begin{align}
	\label{eqn:v}
	d v+ \nu \partial_x^4 vdt +\partial_x^2 vd t +(v+z)\partial_x(v+z)  -azdt = 0,\quad v(0)=u_0.
\end{align}
Since $v$ is decoupled from the noise term $B(u)dW$, we may employ a pathwise energy estimate approach to establish an analogous regularity property for $v$. This is summarized through Lemma \ref{lemma_Holder_v}. Altogether, we can conclude the proof of Lemma \ref{lemma_Holder} upon recovering $u=v+z$.

We now state the H\"older regularity for the process $v$ through Lemma \ref{lemma_Holder_v} below, whose proof is deferred to Appendix \ref{auxiliary}.

\bigskip 

\begin{lemma} \label{lemma_Holder_v}
	Given integers $m\ge 0$ and $q\ge 1$, suppose that  let $u_0\in \cap_{i=0}^{m+2} L^{2^{m+3-i}q}(\Omega;H^i(D))$ and $B$ satisfies \eqref{Assump_Lipschitz} for $i=0,\dots,m+2$.
	Then, for all $0\le s\le t\le T$, the following holds
	\begin{align} \label{ineq:|v(t)-v(s)|^(2q)_Hn}
		\E\left[\|\partial_x^m (v(t) -v(s)) \|_{L^2}^{2q}\right] \le C(t-s)^q.
	\end{align}
	for a positive constant $C=C(m,q,T,u_0)$.
	
\end{lemma}

\medskip

Assuming the result of Lemma \ref{lemma_Holder_v}, we are in a position to conclude Lemma \ref{lemma_Holder}, giving the H\"older properties for the solution $u(t)$. Since the proof is short, we include it here for the sake of completeness.

\begin{proof}[Proof of Lemma \ref{lemma_Holder}]
	Given integers $m\ge 0$ and $q\ge 2$, we invoke Minkowsi's inequality to deduce
	\begin{align*}
		\|\partial_x^m (u(t)-u(s))\|^q_{L^2} \le C\big(\|\partial_x^m (v(t)-v(s))\|^q_{L^2} +\|\partial_x^m (z(t)-z(s))\|^q_{L^2} \big).
	\end{align*}
	In light of Lemmas \ref{lemma_Holder_v} and \ref{lemma_Holder_z}, we conclude the H\"older property \eqref{ineq:|u(t)-u(s)|^(2q)_Hn}, as claimed.
	
\end{proof}

\section{Fully discrete finite element method} \label{sec3}
\subsection{Formulation and stability of the fully discrete finite element method}\,	{Let $N$ be a positive integer and let
	\[
	0 = x_0 < x_1 < \cdots < x_N = L
	\]
	be a uniform partition of $D=[0,L]$, with mesh size $h = x_{j+1}-x_j$ for all $j = 0, ..., N-1$.
	For an integer $r \ge 4$, we define the finite element space
	\begin{align}
		V_h := \Big\{ \varphi \in C^{r-3}_{\mathrm{per}}([0,L]) :\;
		\varphi|_{[x_{i-1},x_i]} \in \mathbb{P}_{r-1}([x_{i-1},x_i]), \quad 1 \le i \le N \Big\},
	\end{align}
	where $\mathbb{P}_{r-1}$ denotes the space of polynomials of degree at most $r-1$, and
	\[
	C^{r-3}_{\mathrm{per}}([0,L]) := \left\{ \varphi \in C^{r-3}([0,L]) : 
	\varphi^{(m)}(0) = \varphi^{(m)}(L),\; 0 \le m \le r-3 \right\}.
	\]}
We recall the $L^2$ projection $P_h$ from $L^2(D)$ into $V_h$, which is defined as
\begin{align}\label{L2_projection}
	\bigl(v - P_h v, \phi\bigr)  = 0\qquad\forall\phi \in V_h.
\end{align}
In addition, we also have the following inequalities: 
\begin{align}\label{projection_ineq}
	\|v -  P_h v\|_{L^2} + h \|\p_x(v - P_h v)\|_{L^2} + h^2\|\p_x^2(v - P_h v)\|_{L^2} \leq C_ph^{r}\|v\|_{H^{r}},
\end{align}
for all $v \in H^{r}(D)$.

{Additionally, let $M$ be a positive integer and let
	\[
	0 = t_0 < t_1 < \cdots < t_M = T
	\]
	be a uniform partition of the time interval $[0,T]$, with the time step size $k = t_{n+1}-t_n$ for all $n = 0, ..., M-1$.}
	
	Now, we present the main algorithm of the paper.

\bigskip

{\bf Main Algorithm}

{Let $u_h^0 = P_h u_0$}.	For $n= 0, 1, ..., M-1$, we seek $u_h^{n+1} \in V_h$, with $r\geq 4$ such that
\begin{align}\label{Scheme_Standard}
	&\bigl(u_h^{n+1} - u_h^{n}, \phi_h\bigr) + \nu k\bigl(\p_x^2 u_h^{n+1}, \p^2_x \phi_h\bigr)  - k(\p_xu_h^{n+1}, \p_x\phi_h)  + k\bigl(u_h^{n+1}\p_xu_h^{n+1}, \phi_h\bigr)\\\notag
	&\qquad= (B(u_h^n)\Delta W_n,\phi_h) \qquad\forall \phi_h \in V_h,
\end{align}
{where $\Delta W_n = W(t_{n+1}) - W(t_n) \sim \mathcal{N}(0, k)$.}

\bigskip

{\begin{remark}
		At each time step, scheme \eqref{Scheme_Standard} leads to a nonlinear problem in the finite-dimensional space $V_h$. The existence of a solution $u_h^{n+1} \in V_h$ can be established by standard finite-dimensional arguments, for instance via Brouwer's fixed-point theorem. The proof follows along the same lines as in the analysis of finite element approximations of nonlinear stochastic PDEs with similar structure, such as the stochastic Navier--Stokes equations; see, e.g., \cite[Lemma 3.1]{brzezniak2013finite}.	Uniqueness follows by considering the difference of two discrete solutions corresponding to different initial data and deriving a suitable stability estimate. Since these arguments are somewhat lengthy and standard, we omit the details for brevity.
\end{remark}}

Next, we state and prove stability estimates for $\{u_h^n\}$. 
\begin{lemma}\label{lemma_2ndmoment_discrete}  Suppose that $u_0 \in L^2(\Omega; L^2(D))$ and that $B$ satisfies \eqref{Assump_Lipschitz}. Then, the fully discrete solution $\{u_h^n\}$ satisfies
	\begin{align}\label{eq3.55}
		\mE\left[\sup_{1\leq n \leq M}\|u_h^n\|^2_{L^2}\right] + k\sum_{n=1}^{M}\mE[\|\p_x^2u_h^{n}\|^2_{L^2}]\leq C_1,
	\end{align}
	where $C_1 = C(u_0,T,c_B)>0$.
\end{lemma}
\begin{proof}
	Taking $\phi_h = u_h^{n+1}$ in \eqref{Scheme_Standard} and using the identity $2a(a-b) = a^2 - b^2 + (a-b)^2$, we obtain
	\begin{align}\label{eq3.4}
		&	\frac12\left[\|u_h^{n+1}\|^2_{L^2} - \|u_h^n\|^2_{L^2} + \|u_h^{n+1} - u_h^n\|^2_{L^2}\right] + \nu k\|\p_x^2 u_h^{n+1}\|^2_{L^2}\\\nonumber
		&= k\bigl(\p_x u_h^{n+1}, \p_x u_h^{n+1}\bigr) -  k \bigl(u_h^{n+1}\p_x u_h^{n+1}, u_h^{n+1}\bigr) + (B(u_h^n)\Delta W_n,u_h^{n+1}) \\\nonumber
		&:= I_1 + I_2 + I_3.
	\end{align}
	
	Employing integration by parts, we immediately see that $I_2 =0$. Also, integration by parts and the Cauchy-Schwarz inequality yield
	\begin{align}\label{eq3.5}
		I_1 = -k\bigl(\p_x^2 u_h^{n+1}, u_h^{n+1}\bigr) \leq \frac{\nu k}{4}\|\p_x^2 u_h^{n+1}\|^2_{L^2} + \frac{k}{\nu} \|u_h^{n+1}\|^2_{L^2}.
	\end{align}
	Next, we control $I_3$ by using Cauchy-Schwarz's inequality and \eqref{Assump_Lipschitz} with $m=0$ as follows:
	\begin{align}\label{eq3.6}
		I_3 &= (B(u_h^n)\Delta W_n,u_h^{n+1} - u_h^n) +  (B(u_h^n)\Delta W_n,u_h^{n})\\\nonumber
		&\leq \frac{1}{4}\|u_h^{n+1} - u_h^n\|^2_{L^2} + \|B(u_h^n)\Delta W_n\|^2_{L^2} +  (B(u_h^n)\Delta W_n,u_h^{n})\\\nonumber
		&\leq \frac{1}{4}\|u_h^{n+1} - u_h^n\|^2_{L^2} + \|B(u_h^n)\|^2_{L^2} |\Delta W_n|^2+  (B(u_h^n)\Delta W_n,u_h^{n})\\\nonumber
		&\leq \frac{1}{4}\|u_h^{n+1} - u_h^n\|^2_{L^2} + c_B\|u_h^n\|^2_{L^2} |\Delta W_n|^2 +  (B(u_h^n)\Delta W_n,u_h^{n}).
	\end{align}
	
	Substituting all the estimates from $I_1, I_2, I_3$ into \eqref{eq3.4}, we arrive at
	\begin{align}\label{eq3.10}
		&	\frac12\left[\|u_h^{n+1}\|^2_{L^2} - \|u_h^n\|^2_{L^2} + \frac12\|u_h^{n+1} - u_h^n\|^2_{L^2}\right] + \frac{3\nu k}{4}\|\p_x^2 u_h^{n+1}\|^2_{L^2}\\\nonumber
		&\leq \frac{k}{\nu}\|u_h^{n+1}\|^2_{L^2} + c_B\|u_h^n\|^2_{L^2}|\Delta W_n|^2 +  (B(u_h^n)\Delta W_n,u_h^{n}).
	\end{align}
	Applying the summation $\sum_{n=0}^{\ell}$, for $0\leq \ell <M$, we get
	\begin{align}\label{eq3.11}
		&	\frac12\|u_h^{\ell+1}\|^2_{L^2} + \frac14\sum_{n=0}^{\ell}\|u_h^{n+1} - u_h^n\|^2_{L^2} + \frac{3\nu k}{4}\sum_{n=0}^{\ell}\|\p_x^2 u_h^{n+1}\|^2_{L^2}\\\nonumber
		&\leq \frac12\|u^0_h\|^2_{L^2} + \frac{k}{\nu}\sum_{n=0}^{\ell}\|u_h^{n+1}\|^2_{L^2} + c_B\sum_{n=0}^{\ell}\|u_h^{n}\|^2_{L^2}|\Delta W_n|^2 + \sum_{n = 0}^{\ell}(B(u_h^n)\Delta W_n,u_h^{n}).
	\end{align}
	Then, there exists $k_0>0$ such that for all $k \in (0, k_0]$
	\begin{align}\label{eq3.18}
		&	\|u_h^{\ell+1}\|^2_{L^2} + \sum_{n=0}^{\ell}\|u_h^{n+1} - u_h^n\|^2_{L^2} + 3\nu k\sum_{n=0}^{\ell}\|\p_x^2 u_h^{n+1}\|^2_{L^2}\\\nonumber
		&\leq 2\|u^0_h\|^2_{L^2} + \frac{k}{\nu}\sum_{n=0}^{\ell}\|u_h^{n}\|^2_{L^2} + 4c_B\sum_{n=0}^{\ell}\|u_h^{n}\|^2_{L^2}|\Delta W_n|^2 + 4\sum_{n = 0}^{\ell}(B(u_h^n)\Delta W_n,u_h^{n}).
	\end{align}
	
	Next, applying the expectation to \eqref{eq3.18}, and using the fact that
	\begin{align*}
		\mE\left[\|u_h^n \|^2_{L^2} |\Delta W_n|^2\right] = \mE\left[\|u_h^n\|^2_{L^2}\right]\mE\left[|\Delta W_n|^2\right] = k \mE\left[\|u_h^n\|^2_{L^2}\right]
	\end{align*}
	and that\begin{align*}
		\mE\left[(B(u_h^n)\Delta W_n,u_h^{n})\right] = 0,
	\end{align*}
	we obtain
	\begin{align}\label{eq3.19}
		&	\mE\left[\|u_h^{\ell+1}\|^2_{L^2} \right]+ \sum_{n=0}^{\ell}\mE\left[\|u_h^{n+1} - u_h^n\|^2_{L^2}\right] + 3\nu k\sum_{n=0}^{\ell}\mE\left[\|\p_x^2 u_h^{n+1}\|^2_{L^2}\right]\\\nonumber
		&\leq 2\mE\left[\|u^0_h\|^2_{L^2}\right] + 4\left(c_B + \frac{1}{\nu}\right)k\sum_{n=0}^{\ell}\mE\left[\|u_h^{n}\|^2_{L^2}\right].
	\end{align}
	Next, applying the deterministic Gronwall inequality on \eqref{eq3.19}, we get
	\begin{align}\label{eq_3.19}
		\mE[\|u_h^{\ell+1}\|^2_{L^2}] + \sum_{n=0}^{\ell}\mE[\|u_h^{n+1} - u_h^n\|^2_{L^2}] + 3\nu k\sum_{n=0}^{\ell}\mE[\|\p_x^2 u_h^{n+1}\|^2_{L^2}]\leq 2\mE[\|u^0_h\|^2_{L^2}]e^{CT}.
	\end{align}
	Now, we use \eqref{eq_3.19} to establish \eqref{eq3.55}. To the end, taking $\sup_{0\leq \ell \leq M-1}$ and then expectations on both sides of \eqref{eq3.18}, we obtain 
	\begin{align}\label{eq3.144}
		&\mE\left[\sup_{0\leq \ell \leq M-1}\|u_h^{\ell+1}\|^2_{L^2}\right] + \sum_{n=0}^{M}\mE\left[\|u_h^{n+1} - u_h^n\|^2_{L^2}\right] + 3\nu k\sum_{n=0}^{M}\mE\left[\|\p_x^2 u_h^{n+1}\|^2_{L^2}\right]\\\nonumber
		&\leq 2\mE\left[\|u^0_h\|^2_{L^2}\right] + 4\left(c_B + \frac{1}{\nu}\right)k\sum_{n=0}^{M}\mE\left[\|u_h^{n}\|^2_{L^2}\right]  + 4\mE\left[\sup_{0 \leq \ell \leq M-1}\sum_{n = 0}^{\ell}(B(u_h^n)\Delta W_n,u_h^{n})\right].
	\end{align}
	{	To control the last term of \eqref{eq3.144}, we apply the Burkholder-Davis-Gundy inequality as follows
		\begin{align}\label{eq3.1555}
			\mE\left[\sup_{0 \leq \ell \leq M-1}\sum_{n = 0}^{\ell}(B(u_h^n)\Delta W_n,u_h^{n})\right] &\leq C\mE\left[\left(k\sum_{n=0}^{M-1}\|B(u_h^n)\|^2_{L^2}\|u_h^n\|^2_{L^2}\right)^{\frac12}\right]\\\nonumber
			&\leq C\mE\left[\sup_{0 \leq n \leq M-1}\|u_h^n\|_{L^2}\left(k\sum_{n=0}^{M-1}\|B(u_h^n)\|^2_{L^2}\right)^{\frac12}\right]\\\nonumber
			&\leq \frac12 \mE\left[\|u_h^0\|^2_{L^2}\right] + \frac12\mE\left[\sup_{1 \leq n \leq M-1}\|u^n_h\|^2_{L^2}\right] \\\nonumber
			&\qquad + C\mE\left[k\sum_{n=0}^{M-1}\|B(u_h^n)\|^2_{L^2}\right]\\\nonumber
			&\leq \frac12 \mE\left[\|u_h^0\|^2_{L^2}\right] + \frac12\mE\left[\sup_{1 \leq n \leq M-1}\|u^n_h\|^2_{L^2}\right] \\\nonumber
			&\qquad + C\mE\left[k\sum_{n=0}^{M-1}\|u_h^n\|^2_{L^2}\right],
		\end{align}
		where the second term on the right-hand side of \eqref{eq3.1555} is absorbed into the left side of \eqref{eq3.144}.}
	
	Using condition \eqref{Assump_Lipschitz} with $m=0$, the estimate \eqref{eq_3.19}, and \eqref{eq3.1555}, we obtain
	\begin{align}\label{eq3.155}
		&\mE\left[\sup_{0\leq \ell \leq M-1}\|u_h^{\ell+1}\|^2_{L^2}\right] + \sum_{n=0}^{M}\mE\left[\|u_h^{n+1} - u_h^n\|^2_{L^2}\right] + \nu k\sum_{n=0}^{M}\mE\left[\|\p_x^2 u_h^{n+1}\|^2_{L^2}\right]\\\nonumber
		&\leq 5\mE\left[\|u^0_h\|^2_{L^2}\right] + 8\left(c_B + \frac{1}{\nu}\right)k\sum_{n=0}^{M}\mE\left[\|u_h^{n}\|^2_{L^2}\right]  + 8c_B\mE\left[k\sum_{n = 0}^{M}\|u_h^n\|^2_{L^2}\right]\\\nonumber
		&\leq 5\mE\left[\|u^0_h\|^2_{L^2}\right] + 16\left(2c_B + \frac{1}{\nu}\right)T\mE[\|u^0_h\|^2_{L^2}]e^{CT}:= C_1.
	\end{align}
	The proof is thus complete.

\end{proof}

\begin{lemma}\label{lemma_highmoment_discrete} Suppose that $u_0 \in L^{2^p}(\Omega; L^2(D))$ for any $0\leq p \leq 3$ and that $B$ satisifies \eqref{Assump_Lipschitz}. Then, the fully discrete solution $\{u_h^n\}$ satisfies
	\begin{align}\label{eq_3.11}
		\mE\left[\sup_{1\leq n \leq M}\|u_h^n\|^{2^p}_{L^2}\right] \leq C_p,
	\end{align}
	where {$C_p= C(u_0,T, c_B, p)>0$. }
\end{lemma}
\begin{proof} 
	We note that the case $p = 1$ was established in Lemma \ref{lemma_2ndmoment_discrete}. We only need to give the proofs for the cases $p = 2$ and $p=3$. 
	
	Multiplying \eqref{eq3.10} by $\|u_h^{n+1}\|^2_{L^2}$ and using the identity $2a(a-b) = a^2 - b^2 + (a-b)^2$, we obtain
	\begin{align}\label{eq3.12}
		&\frac{1}{4}\left[\|u_h^{n+1}\|^4_{L^2} - \|u_h^n\|^4_{L^2}\right] + \frac14 \left(\|u_h^{n+1}\|^2_{L^2} - \|u_h^n\|^2_{L^2}\right)^2 \\\nonumber
		&\qquad\qquad+ \frac14\|u_h^{n+1} - u_h^n\|^2_{L^2}\|u_h^{n+1}\|^2_{L^2} + \frac{3\nu k}{4}\|\p_x^2 u_h^{n+1}\|^2_{L^2}\|u_h^{n+1}\|^2_{L^2}\\\nonumber
		&\leq \frac{k}{\nu}\|u^{n+1}_h\|^4_{L^2} + c_B^2\|u_h^{n}\|^2_{L^2}|\Delta W_n|^2 \|u_h^{n+1}\|^2_{L^2} + \left(B(u_h^n)\Delta W_n, u_h^n\right)\|u_h^{n+1}\|^2_{L^2}\\\nonumber
		&=\frac{k}{\nu}\|u^{n+1}_h\|^4_{L^2} + c_B^2\|u_h^{n}\|^2_{L^2}|\Delta W_n|^2\left( \|u_h^{n+1}\|^2_{L^2} -  \|u_h^{n}\|^2_{L^2}\right) + c_B^2\|u_h^{n}\|^4_{L^2}|\Delta W_n|^2 \\\nonumber
		&\qquad+ \left(B(u_h^n)\Delta W_n, u_h^n\right)\left(\|u_h^{n+1}\|^2_{L^2} - \|u_h^{n}\|^2_{L^2}\right) + \left(B(u_h^n)\Delta W_n, u_h^n\right)\|u_h^{n}\|^2_{L^2}\\\nonumber
		&\leq \frac{k}{\nu}\|u^{n+1}_h\|^4_{L^2} + 4c_B^4\|u_h^{n}\|^4_{L^2}|\Delta W_n|^4 + 5c_B^2\|u_h^{n}\|^4_{L^2}|\Delta W_n|^2 \\\nonumber
		&\qquad +\frac18\left(\|u_h^{n+1}\|^2_{L^2} - \|u_h^n\|^2_{L^2}\right)^2 + \left(B(u_h^n)\Delta W_n, u_h^n\right)\|u_h^{n}\|^2_{L^2},
	\end{align}
	where the last inequality above is obtained by using the Young inequality and the assumption \eqref{Assump_Lipschitz}. 
	
	Next, applying the summation $\sum_{n = 0}^{\ell}$ to \eqref{eq3.12}, we get
	\begin{align}\label{eq_3.13}
		&\frac{1}{4}\|u_h^{\ell+1}\|^4_{L^2} + \frac18 \sum_{n = 0}^{\ell}\left(\|u_h^{n+1}\|^2_{L^2} - \|u_h^n\|^2_{L^2}\right)^2 \\\nonumber
		&\leq \frac14 \|u_h^0\|^4_{L^2}+ \frac{k}{\nu}\sum_{n = 0}^{\ell}\|u^{n+1}_h\|^4_{L^2}  + 4c_B^4\sum_{n = 0}^{\ell}\|u_h^{n}\|^4_{L^2}|\Delta W_n|^4 \\\nonumber
		&\qquad+ 5c_B^2\sum_{n = 0}^{\ell}\|u_h^{n}\|^4_{L^2}|\Delta W_n|^2 + \sum_{n = 0}^{\ell}  \left(B(u_h^n)\Delta W_n, u_h^n\right)\|u_h^{n}\|^2_{L^2}.
	\end{align}
	
	Then, there exists a constant $k_1>0$ such that for any $k\in (0,k_1]$
	\begin{align}\label{eq319}
		&\|u_h^{\ell+1}\|^4_{L^2} +  \sum_{n = 0}^{\ell}\left(\|u_h^{n+1}\|^2_{L^2} - \|u_h^n\|^2_{L^2}\right)^2 \\\nonumber
		&\leq  2\|u_h^0\|^4_{L^2}+ \frac{8k}{\nu}\sum_{n = 0}^{\ell}\|u^{n}_h\|^4_{L^2}  + 32c_B^4\sum_{n = 0}^{\ell}\|u_h^{n}\|^4_{L^2}|\Delta W_n|^4 \\\nonumber
		&\qquad+ 40c_B^2\sum_{n = 0}^{\ell}\|u_h^{n}\|^4_{L^2}|\Delta W_n|^2 + 8\sum_{n = 0}^{\ell}  \left(B(u_h^n)\Delta W_n, u_h^n\right)\|u_h^{n}\|^2_{L^2}.
	\end{align}
	
	Taking the expectation and then applying the discrete Gronwall inequality to \eqref{eq319}, we obtain 
	\begin{align}\label{eq3.200}
		&\mE\left[\|u_h^{\ell+1}\|^4_{L^2} \right]+  \sum_{n = 0}^{\ell}\mE\left[\left(\|u_h^{n+1}\|^2_{L^2} - \|u_h^n\|^2_{L^2}\right)^2\right] \leq 2\mE\left[\|u_h^0\|^4_{L^2}\right] e^{CT}.
	\end{align}
	
	{Now, we use \eqref{eq3.200} to establish \eqref{eq_3.11} with $p=2$. To the end, taking $\sup_{0\leq \ell \leq M-1}$ and then expectations on both sides of \eqref{eq319}, we obtain 
		\begin{align}\label{eq3.222}
			&\mE\left[\sup_{0 \leq \ell \leq M-1}\|u_h^{\ell+1}\|^4_{L^2}\right] + \mE\left[\sum_{n = 0}^{M-1}\left(\|u_h^{n+1}\|^2_{L^2} - \|u_h^n\|^2_{L^2}\right)^2\right]\\\nonumber
			&\leq  2\mE\left[\|u_h^0\|^4_{L^2}\right]+ \frac{8k}{\nu}\sum_{n = 0}^{M-1}\mE\left[\|u^{n}_h\|^4_{L^2}\right]  + 32c_B^4k^2\sum_{n = 0}^{M-1}\mE\left[\|u_h^{n}\|^4_{L^2}\right] \\\nonumber
			&\qquad+ 40c_B^2k\sum_{n = 0}^{M-1}\mE\left[\|u_h^{n}\|^4_{L^2}\right] + 8\mE\left[\sup_{0 \leq \ell \leq M-1}\sum_{n = 0}^{\ell}  \left(B(u_h^n)\Delta W_n, u_h^n\right)\|u_h^{n}\|^2_{L^2}\right].
		\end{align}
		To control the last term of \eqref{eq3.22}, we use the Burkholder-Davis-Gundy inequality as follows
		\begin{align}\label{eq3.233}
			&\mE\left[\sup_{0 \leq \ell \leq M-1}\sum_{n = 0}^{\ell}  \left(B(u_h^n)\Delta W_n, u_h^n\right)\|u_h^{n}\|^2_{L^2}\right] \\\nonumber
			&\leq C\mE\left[\left(k\sum_{n=0}^{M-1}\|B(u_h^n)\|^2_{L^2}\|u_h^n\|^6_{L^2}\right)^{\frac12}\right]\\\nonumber
			&\leq C\mE\left[\sup_{0 \leq n \leq M-1}\|u_h^n\|^2_{L^2}\left(k\sum_{n=0}^{M-1}\|B(u_h^n)\|^2_{L^2}\|u^n_h\|^2_{L^2}\right)^{\frac12}\right]\\\nonumber
			&\leq \frac12 \mE\left[\|u_h^0\|^4_{L^2}\right] + \frac12\mE\left[\sup_{1 \leq n \leq M-1}\|u^n_h\|^4_{L^2}\right] + C\mE\left[k\sum_{n=0}^{M-1}\|B(u_h^n)\|^2_{L^2}\|u_h^n\|^2_{L^2}\right]\\\nonumber
			&\leq \frac12 \mE\left[\|u_h^0\|^4_{L^2}\right] + \frac12\mE\left[\sup_{1 \leq n \leq M-1}\|u^n_h\|^4_{L^2}\right] + C\mE\left[k\sum_{n=0}^{M-1}\|u_h^n\|^4_{L^2}\right],
		\end{align}
		where the second term on the right-hand side of \eqref{eq3.233} is absorbed into the left side of \eqref{eq3.222}.
		
		Now, substituting \eqref{eq3.233} into \eqref{eq3.222} and then using the estimate \eqref{eq3.200}, we obtain the desired estimate \eqref{eq_3.11} with $p=2$.
	}
	Similarly, we can obtain \eqref{eq_3.11} for $p =3$ by multiplying \eqref{eq3.12} by $\|u_h^{n+1}\|^4_{L^2}$ and proceed as shown above.
	
	The proof is complete.
\end{proof}

Next, we present the error estimates of the scheme \eqref{Scheme_Standard} in two cases: bounded multiplicative noise and general multiplicative noise. 

\subsection{Full expectation error estimates in the case of bounded multiplicative noise}\label{sub-sec3.2}

In this part, we derive strong error estimates for the full discrete solution $\{u_h^n\}$ in the case of bounded noise, i.e., $B$ satisfies the assumption \eqref{Assump_Lineargrowth}. Examples of such a diffusion include the instances $B(u) = \sin(u), \cos(u)$, or $\frac{u^2}{u^2 + 1}$.  With the assumption, the errors will be computed with full expectations, which provide strong convergence in suitable $L^p$ norms.
\subsubsection{Sub-second moment error estimates}
First, we derive the optimal error estimates of $\{u_h^n\}$ in the  $L_{\omega}^{p}L_t^{\infty}L_x^2$- and $L_{\omega}^{p}L_t^{2}H_x^2$-norm for $0<p<2$, which are the sub-second moment error estimates.

\begin{theorem}\label{Theorem_sub_moment} 
	Let $u$ be the variational solution to \eqref{Weak_formulation} and $\{u_h^{n}\}_{n=1}^M$ be generated by \eqref{Scheme_Standard}. Suppose that $B$ satisfies conditions \eqref{Assump_Lineargrowth} and \eqref{Assump_Lipschitz} and that $u_0 \in \cap_{i=0}^{r} L^{2^{r+3-i}}(\Omega; H^i(D))$ for any integers $r \geq 4$. 
	Additionally, for any $0 < q < \frac{99}{100}$, assume that $L_0 <\frac{{\nu}\pi^2}{170 C_eL^2\sqrt{ q}}$ and $\mE\left[\exp\left(16\kappa \|u_0\|^2_{L^2}\right)\right] < \infty$, where $\kappa = \frac{900C_e^2 L^4q}{\nu^2\pi^4}$, and $L_0$ and $C_e$ are respectively the constants as in \eqref{Assump_Lineargrowth} and \eqref{sobolev_ineq}. Then, there holds
	\begin{align}\label{equu310}
		&\left(\mE\left[\max_{1\leq n\leq M}\|u(t_n) - u_h^n\|^{2q}_{L^2}\right]\right)^{\frac{1}{2q}} \\\nonumber
		&\qquad\qquad\qquad+ \left(\mE\left[\left(\nu k\sum_{n = 1}^M\|\p_x^2(u(t_n) - u_h^n)\|^2_{L^2}\right)^q\right]\right)^{\frac{1}{2q}}
		\leq \widehat{C}_1\,\left(k^{\frac12} + h^{r-2}\right),
	\end{align}
	where $\widehat{C}_1 = C(q, u_0, T, C_B)$ is a positive constant.
\end{theorem}
\begin{proof}
	Denote $e^n:= u(t_n) - u_h^n = \theta^n + \varepsilon^n$, where
	\begin{align*}
		\theta^n = u(t_n) - P_h u(t_n),\qquad \varepsilon^n = P_h u(t_n) - u_h^n.
	\end{align*}
	
	Subtracting \eqref{Weak_formulation} from \eqref{Scheme_Standard}, we obtain for all $\phi_h\in V_h$
	\begin{align}\label{eq3.14}
		&\bigl(e^{n+1} - e^n, \phi_h\bigr) + \nu k\bigl(\p_x^2e^{n+1}, \p_x^2 \phi_h\bigr) - k\bigl(\p_x e^{n+1}, \p_x \phi_h\bigr) \\\nonumber
		&= \nu\int_{t_n}^{t_{n+1}} (\p_x^2(u(t_{n+1}) - u(s)), \p_x^2\phi_h)\, ds - \int_{t_n}^{t_{n+1}} (\p_x(u(t_{n+1}) - u(s)), \p_x\phi_h)\, ds \\\nonumber
		&\quad -\int_{t_n}^{t_{n+1}}(u(s)\p_x u(s) - u(t_{n+1})\p_xu(t_{n+1}), \phi_h)\, ds  \\\nonumber
		&\quad- k(u(t_{n+1})\p_xu(t_{n+1}) - u_h^{n+1}\p_xu_h^{n+1}, \phi_h) + \left(\int_{t_n}^{t_{n+1}} (B(u(s)) - B(u_h^n))\, dW(s), \phi_h\right).
	\end{align}

	Using the $L^2$ projection orthogonality \eqref{L2_projection}, we can recast the left-hand side of \eqref{eq3.14} as follow:
	{\begin{align*}
			\bigl(e^{n+1} - e^n, \phi_h\bigr) + \nu k\bigl(\p_x^2e^{n+1}, \p_x^2 \phi_h\bigr) - k\bigl(\p_x e^{n+1}, \p_x \phi_h\bigr) &= (\varepsilon^{n+1} - \varepsilon^n, \phi_h) + \nu k(\p_x^2 \varepsilon^{n+1} , \p_x^2 \phi_h) \\\nonumber&\qquad+\nu k (\p_x^2\theta^{n+1}, \p_x^2\phi_h) - k (\p_x\varepsilon^{n+1}, \p_x\phi_h) \\\nonumber
			&\qquad- k (\p_x\theta^{n+1},\p_x\phi_h).
	\end{align*}}
	
	Substituting this into \eqref{eq3.14}, taking $\phi_h = \varepsilon^{n+1}\in V_h$ and using the identity $2a(a-b) = a^2-b^2 + (a-b)^2$, we obtain
	\begin{align}\label{eq3.16}
		&\frac12\bigl[\|\varepsilon^{n+1}\|^2_{L^2} - \|\varepsilon^n\|^2_{L^2} + \|\varepsilon^{n+1} - \varepsilon^n\|^2_{L^2}\bigr]+ \nu k\|\p_x^2\varepsilon^{n+1}\|^2_{L^2} \\\nonumber
		&= -\nu k (\p_x^2\theta^{n+1}, \p_x^2\varepsilon^{n+1}) + k (\p_x\varepsilon^{n+1}, \p_x\varepsilon^{n+1}) + k (\p_x\theta^{n+1},\p_x\varepsilon^{n+1})\\\nonumber
		&\quad+\nu\int_{t_n}^{t_{n+1}} (\p_x^2(u(t_{n+1}) - u(s)), \p_x^2\varepsilon^{n+1})\, ds - \int_{t_n}^{t_{n+1}} (\p_x(u(t_{n+1}) - u(s)), \p_x\varepsilon^{n+1})\, ds \\\nonumber
		&\quad -\int_{t_n}^{t_{n+1}}(u(s)\p_x u(s) - u(t_{n+1})\p_xu(t_{n+1}), \varepsilon^{n+1})\, ds  \\\nonumber
		&\quad- k(u(t_{n+1})\p_xu(t_{n+1}) - u_h^{n+1}\p_xu_h^{n+1}, \varepsilon^{n+1}) \\\nonumber
		&\quad+ \left(\int_{t_n}^{t_{n+1}} (B(u(s)) - B(u_h^n))\, dW(s), \varepsilon^{n+1}\right)\\\nonumber
		&:= Z_1 + Z_2 + Z_3 + Z_4 + Z_5 + Z_6 + Z_7 + Z_8.
	\end{align}
	
	Now, we proceed to estimate each term on the right-hand side of \eqref{eq3.16}. Concerning $Z_1$, using the Cauchy-Schwarz inequality and \eqref{projection_ineq}, we have
	\begin{align*}
		Z_1 &\leq 2\nu k \|\p_x^2 \theta^{n+1}\|^2_{L^2}  + \frac{\nu k}{8}\|\p_x^2\varepsilon^{n+1}\|^2_{L^2}\\\nonumber
		&\leq 2\nu k C_p^2h^{2(r-2)}\|u(t_{n+1})\|^2_{H^{r}}  + \frac{\nu k}{8}\|\p_x^2\varepsilon^{n+1}\|^2_{L^2}.
	\end{align*}
	Next, using integration by parts and \eqref{projection_ineq}, we estimate $Z_2$ and $Z_3$ as follows:
	\begin{align*}
		Z_2 &= - k(\varepsilon^{n+1}, \p_x^2\varepsilon^{n+1}) \leq \frac{\nu k}{8} \|\p_x^2\varepsilon^{n+1}\|^2_{L^2} + \frac{2k}{\nu}\|\varepsilon^{n+1}\|^2_{L^2}.\\\nonumber
		Z_3 &= -k (\theta^{n+1}, \p_x^2\varepsilon^{n+1}) \\\nonumber
		&\leq \frac{\nu k}{8} \|\p_x^2\varepsilon^{n+1}\|^2_{L^2} + \frac{2k}{\nu}\|\theta^{n+1}\|^2_{L^2}  \leq \frac{\nu k}{8} \|\p_x^2\varepsilon^{n+1}\|^2_{L^2} + \frac{2k}{\nu}C_p^2 h^{2r}\|u(t_{n+1})\|^2_{H^{r}}.
	\end{align*}
	Turning to $Z_4$, we invoke the Cauchy-Schwarz inequality to infer
	\begin{align*}
		Z_4 &\leq \frac{\nu k}{8}\|\p_x^2\varepsilon^{n+1}\|^2_{L^2} + 2\nu\int_{t_n}^{t_{n+1}}\|\p_x^2(u(t_{n+1}) - u(s))\|^2_{L^2}\, ds.
	\end{align*}
	With regard to $Z_5$, once again, we employ  integration by parts to obtain
	\begin{align*}
		Z_5 &=\int_{t_n}^{t_{n+1}} (\p_x^2(u(t_{n+1}) - u(s)), \varepsilon^{n+1})\, ds\\\nonumber
		&\leq k\|\varepsilon^{n+1}\|^2_{L^2} + \frac{1}{4} \int_{t_n}^{t_{n+1}}\|\p_x^2(u(t_{n+1}) - u(s))\|^2_{L^2}\, ds.
	\end{align*}
	Next, we estimate $Z_6$ by adding and subtracting the term $u(t_{n+1})\p_xu(s)$ and using \eqref{sobolev_ineq} as follow:
	\begin{align*}
		Z_6 &=-\int_{t_n}^{t_{n+1}}\bigl([u(s)-u(t_{n+1})]\p_x u(s) + u(t_{n+1})\p_x[u(t_{n+1}) - u(s)], \varepsilon^{n+1}\bigr)\, ds  \\\nonumber
		&\leq k\|\varepsilon^{n+1}\|^2_{L^2} + \frac14\int_{t_n}^{t_{n+1}}\|[u(s)-u(t_{n+1})]\p_x u(s) + u(t_{n+1})\p_x[u(t_{n+1}) - u(s)]\|^2_{L^2}\, ds\\\nonumber
		&\leq k\|\varepsilon^{n+1}\|^2_{L^2} + \frac12\int_{t_n}^{t_{n+1}}\left(\|[u(s)-u(t_{n+1})]\p_x u(s)\|^2_{L^2} + \|u(t_{n+1})\p_x[u(t_{n+1}) - u(s)]\|^2_{L^2}\right)\, ds\\\nonumber
		&\leq k\|\varepsilon^{n+1}\|^2_{L^2} + \frac12\int_{t_n}^{t_{n+1}}\left(\|u(s)-u(t_{n+1})\|^2_{L^{\infty}}\|\p_x u(s)\|^2_{L^2} \right.\\\nonumber
		&\qquad+ \left.\|u(t_{n+1})\|^2_{L^{\infty}}\|\p_xu(t_{n+1}) - u(s)\|^2_{L^2}\right)\, ds\\\nonumber
		&\leq k\|\varepsilon^{n+1}\|^2_{L^2} + \frac{C_e^2}{2}\int_{t_n}^{t_{n+1}}\left(\|\p_xu(s)\|^2_{L^2} + \|\p_x u(t_{n+1})\|^2_{L^2}\right)\|\p_xu(t_{n+1}) - u(s)\|^2_{L^2}\, ds.
	\end{align*}
	To estimate $Z_7$, we add and subtract the term $u_h^{n+1}\p_xu(t_{n+1})$ to get
	\begin{align*}
		Z_7 &= -k\bigl(e^{n+1}\p_x u(t_{n+1}) + u_h^{n+1}\p_xe^{n+1}, \varepsilon^{n+1}\bigr)\\\nonumber
		&= -k\bigl(\theta^{n+1}\p_x u(t_{n+1}) + u_h^{n+1}\p_x\theta^{n+1}, \varepsilon^{n+1}\bigr) -k\bigl(\varepsilon^{n+1}\p_x u(t_{n+1}) + u_h^{n+1}\p_x\varepsilon^{n+1}, \varepsilon^{n+1}\bigr)\\\nonumber
		&:= Z_{7,1} + Z_{7,2}.
	\end{align*}
	
	Using \eqref{projection_ineq}, we obtain
	\begin{align*}
		Z_{7,1} &\leq k\|\varepsilon^{n+1}\|^2_{L^2} + \frac{k}{2}\left[\|\theta^{n+1}\p_x u(t_{n+1})\|^2_{L^2} + \|u_h^{n+1}\p_x\theta^{n+1}\|^2_{L^2}\right]\\\nonumber
		&\leq k\|\varepsilon^{n+1}\|^2_{L^2} + \frac{k}{2}\left[\|\theta^{n+1}\|^2_{L^{\infty}}\|\p_x u(t_{n+1})\|^2_{L^2} + \|u_h^{n+1}\|^2_{L^2}\|\p_x\theta^{n+1}\|^2_{L^{\infty}}\right]\\\nonumber
		&\leq k\|\varepsilon^{n+1}\|^2_{L^2} + \frac{C_e^2k}{2}\left[\|\theta^{n+1}\|^2_{H^1}\|\p_x u(t_{n+1})\|^2_{L^2} + \|u_h^{n+1}\|^2_{L^2}\|\theta^{n+1}\|^2_{H^2}\right]\\\nonumber
		&\leq k\|\varepsilon^{n+1}\|^2_{L^2} + \frac{C_e^2k}{2}C_p^2\left[h^{2(r-1)}\|u(t_{n+1})\|^2_{H^{r}}\|\p_x u(t_{n+1})\|^2_{L^2} + h^{2(r-2)}\|u_h^{n+1}\|^2_{L^2}\|u(t_{n+1})\|^2_{H^{r}}\right].
	\end{align*}
	With regard to $Z_{7,2}$, using the fact that $\bigl(\varepsilon^{n+1}\p_x\varepsilon^{n+1},\varepsilon^{n+1}\bigr) = 0$, we get
	\begin{align*}
		\bigl(u_h^{n+1}\p_x\varepsilon^{n+1},\varepsilon^{n+1}\bigr) = \bigl(P_h[u(t_{n+1})]\p_x\varepsilon^{n+1},\varepsilon^{n+1}\bigr).
	\end{align*}
	In turn, we can update $Z_{7,2}$ as follows:
	\begin{align*}
		Z_{7,2} &= -k\bigl(\varepsilon^{n+1}\p_x u(t_{n+1}) + P_h[u(t_{n+1})]\p_x\varepsilon^{n+1}, \varepsilon^{n+1}\bigr).
	\end{align*}
	As a consequence, using integration by parts, the stability of the $L^2$ projection and \eqref{sobolev_ineq}, we get
	\begin{align*}
		Z_{7,2} &= 2k\bigl(u(t_{n+1}), \p_x\varepsilon^{n+1}\varepsilon^{n+1}\bigr) -k\bigl(P_h[u(t_{n+1})]\p_x\varepsilon^{n+1}, \varepsilon^{n+1}\bigr)\\\nonumber
		&\leq 2k\|u(t_{n+1})\|_{L^2}\|\varepsilon^{n+1}\|_{L^2}\|\p_x \varepsilon^{n+1}\|_{L^{\infty}} + k\|P_h[u(t_{n+1})]\|_{L^2}\|\varepsilon^{n+1}\|_{L^2}\|\p_x \varepsilon^{n+1}\|_{L^{\infty}}\\\nonumber
		&\leq 3k\|u(t_{n+1})\|_{L^2}\|\varepsilon^{n+1}\|_{L^2}\|\p_x \varepsilon^{n+1}\|_{L^{\infty}}\\\nonumber
		&\leq 3kC_e\|u(t_{n+1})\|_{L^2}\|\varepsilon^{n+1}\|_{L^2}\|\p_x^2 \varepsilon^{n+1}\|_{L^{2}}\\\nonumber
		&\leq \frac{\nu k}{8}\|\p_x^2 \varepsilon^{n+1}\|^2_{L^{2}} + \frac{18 k C_e^2}{\nu}\|u(t_{n+1})\|^2_{L^2}\|\varepsilon^{n+1}\|^2_{L^2}.
	\end{align*}
	
	Now, we estimate the noise term $Z_8a$ as follows:
	\begin{align*}
		Z_{8} &=  \left(\int_{t_n}^{t_{n+1}} (B(u(s)) - B(u_h^n))\, dW(s), \varepsilon^{n+1} - \varepsilon^n\right) +  \left(\int_{t_n}^{t_{n+1}} (B(u(s)) - B(u_h^n))\, dW(s),  \varepsilon^n\right)\\\nonumber
		&=  \left(\int_{t_n}^{t_{n+1}} (B(u(s)) - B(u(t_{n})))\, dW(s), \varepsilon^{n+1} - \varepsilon^n\right) +\bigl((B(u(t_n)) - B(u_h^n))\Delta W_n, \varepsilon^{n+1} - \varepsilon^n\bigr) \\\nonumber
		&\qquad+  \left(\int_{t_n}^{t_{n+1}} (B(u(s)) - B(u^n))\, dW(s),  \varepsilon^n\right)\\\nonumber
		&\leq  2\left\|\int_{t_n}^{t_{n+1}} (B(u(s)) - B(u(t_{n})))\, dW(s)\right\|^2_{L^2} +2\left\|(B(u(t_n)) - B(u_h^n))\Delta W_n\right\|^2_{L^2} + \frac14\|\varepsilon^{n+1} - \varepsilon^{n}\|^2_{L^2} \\\nonumber
		&\qquad+  \left(\int_{t_n}^{t_{n+1}} (B(u(s)) - B(u^n))\, dW(s),  \varepsilon^n\right).
	\end{align*}
	
	In order to apply the stochastic Gronwall inequality in Lemma \ref{Stochastic_Gronwall}, we add and subtract the term $2k\|B(u(t_n)) - B(u_h^n)\|^2_{L^2}$ on the right-hand side of $Z_8$, and then using the assumption \eqref{Assump_Lipschitz} as follows:
	\begin{align*}
		Z_8 &\leq  2\left\|\int_{t_n}^{t_{n+1}} (B(u(s)) - B(u(t_{n})))\, dW(s)\right\|^2_{L^2}  + \frac14\|\varepsilon^{n+1} - \varepsilon^{n}\|^2_{L^2} \\\nonumber
		&\qquad+2\left\|(B(u(t_n)) - B(u_h^n))\Delta W_n\right\|^2_{L^2} - 2k\|B(u(t_n)) - B(u_h^n)\|^2_{L^2} + 2k\|B(u(t_n)) - B(u_h^n)\|^2_{L^2}\\\nonumber
		&\qquad+  \left(\int_{t_n}^{t_{n+1}} (B(u(s)) - B(u^n))\, dW(s),  \varepsilon^n\right)\\\nonumber
		&\leq  2\left\|\int_{t_n}^{t_{n+1}} (B(u(s)) - B(u(t_{n})))\, dW(s)\right\|^2_{L^2}  + \frac14\|\varepsilon^{n+1} - \varepsilon^{n}\|^2_{L^2} \\\nonumber
		&\qquad+2\left\|(B(u(t_n)) - B(u_h^n))\Delta W_n\right\|^2_{L^2} - 2k\|B(u(t_n)) - B(u_h^n)\|^2_{L^2} + 2kC_B^2\|e^{n}\|^2_{L^2}\\\nonumber
		&\qquad+  \left(\int_{t_n}^{t_{n+1}} (B(u(s)) - B(u^n))\, dW(s),  \varepsilon^n\right)\\\nonumber
		&\leq  2\left\|\int_{t_n}^{t_{n+1}} (B(u(s)) - B(u(t_{n})))\, dW(s)\right\|^2_{L^2}  + \frac14\|\varepsilon^{n+1} - \varepsilon^{n}\|^2_{L^2} \\\nonumber
		&\qquad+2\left\|(B(u(t_n)) - B(u_h^n))\Delta W_n\right\|^2_{L^2} - 2k\|B(u(t_n)) - B(u_h^n)\|^2_{L^2} \\\nonumber
		&\qquad+ 2kC_B^2\|\varepsilon^{n}\|^2_{L^2} + 2kC_B^2C_p^2 h^{2r}\|u(t_n)\|^2_{H^{r}}\\\nonumber
		&\qquad+  \left(\int_{t_n}^{t_{n+1}} (B(u(s)) - B(u^n))\, dW(s),  \varepsilon^n\right).
	\end{align*}

	Now, substituting all of the estimates on $Z_1, ..., Z_8$ into \eqref{eq3.16} and absorbing the like terms to the left-hand side of \eqref{eq3.16}, we obtain
	\begin{align}\label{eq_3.14}
		&\frac{1}{2}\left[\|\varepsilon^{n+1}\|^2_{L^2} - \|\varepsilon^n\|^2_{L^2}\right] + \frac14\|\varepsilon^{n+1} - \varepsilon^n\|^2_{L^2} + \frac{3\nu k}{8}\|\p_x^2\varepsilon^{n+1}\|^2_{L^2} \\\nonumber
		&\leq 2kC_p^2h^{2(r-2)}\left(\nu + \frac{1}{\nu}\right)\|u(t_{n+1})\|^2_{H^{r}} + 2kC_B^2C_p^2 h^{2r}\|u(t_n)\|^2_{H^{r}}\\\nonumber&\qquad+ \frac{C_e^2k}{2}C_p^2\left[h^{2(r-1)}\|u(t_{n+1})\|^2_{H^{r}}\|\p_x u(t_{n+1})\|^2_{L^2} + h^{2(r-2)}\|u_h^{n+1}\|^2_{L^2}\|u(t_{n+1})\|^2_{H^{r}}\right]\\\nonumber
		&\qquad+ 2\nu\int_{t_n}^{t_{n+1}}\|\p_x^2(u(t_{n+1}) - u(s))\|^2_{L^2}\, ds + \frac{1}{4} \int_{t_n}^{t_{n+1}}\|\p_x^2(u(t_{n+1}) - u(s))\|^2_{L^2}\, ds\\\nonumber
		&\qquad + \frac{C_e^2}{2}\int_{t_n}^{t_{n+1}}\left(\|\p_xu(s)\|^2_{L^2} + \|\p_x u(t_{n+1})\|^2_{L^2}\right)\|\p_xu(t_{n+1}) - u(s)\|^2_{L^2}\, ds \\\nonumber
		&\qquad+ 2\left\|\int_{t_n}^{t_{n+1}} (B(u(s)) - B(u(t_{n})))\, dW(s)\right\|^2_{L^2} \\\nonumber
		&\qquad + k\left( 3 + \frac{2}{\nu} + \frac{18C_e^2}{\nu}\|u(t_{n+1})\|^2_{L^2} \right)\|\varepsilon^{n+1}\|^2_{L^2} + 2kC_B^2\|\varepsilon^{n}\|^2_{L^2}\\\nonumber
		&\qquad +2\left\|(B(u(t_n)) - B(u_h^n))\Delta W_n\right\|^2_{L^2} - 2k\|B(u(t_n)) - B(u_h^n)\|^2_{L^2} \\\nonumber
		&\qquad+  \left(\int_{t_n}^{t_{n+1}} (B(u(s)) - B(u^n))\, dW(s),  \varepsilon^n\right)\\\nonumber
			&\leq 2kC_p^2h^{2(r-2)}\left(\nu + \frac{1}{\nu}\right)\|u(t_{n+1})\|^2_{H^{r}} + 2kC_B^2C_p^2 h^{2r}\|u(t_n)\|^2_{H^{r}}\\\nonumber&\qquad+ \frac{C_e^2k}{2}C_p^2\left[h^{2(r-1)}\|u(t_{n+1})\|^2_{H^{r}}\|\p_x u(t_{n+1})\|^2_{L^2} + h^{2(r-2)}\|u_h^{n+1}\|^2_{L^2}\|u(t_{n+1})\|^2_{H^{r}}\right]\\\nonumber
		&\qquad+ 2\nu\int_{t_n}^{t_{n+1}}\|\p_x^2(u(t_{n+1}) - u(s))\|^2_{L^2}\, ds + \frac{1}{4} \int_{t_n}^{t_{n+1}}\|\p_x^2(u(t_{n+1}) - u(s))\|^2_{L^2}\, ds\\\nonumber
		&\qquad + \frac{C_e^2}{2}\int_{t_n}^{t_{n+1}}\left(\|\p_xu(s)\|^2_{L^2} + \|\p_x u(t_{n+1})\|^2_{L^2}\right)\|\p_xu(t_{n+1}) - u(s)\|^2_{L^2}\, ds \\\nonumber
		&\qquad+ 2\left\|\int_{t_n}^{t_{n+1}} (B(u(s)) - B(u(t_{n})))\, dW(s)\right\|^2_{L^2} \\\nonumber
		&\qquad + \int_{t_{n}}^{t_{n+1}}\left( 3 + \frac{2}{\nu} + \frac{36C_e^2}{\nu}\|u(t_{n+1}) - u(s)\|^2_{L^2} \right)\, ds\|\varepsilon^{n+1}\|^2_{L^2} + 2kC_B^2\|\varepsilon^{n}\|^2_{L^2} \\\nonumber
		&\qquad + \int_{t_{n}}^{t_{n+1}}\left( 3 + \frac{2}{\nu} + \frac{36C_e^2}{\nu}\| u(s)\|^2_{L^2} \right)\, ds\|\varepsilon^{n+1}\|^2_{L^2}\\\nonumber
		&\qquad +2\left\|(B(u(t_n)) - B(u_h^n))\Delta W_n\right\|^2_{L^2} - 2k\|B(u(t_n)) - B(u_h^n)\|^2_{L^2} \\\nonumber
		&\qquad+  \left(\int_{t_n}^{t_{n+1}} (B(u(s)) - B(u^n))\, dW(s),  \varepsilon^n\right).
	\end{align}
	
	Applying the summation $\sum_{n=0}^{\ell}$ for any $0\leq \ell <M$, we obtain
	\begin{align}\label{eq3.15}
		&	\|\varepsilon^{\ell+1}\|^2_{L^2} + \frac12\sum_{n=0}^{\ell}\|\varepsilon^{n+1} - \varepsilon^n\|^2_{L^2} + \frac{3\nu k}{4}\sum_{n=0}^{\ell} \|\p_x^2 \varepsilon^{n+1}\|^2_{L^2} \\\nonumber
		&\leq F_{\ell} + M_{\ell} + \sum_{n=0}^{\ell} G_n\|\varepsilon^{n}\|^2_{L^2},
	\end{align}
	
	where 
	\begin{align*}
		F_{\ell} &:= \sum_{n=0}^{\ell}\left[4kC_p^2h^{2(r-2)}\left(\nu + \frac{h^4}{\nu}\right)\|u(t_{n+1})\|^2_{H^{r}} + 4kC_B^2C_p^2 h^{2r}\|u(t_n)\|^2_{H^{r}}\right.\\\nonumber&\qquad\left.+ C_e^2kC_p^2\left[h^{2(r-1)}\|u(t_{n+1})\|^2_{H^{r}}\|\p_x u(t_{n+1})\|^2_{L^2} + h^{2(r-2)}\|u_h^{n+1}\|^2_{L^2}\|u(t_{n+1})\|^2_{H^{r}}\right]\right.\\\nonumber
		&\qquad\left.+ 4\nu\int_{t_n}^{t_{n+1}}\|\p_x^2(u(t_{n+1}) - u(s))\|^2_{L^2}\, ds + \frac{1}{2} \int_{t_n}^{t_{n+1}}\|\p_x^2(u(t_{n+1}) - u(s))\|^2_{L^2}\, ds\right.\\\nonumber
		&\qquad\left. + C_e^2\int_{t_n}^{t_{n+1}}\left(\|\p_xu(s)\|^2_{L^2} + \|\p_x u(t_{n+1})\|^2_{L^2}\right)\|\p_xu(t_{n+1}) - u(s)\|^2_{L^2}\, ds\right. \\\nonumber
		&\qquad\left.+ 4\left\|\int_{t_n}^{t_{n+1}} (B(u(s)) - B(u(t_{n})))\, dW(s)\right\|^2_{L^2} \right.\\\nonumber
		&\qquad\left.+ \int_{t_{n}}^{t_{n+1}}\left( 3 + \frac{2}{\nu} + \frac{36C_e^2}{\nu}\|u(t_{n+1}) - u(s)\|^2_{L^2} \right)\, ds\|\varepsilon^{n+1}\|^2_{L^2}\right] \\\nonumber
		&\qquad + \int_{t_{\ell}}^{t_{\ell+1}}\left( 3 + \frac{2}{\nu} + \frac{36C_e^2}{\nu}\| u(s)\|^2_{L^2} \right)\, ds\|\varepsilon^{\ell+1}\|^2_{L^2} ,\\\nonumber
		M_{\ell} &:= \sum_{n=0}^{\ell} Z_{n},\\\nonumber
		Z_n&:= 4\left\|(B(u(t_n)) - B(u_h^n))\Delta W_n\right\|^2_{L^2} - 2k\|B(u(t_n)) - B(u_h^n)\|^2_{L^2} \\\nonumber
		&\qquad+ 2 \left(\int_{t_n}^{t_{n+1}} (B(u(s)) - B(u^n))\, dW(s),  \varepsilon^n\right),\\\nonumber
		G_n&:=  2\int_{t_{n}}^{t_{n+1}}\left( 3 + \frac{2}{\nu} +2C_B^2 + \frac{36C_e^2}{\nu}\|u(s)\|^2_{L^2} \right)\,ds.
	\end{align*}
	
	Suppose that $\{M_{\ell};\, \ell\geq 0\}$ is a Martingale  (this fact will be verified at the end of the proof), using the stochastic Gronwall's inequality \eqref{ineq2.4} with $\alpha = \frac{100}{99}, \beta = 100$,  $0<q <\frac{99}{100}$ to \eqref{eq3.15}, we obtain
	\begin{align}\label{eq3.13}
		&\Bigl(\mE\bigl[\sup_{0 \leq \ell \leq M}\|\varepsilon^{\ell}\|^{2q}_{L^2}\bigr]\Bigr)^{\frac{1}{2q}} + \left(\mE\left[\left(\frac{3\nu}{4}k\sum_{n = 0}^M\|\p_x^2\varepsilon^{n}\|^2_{L^2}\right)^q\right]\right)^{\frac{1}{2q}}\\\nonumber
		&\qquad \quad \leq \Bigl(1+ \frac{1}{1- \alpha q}\Bigr)^{\frac{1}{2\alpha q}} \left(\mE\left[\exp\left( \beta q\sum_{n = 0}^{M-1}{G}_n\right)\right]\right)^{\frac{1}{2\beta q}} \,\Bigl(\mE\Bigl[\sup_{0 \leq \ell < M} F_{\ell}\Bigr]\Bigr)^{\frac12}.
	\end{align}
	
	Now, we proceed to estimate the right-hand side of \eqref{eq3.13}. Using  Lemma \ref{lemma_expo_moment_H}, we control the second term as follows:
	\begin{align*}
		\mE\left[\exp\left( \beta q\sum_{n = 0}^{M-1}{G}_n\right)\right] &= \mE\left[\exp\left( \beta q\sum_{n = 0}^{M-1}2\int_{t_{n}}^{t_{n+1}}\left( 3 + \frac{2}{\nu} +2C_B^2 + \frac{36C_e^2}{\nu}\|u(s)\|^2_{L^2} \right)\,ds\right)\right]\\\nonumber
		&= \exp\left(2\beta q\left(3+\frac{2}{\nu} + 2C_B^2\right)T\right) \mE\left[\exp\left( \frac{72C_e^2\beta q}{\nu}\int_{0}^T\|u(s)\|^2_{L^2}\, ds \right)\right],
	\end{align*}
	which, together with the inequality $\|u(t)\|_{L^2} \leq \frac{L^2}{4\pi^2}\|\partial_{x}^2 u(t)\|_{L^2}$, implies that
	\begin{align*}
		\mE\left[\exp\left( \beta q\sum_{n = 0}^{M-1}{G}_n\right)\right] 
		&\leq \exp\left(2\beta q\left(3+\frac{2}{\nu} + 2C_B^2\right)T\right)  \mE\left[\exp\left( \frac{9C_e^2L^4\beta q}{2\pi^4\nu}\int_{0}^T\|\partial_x^2u(s)\|^2_{L^2}\, ds \right)\right]\\\nonumber
		&\leq \exp\left(2\beta q\left(3+\frac{2}{\nu} + 2C_B^2\right)T\right) \mE\left[\exp\left\{ 8\kappa\|u_0\|^2_{L^2}+\kappa c_0T+\kappa \tilde{c}_0 \right\}\right]:= \widehat{C}_{0},
	\end{align*}
	where $\kappa = \frac{9C_e^2 L^4\beta q}{\nu^2\pi^4} < \frac{1}{32L_0^2}$ by using the hypothesis $L_0 <\frac{{\nu}\pi^2}{170 C_eL^2\sqrt{ q}}$.
	
	Next, we estimate $\Bigl(\mE\Bigl[\sup_{0 \leq \ell < M} F_{\ell}\Bigr]\Bigr)^{\frac12}$ as follows: firstly, we have
	\begin{align*}
		\mE\Bigl[\sup_{0 \leq \ell < M} F_{\ell}\Bigr] &\leq  \sum_{n=0}^{M-1}\mE\left[4kC_p^2h^{2(r-1)}\left(\nu + \frac{h^4}{\nu}\right)\|u(t_{n+1})\|^2_{H^{r+1}} + 4kC_B^2C_p^2 h^{2(r+1)}\|u(t_n)\|^2_{H^{r+1}}\right.\\\nonumber&\qquad\left.+ C_e^2kC_p^2\left(h^{2r}\|\p_x u(t_{n+1})\|^2_{L^2} + h^{2(r-1)}\|u_h^{n+1}\|^2_{L^2}\right)\|u(t_{n+1})\|^2_{H^{r+1}}\right]\\\nonumber
		&\qquad + \sum_{n = 0}^{M-1}\mE\left[4\nu\int_{t_n}^{t_{n+1}}\|\p_x^2(u(t_{n+1}) - u(s))\|^2_{L^2}\, ds\right. \\\nonumber
		&\qquad\left.+ \frac{1}{2} \int_{t_n}^{t_{n+1}}\|\p_x^2(u(t_{n+1}) - u(s))\|^2_{L^2}\, ds\right.\\\nonumber
		&\qquad\left. + C_e^2\int_{t_n}^{t_{n+1}}\left(\|\p_xu(s)\|^2_{L^2} + \|\p_x u(t_{n+1})\|^2_{L^2}\right)\|\p_xu(t_{n+1}) - u(s)\|^2_{L^2}\, ds\right] \\\nonumber
		&\qquad+\sum_{n=0}^{M-1}\left[ 4\left\|\int_{t_n}^{t_{n+1}} (B(u(s)) - B(u(t_{n})))\, dW(s)\right\|^2_{L^2}\right] \\\nonumber
		&\qquad+ \left[\sum_{n=0}^{M-1}\int_{t_{n}}^{t_{n+1}}\left( 3 + \frac{2}{\nu} + \frac{36C_e^2}{\nu}\|u(t_{n+1}) - u(s)\|^2_{L^2} \right)\, ds\|\varepsilon^{n+1}\|^2_{L^2}\right.\\\nonumber
&\qquad		\left.+ \int_{t_{\ell}}^{t_{\ell+1}}\left( 3 + \frac{2}{\nu} + \frac{36C_e^2}{\nu}\| u(s)\|^2_{L^2} \right)\, ds\|\varepsilon^{\ell+1}\|^2_{L^2}\right]\\\nonumber
		&:= T_1 + T_2 + T_3 + T_4.
	\end{align*}
	
	In view of Lemma \ref{lemma_poly_moment_H^s} and Lemma \ref{lemma_highmoment_discrete}, we obtain
	\begin{align*}
		T_1 &\leq kC_p^2h^{2(r-2)}\sum_{n=0}^{M-1}\left\{4\left(\nu + \frac{1}{\nu}\right)\mE\left[\|u(t_{n+1})\|^2_{H^{r}}\right] + 4C_B^2\mE\left[\|u(t_n)\|^2_{H^{r}}\right]\right.\\\nonumber&\qquad\left.+ C_e^2\left\{\mE\left[\left(\|\p_x u(t_{n+1})\|^2_{L^2} + \|u_h^{n+1}\|^2_{L^2}\right)^2\right]\right\}^{\frac12}\, \left\{\mE\left[\|u(t_{n+1})\|^4_{H^{r}}\right]\right\}^{\frac12}\right\}\\\nonumber
		&\leq C_p^2T\left\{4C_{r,2}\left(\nu + \frac{1}{\nu}\right) + 4C_B^2C_{r,2} + C_e^2\left(C_{1,4} + C_1\right)\right\}h^{2(r-2)}\\\nonumber
		&:= C_{T_1} h^{{2(r-2)}}.
	\end{align*}
	
	Using Lemma \ref{lemma_Holder} and Lemma \ref{lemma_poly_moment_H^s}, we get
	\begin{align*}
		T_2 &=\left(4\nu+\frac12\right)\sum_{n = 0}^{M-1}\int_{t_n}^{t_{n+1}}\mE\left[\|\p_x^2(u(t_{n+1}) - u(s))\|^2_{L^2}\right]\, ds\\\nonumber
		&\quad+ C_e^2\sum_{n = 0}^{M-1}\int_{t_n}^{t_{n+1}}\mE\left[\left(\|\p_xu(s)\|^2_{L^2} + \|\p_x u(t_{n+1})\|^2_{L^2}\right)\|\p_xu(t_{n+1}) - u(s)\|^2_{L^2}\right]\, ds\\\nonumber
		&\leq  \left(4\nu+\frac12\right)TK_{2,1}k \\\nonumber
		&\quad+ C_e^2\sum_{n = 0}^{M-1}\int_{t_n}^{t_{n+1}}\left\{\mE\left[\left(\|\p_xu(s)\|^2_{L^2} + \|\p_x u(t_{n+1})\|^2_{L^2}\right)^2\right]\right\}^{\frac12}\left\{\mE\left[\|\p_xu(t_{n+1}) - u(s)\|^4_{L^2}\right]\right\}^{\frac12}\, ds\\\nonumber
		&\leq \left(4\nu+\frac12\right)TK_{2,1}k + 4C_e^2TC_{1,4} K_{1,2} k = \left( \left(4\nu+\frac12\right)TK_{2,1} + 4C_e^2TC_{1,4} K_{1,2} \right)k := C_{T_2}k.
	\end{align*}
	
	To estimate $T_3$, we recall the It\^o isometry, the assumption \eqref{Assump_Lipschitz}, and Lemma \ref{lemma_Holder}:
	\begin{align*}
		T_3 &=4\sum_{n = 0}^{M-1} \mE\left[\int_{t_n}^{t_{n+1}} \|B(u(s)) - B(u(t_n))\|^2_{L^2}\, ds\right]\\\nonumber
		&\leq 4C_B\sum_{n = 0}^{M-1} \mE\left[\int_{t_n}^{t_{n+1}} \|u(s) - u(t_n)\|^2_{L^2}\, ds\right] \leq 4C_BT K_{0,1}k = C_{T_3}k.
	\end{align*}
	
	In light of Lemma \ref{lemma_poly_moment_H^s} and Lemma \ref{lemma_highmoment_discrete}, we also get
	\begin{align*}
		T_4 &=\mE\left[\sum_{n=0}^{M-1}\int_{t_{n}}^{t_{n+1}}\left( 3 + \frac{2}{\nu} + \frac{36C_e^2}{\nu}\|u(t_{n+1}) - u(s)\|^2_{L^2} \right)\, ds\|\varepsilon^{n+1}\|^2_{L^2}\right.\\\nonumber
		&\qquad		\left.+ \int_{t_{\ell}}^{t_{\ell+1}}\left( 3 + \frac{2}{\nu} + \frac{36C_e^2}{\nu}\| u(s)\|^2_{L^2} \right)\, ds\|\varepsilon^{\ell+1}\|^2_{L^2}\right]
		\leq C_{T_4}k.
	\end{align*}
	
	Collecting all the estimates from $T_1, ..., T_4$ into \eqref{eq3.13} we arrive at
	\begin{align*}
		&\Bigl(\mE\bigl[\sup_{0 \leq \ell \leq M}\|\varepsilon^{\ell}\|^{2q}_{L^2}\bigr]\Bigr)^{\frac{1}{2q}} + \left(\mE\left[\left(\frac{3\nu}{4}k\sum_{n = 0}^M\|\p_x^2\varepsilon^{n}\|^2_{L^2}\right)^q\right]\right)^{\frac{1}{2q}} \leq \widehat{C}_1\left(k + h^{2(r-2)}\right),
	\end{align*}
	where 
	\begin{align*}
		\widehat{C}_1 =\max\left\{C_{T_1}, C_{T_2}, C_{T_3}, C_{T_4}\right\}.
	\end{align*}
	
	Lastly, it remains to verify that $\{M_{\ell}\}_{\ell \geq 0}$ is a Martingale. To the end, we first use the It\^o's isometry,  the assumption \eqref{Assump_Lipschitz}, and the Burkholder-Davis-Gundy inequality to get
	\begin{align}\label{eq_3.16}
		\mE[|M_{\ell}|] &\leq \sum_{n = 0}^{\ell} \mE[|Z_n|]\\\nonumber
		&\leq 8k\sum_{n=0}^{\ell} \mE\bigl[\|B(u(t_n)) - B(u_h^n)\|^2_{L^2}\bigr] \\\nonumber
		&\qquad+ 2\mE\left[\left|\sum_{n = 0}^{\ell} \left(\int_{t_{n}}^{t_{n+1}}\bigl(B(u(s)) - B(u_h^n)\bigr)\, dW(s), \varepsilon^n\right)\right|\right]\\\nonumber
		&\leq 8C_BT(C_{0,2} + C_1) + \left(\mE\left[\sum_{n = 0}^{M-1}\int_{t_{n}}^{t_{n+1}} \|B(u(s)) - B(u_h^n)\|^2_{L^2}\|\varepsilon^n\|^2_{L^2}\, ds\right]\right)^{\frac12} \\\nonumber
		&\leq 8C_BT(C_{0,2} + C_1) + C_BT(C_{0,4} + C_2)(C_{0,4} + C_2) <\infty.
	\end{align}
	
	In addition, for any $0\leq n \leq M-1$, using the Martingale property of the It\^o integrals, we have
	\begin{align*}
		\mE[Z_n] &=  4\mE\bigl[\|B(u(t_n)) - B(u_h^n)\|^2_{L^2}|\Delta W_{n}|^2 \bigr]  { -  4k\mE\bigl[ \|B(u(t_n)) - B(u_h^n)\|^2_{L^2} \bigr] } \\\nonumber
		&\qquad \quad  + 2\mE\biggl[\biggl(\int_{t_n}^{t_{n+1}}\bigl(B(u(s))- B(u_h^n)\bigr)\,dW(s),\varepsilon^n\biggr)\biggr]\\\nonumber
		&= 4k \mE \bigl[\|B(u(t_n)) - B(u_h^n)\|^2_{L^2} \bigr] 
		{ - 4k\mE \bigl[ \|B(u(t_n)) - B(u_h^n)\|^2_{L^2} \bigr]  } + 0  =0.
	\end{align*}
	Then,  the conditional expectation of $M_{\ell}$ given $\{M_n\}_{n=0}^{\ell-1}$ is
	\begin{align}\label{eq3.17}
		\mE\bigl[M_{\ell} | M_0, M_1, \cdots, M_n \bigr]  
		&= \mE\bigl[Z_0+Z_1 + \cdots + Z_n| M_0, M_1, \cdots, M_n \bigr] \\\nonumber
		&\qquad+ \mE \bigl[Z_{n+1} + \cdots+Z_{\ell}| M_0, M_1,\cdots, M_{n} \bigr]\\\nonumber
		&= M_{n} + \mE\bigl[Z_{n+1} + \cdots + Z_{\ell} \bigr] = M_n.
	\end{align}
	
	Thus,  we conclude that $\bigl\{M_{\ell}; \ell \geq 0\bigr\}$ is a Martingale using \eqref{eq_3.16} and \eqref{eq3.17}.
	
	The proof is finished by combining the triangle inequality, \eqref{equu310}, and \eqref{projection_ineq}.

\end{proof}

\subsubsection{Higher moment error estimates}
The sub-second moment error estimates obtained in Theorem \ref{Theorem_sub_moment} imply a strong convergence in the $L_{\omega}^{p}L_t^{\infty}L_x^2$- and $L_{\omega}^{p}L_t^{2}H_x^2$-norm for $0<p<2$.  We note that these sub-second moment estimates are consequences of using the stochastic Gronwall inequality \eqref{ineq2.4}. However, we will demonstrate below that a bootstrap argument can overcome such a limitation to obtain higher moment estimates. {In turn, they allow for establishing a strong convergence in the $L_{\omega}^{p}L_t^{\infty}L_x^2$-norm for $0<p <4$, which is the goal of this subsection}.

\begin{theorem}\label{Theorem_higher_moment} 
	Let $u$ be the variational solution to \eqref{Weak_formulation} and $\{u_h^{n}\}_{n=1}^M$ be generated by \eqref{Scheme_Standard}. Let $u_0 \in \cap_{i=0}^{r} L^{2^{r+4-i}}(\Omega; H^i(D))$ for any integers $r \geq 4$. Assume that $B$ satisfies conditions \eqref{Assump_Lineargrowth} and \eqref{Assump_Lipschitz}. 
	Additionally, for any $0 < q < \frac{99}{100}$, assume that $L_0 <\frac{{\nu}\pi^2}{340 C_eL^2\sqrt{ q}}$ and $\mE\left[\exp\left(16\kappa \|u_0\|^2_{L^2}\right)\right] < \infty$,  {where $\kappa =\frac{3600 C_e^2L^4q}{\nu^2\pi^4}$, and $L_0$ and $C_e$ are respectively the constants as in \eqref{Assump_Lineargrowth} and \eqref{sobolev_ineq}}. Then, there holds
	\begin{align}\label{eq3.20}
		\left(\mE\left[\max_{1\leq n\leq M}\|u(t_n) - u_h^n\|^{4q}_{L^2}\right]\right)^{\frac{1}{4q}} 
		\leq \widehat{C}_2\left(k^{\frac12} + h^{r-2}\right),
	\end{align}
	for some constant  $\widehat{C}_2 = C(q,u_0,T,C_B)>0$.
\end{theorem}
\begin{proof} For the sake of simplicity, we only give the proof for $q = \frac12$. The proof for other $ q$'s should be similar and straightforward. First, testing \eqref{eq_3.14} with $\|\varepsilon^{n+1}\|^2_{L^2}$ and using the identity $2a(a-b) = a^2 - b^2 + (a-b)^2$ we obtain
	\begin{align}\label{eq3.21}
		&\frac{1}{4}\left[\|\varepsilon^{n+1}\|^4_{L^2} - \|\varepsilon^n\|^4_{L^2}\right] +\frac14 \left(\|\varepsilon^{n+1}\|^2_{L^2} - \|\varepsilon^n\|^2_{L^2}\right)^2 + \frac{3\nu k}{8}\|\p_x^2\varepsilon^{n+1}\|^2_{L^2}\|\varepsilon^{n+1}\|^2_{L^2} \\\nonumber
		&\leq 2kC_p^2h^{2(r-2)}\left[\left(\nu + \frac{1}{\nu}\right)\|u(t_{n+1})\|^2_{H^{r}} + C_B^2 \|u(t_n)\|^2_{H^{r}}\right]\|\varepsilon^{n+1}\|^2_{L^2}\\\nonumber&\qquad+ \frac{C_e^2k}{2}C_p^2h^{2(r-2)}\left[\|\p_x u(t_{n+1})\|^2_{L^2} + \|u_h^{n+1}\|^2_{L^2}\right]\|u(t_{n+1})\|^2_{H^{r}}\|\varepsilon^{n+1}\|^2_{L^2}\\\nonumber
		&\qquad+ \left(2\nu + \frac14\right)\int_{t_n}^{t_{n+1}}\|\p_x^2(u(t_{n+1}) - u(s))\|^2_{L^2}\|\varepsilon^{n+1}\|^2_{L^2}\, ds \\\nonumber
		&\qquad + \frac{C_e^2}{2}\int_{t_n}^{t_{n+1}}\left(\|\p_xu(s)\|^2_{L^2} + \|\p_x u(t_{n+1})\|^2_{L^2}\right)\|\p_xu(t_{n+1}) - u(s)\|^2_{L^2}\|\varepsilon^{n+1}\|^2_{L^2}\, ds \\\nonumber
		&\qquad+ 2\left\|\int_{t_n}^{t_{n+1}} (B(u(s)) - B(u(t_{n})))\, dW(s)\right\|^2_{L^2}\|\varepsilon^{n+1}\|^2_{L^2} \\\nonumber
		&\qquad + k\left( 3 + \frac{2}{\nu} + \frac{18C_e^2}{\nu}\|u(t_{n+1})\|^2_{L^2} \right)\|\varepsilon^{n+1}\|^4_{L^2} + 2kC_B^2\|\varepsilon^{n}\|^2_{L^2}\|\varepsilon^{n+1}\|^2_{L^2}\\\nonumber
		&\qquad +2\left\|(B(u(t_n)) - B(u_h^n))\Delta W_n\right\|^2_{L^2}\|\varepsilon^{n+1}\|^2_{L^2} - 2k\|B(u(t_n)) - B(u_h^n)\|^2_{L^2}\|\varepsilon^{n+1}\|^2_{L^2} \\\nonumber
		&\qquad+  \left(\int_{t_n}^{t_{n+1}} (B(u(s)) - B(u_h^n))\, dW(s),  \varepsilon^n\right)\|\varepsilon^{n+1}\|^2_{L^2}\\\nonumber
		&:= X_1 + ... + X_{10}.
	\end{align}
	
	Using Young's inequality, we obtain
	\begin{align*}
		&X_1 + X_2 + X_3 + X_4 + X_6+ X_7 + X_9\\\nonumber
		&= 2kC_p^2h^{2(r-2)}\left[\left(\nu + \frac{1}{\nu}\right)\|u(t_{n+1})\|^2_{H^{r}} + C_B^2 \|u(t_n)\|^2_{H^{r}}\right]\|\varepsilon^{n+1}\|^2_{L^2}\\\nonumber&\qquad+ \frac{C_e^2k}{2}C_p^2h^{2(r-2)}\left[\|\p_x u(t_{n+1})\|^2_{L^2} + \|u_h^{n+1}\|^2_{L^2}\right]\|u(t_{n+1})\|^2_{H^{r}}\|\varepsilon^{n+1}\|^2_{L^2}\\\nonumber
		&\qquad+ \left(2\nu + \frac14\right)\int_{t_n}^{t_{n+1}}\|\p_x^2(u(t_{n+1}) - u(s))\|^2_{L^2}\|\varepsilon^{n+1}\|^2_{L^2}\, ds \\\nonumber
		&\qquad + \frac{C_e^2}{2}\int_{t_n}^{t_{n+1}}\left(\|\p_xu(s)\|^2_{L^2} + \|\p_x u(t_{n+1})\|^2_{L^2}\right)\|\p_xu(t_{n+1}) - u(s)\|^2_{L^2}\|\varepsilon^{n+1}\|^2_{L^2}\, ds \\\nonumber
		&\qquad+ 2kC_B^2\|\varepsilon^{n}\|^2_{L^2}\|\varepsilon^{n+1}\|^2_{L^2}
		- 2k\|B(u(t_n)) - B(u_h^n)\|^2_{L^2}\|\varepsilon^{n+1}\|^2_{L^2} \\\nonumber
		&\qquad + k\left( 3 + \frac{2}{\nu} + \frac{18C_e^2}{\nu}\|u(t_{n+1})\|^2_{L^2} \right)\|\varepsilon^{n+1}\|^4_{L^2} \\\nonumber
		&\leq 6k\|\varepsilon^{n+1}\|^4_{L^2} + 2kC_p^4 h^{4(r-2)}\left[\left(\nu + \frac{1}{\nu}\right)^2\|u(t_{n+1})\|^4_{H^{r}} + C_B^4 \|u(t_n)\|^4_{H^{r}}\right]\\\nonumber
		&\qquad + \frac{C_e^4 C_p^4}{8}h^{4(r-2)} k\left[\|\p_x u(t_{n+1})\|^4_{L^2} + \|u_h^{n+1}\|^4_{L^2}\right]\|u(t_{n+1})\|^4_{H^{r}}\\\nonumber
		&\qquad + \frac14\left(2\nu + \frac14\right)^2\int_{t_n}^{t_{n+1}}\|\p_x^2(u(t_{n+1}) - u(s))\|^4_{L^2}\, ds\\\nonumber
		&\qquad + \frac{C_e^4}{8}\int_{t_n}^{t_{n+1}}\left(\|\p_xu(s)\|^4_{L^2} + \|\p_x u(t_{n+1})\|^4_{L^2}\right)\|\p_xu(t_{n+1}) - u(s)\|^4_{L^2}\, ds  + kC_B^4\|\varepsilon^n\|^4_{L^2}\\\nonumber
		&\qquad+ k\left( 3 + \frac{2}{\nu} + \frac{18C_e^2}{\nu}\|u(t_{n+1})\|^2_{L^2} \right)\|\varepsilon^{n+1}\|^4_{L^2} \\\nonumber
		&= 2kC_p^4 h^{4(r-2)}\left[\left(\nu + \frac{1}{\nu}\right)^2\|u(t_{n+1})\|^4_{H^{r}} + C_B^4 \|u(t_n)\|^4_{H^{r}}\right]\\\nonumber
		&\qquad + \frac{C_e^4 C_p^4}{8}h^{4(r-2)} k\left[\|\p_x u(t_{n+1})\|^4_{L^2} + \|u_h^{n+1}\|^4_{L^2}\right]\|u(t_{n+1})\|^4_{H^{r}}\\\nonumber
		&\qquad + \frac14\left(2\nu + \frac14\right)^2\int_{t_n}^{t_{n+1}}\|\p_x^2(u(t_{n+1}) - u(s))\|^4_{L^2}\, ds\\\nonumber
		&\qquad + \frac{C_e^4}{8}\int_{t_n}^{t_{n+1}}\left(\|\p_xu(s)\|^4_{L^2} + \|\p_x u(t_{n+1})\|^4_{L^2}\right)\|\p_xu(t_{n+1}) - u(s)\|^4_{L^2}\, ds  + kC_B^4\|\varepsilon^n\|^4_{L^2}\\\nonumber
		&\qquad+ k\left( 9 + \frac{2}{\nu} + \frac{18C_e^2}{\nu}\|u(t_{n+1})\|^2_{L^2} \right)\|\varepsilon^{n+1}\|^4_{L^2}.
	\end{align*}
	Next, to control $X_5, X_8$, and $X_{10}$, we proceed as follows:
	\begin{align*}
		X_5 + X_8 + X_{10}&= 2\left\|\int_{t_n}^{t_{n+1}} (B(u(s)) - B(u(t_{n})))\, dW(s)\right\|^2_{L^2}\|\varepsilon^{n+1}\|^2_{L^2}\\\nonumber
		&\qquad +2\left\|(B(u(t_n)) - B(u_h^n))\Delta W_n\right\|^2_{L^2}\|\varepsilon^{n+1}\|^2_{L^2}\\\nonumber
		&\qquad+  \left(\int_{t_n}^{t_{n+1}} (B(u(s)) - B(u_h^n))\, dW(s),  \varepsilon^n\right)\|\varepsilon^{n+1}\|^2_{L^2}\\\nonumber
		&= 2\left\|\int_{t_n}^{t_{n+1}} (B(u(s)) - B(u(t_{n})))\, dW(s)\right\|^2_{L^2}\left(\|\varepsilon^{n+1}\|^2_{L^2} - \|\varepsilon^{n}\|^2_{L^2}\right)\\\nonumber
		&\qquad +2\left\|(B(u(t_n)) - B(u_h^n))\Delta W_n\right\|^2_{L^2}\left(\|\varepsilon^{n+1}\|^2_{L^2} - \|\varepsilon^{n}\|^2_{L^2}\right)\\\nonumber
		&\qquad+  \left(\int_{t_n}^{t_{n+1}} (B(u(s)) - B(u_h^n))\, dW(s),  \varepsilon^n\right)\left(\|\varepsilon^{n+1}\|^2_{L^2} - \|\varepsilon^{n}\|^2_{L^2}\right)\\\nonumber
		&\qquad + 2\left\|\int_{t_n}^{t_{n+1}} (B(u(s)) - B(u(t_{n})))\, dW(s)\right\|^2_{L^2}\|\varepsilon^{n}\|^2_{L^2}\\\nonumber
		&\qquad +2\left\|(B(u(t_n)) - B(u_h^n))\Delta W_n\right\|^2_{L^2}\|\varepsilon^{n}\|^2_{L^2}\\\nonumber
		&\qquad+  \left(\int_{t_n}^{t_{n+1}} (B(u(s)) - B(u_h^n))\, dW(s),  \varepsilon^n\right)\|\varepsilon^{n}\|^2_{L^2}.
	\end{align*}
	Using Young's inequality on the first three terms on the right-hand side of $X_5 + X_8 + X_{10}$, we arrive at
	\begin{align*}
		X_5 + X_8 + X_{10} 
		&\leq 16\left\|\int_{t_n}^{t_{n+1}} (B(u(s)) - B(u(t_{n})))\, dW(s)\right\|^4_{L^2} +16\left\|(B(u(t_n)) - B(u_h^n))\Delta W_n\right\|^4_{L^2}\\\nonumber
		&\qquad+  4\left\|\int_{t_n}^{t_{n+1}} (B(u(s)) - B(u_h^n))\, dW(s)\right\|^2_{L^2} \|\varepsilon^n\|^2_{L^2}\\\nonumber
		&\qquad +2\left\|(B(u(t_n)) - B(u_h^n))\Delta W_n\right\|^2_{L^2}\|\varepsilon^{n}\|^2_{L^2}\\\nonumber
		&\qquad+  \left(\int_{t_n}^{t_{n+1}} (B(u(s)) - B(u_h^n))\, dW(s),  \varepsilon^n\right)\|\varepsilon^{n}\|^2_{L^2}  + \frac{3}{16}\left(\|\varepsilon^{n+1}\|^2_{L^2} - \|\varepsilon^n\|^2_{L^2}\right)^{2},
	\end{align*}
	which implies that
	\begin{align*}
		X_5 + X_8 + X_{10} 
		&\leq 16\left\|\int_{t_n}^{t_{n+1}} (B(u(s)) - B(u(t_{n})))\, dW(s)\right\|^4_{L^2} +16\left\|(B(u(t_n)) - B(u_h^n))\Delta W_n\right\|^4_{L^2}\\\nonumber
		&\qquad - 48k^2\left\|B(u(t_n)) - B(u_h^n)\right\|^4_{L^2}\\\nonumber
		&\qquad+  4\left\|\int_{t_n}^{t_{n+1}} (B(u(s)) - B(u_h^n))\, dW(s)\right\|^2_{L^2} \|\varepsilon^n\|^2_{L^2} \\\nonumber
		&\qquad\qquad- 4\int_{t_{n}}^{t_{n+1}}\|B(u(s)) - B(u^n_h)\|^2_{L^2}\|\varepsilon^n\|^2_{L^2}\, ds\\\nonumber
		&\qquad +2\left\|(B(u(t_n)) - B(u_h^n))\Delta W_n\right\|^2_{L^2}\|\varepsilon^{n}\|^2_{L^2} - 2k\|B(u(t_n)) - B(u_h^n)\|^2_{L^2}\|\varepsilon^n\|^2_{L^2}\\\nonumber
		&\qquad+  \left(\int_{t_n}^{t_{n+1}} (B(u(s)) - B(u_h^n))\, dW(s),  \varepsilon^n\right)\|\varepsilon^{n}\|^2_{L^2}\\\nonumber
		&\qquad + 4\int_{t_{n}}^{t_{n+1}}\|B(u(s)) - B(u^n_h)\|^2_{L^2}\|\varepsilon^n\|^2_{L^2}\, ds + 2k\|B(u(t_n)) - B(u_h^n)\|^2_{L^2}\|\varepsilon^n\|^2_{L^2}\\\nonumber
		&\qquad +48k^2\left\|B(u(t_n)) - B(u_h^n)\right\|^4_{L^2}+ \frac{3}{16}\left(\|\varepsilon^{n+1}\|^2_{L^2} - \|\varepsilon^n\|^2_{L^2}\right)^{2}\\\nonumber
		&\leq 16\left\|\int_{t_n}^{t_{n+1}} (B(u(s)) - B(u(t_{n})))\, dW(s)\right\|^4_{L^2} \\\nonumber
		&\qquad+ \left[16\left\|(B(u(t_n)) - B(u_h^n))\Delta W_n\right\|^4_{L^2} - 48k^2\|B(u(t_n)) - B(u^n_h)\|^4_{L^2}\right]\\\nonumber
		&\qquad+  \left[4\left\|\int_{t_n}^{t_{n+1}} (B(u(s)) - B(u_h^n))\, dW(s)\right\|^2_{L^2} \|\varepsilon^n\|^2_{L^2} \right.\\\nonumber
		&\qquad\qquad\left.- 4\int_{t_{n}}^{t_{n+1}}\|B(u(s)) - B(u^n_h)\|^2_{L^2}\|\varepsilon^n\|^2_{L^2}\, ds\right]\\\nonumber
		&\qquad +\left[2\left\|(B(u(t_n)) - B(u_h^n))\Delta W_n\right\|^2_{L^2}\|\varepsilon^{n}\|^2_{L^2} - 2k\|B(u(t_n)) - B(u_h^n)\|^2_{L^2}\|\varepsilon^n\|^2_{L^2}\right]\\\nonumber
		&\qquad+  \left(\int_{t_n}^{t_{n+1}} (B(u(s)) - B(u_h^n))\, dW(s),  \varepsilon^n\right)\|\varepsilon^{n}\|^2_{L^2}\\\nonumber
		&\qquad + 16 C_B^2\int_{t_{n}}^{t_{n+1}}\|u(s) - u(t_n)\|^4_{L^2}\, ds + 20C_B^2k\|\varepsilon^n\|^4_{L^2} + 384 C_B^4C_p^4 k^2 h^{4r}\|u(t_{n})\|^4_{H^{r}} \\\nonumber
		&\qquad+ 384 C_B^4k^2 \|\varepsilon^{n}\|^4_{L^2} + \frac{3}{16}\left(\|\varepsilon^{n+1}\|^2_{L^2} - \|\varepsilon^n\|^2_{L^2}\right)^{2},
	\end{align*}
	where the first inequality above is obtained by adding the terms $\pm 48k^2\|B(u(t_n)) - B(u_h^n)\|^4_{L^2}\pm4\int_{t_{n}}^{t_{n+1}}\|B(u(s)) - B(u^n_h)\|^2_{L^2}\|\varepsilon^n\|^2_{L^2}\, ds \pm 2k\|B(u(t_n)) - B(u_h^n)\|^2_{L^2}\|\varepsilon^n\|^2_{L^2}$. It should be noted that the modification will guarantee the use of the stochastic Gronwall inequality later.
	
	Now, collecting all the estimates from $X_1, ..., X_{10}$ and substituting them to \eqref{eq3.21} and absorbing the like terms to the left-hand side of \eqref{eq3.21}, and then applying the summation $\sum_{n = 0}^{\ell}$ for any $0\leq \ell \leq M-1$, we obtain the following inequality:
	\begin{align}\label{eq3.22}
		&	\frac{1}{4}\|\varepsilon^{\ell+1}\|^4_{L^2} + \frac{1}{16}\sum_{n = 0}^{\ell}\left(\|\varepsilon^{n+1}\|^2_{L^2} - \|\varepsilon^n\|^2_{L^2}\right)^2 +  \frac{3\nu k}{8}\sum_{n = 0}^{\ell}\|\p_x^2\varepsilon^{n+1}\|^2_{L^2}\|\varepsilon^{n+1}\|^2_{L^2} \\\nonumber
		&\leq  2kC_p^4 h^{4(r-2)}\sum_{n = 0}^{\ell}\left[\left(\nu + \frac{1}{\nu}\right)^2\|u(t_{n+1})\|^4_{H^{r}} + C_B^4 \|u(t_n)\|^4_{H^{r}}\right]\\\nonumber
		&\qquad + \frac{C_e^4 C_p^4}{8}h^{4(r-2)} k\sum_{n = 0}^{\ell}\left[\|\p_x u(t_{n+1})\|^4_{L^2} + \|u_h^{n+1}\|^4_{L^2}\right]\|u(t_{n+1})\|^4_{H^{r}}\\\nonumber
		&\qquad + \frac14\left(2\nu + \frac14\right)^2\sum_{n = 0}^{\ell}\int_{t_n}^{t_{n+1}}\|\p_x^2(u(t_{n+1}) - u(s))\|^4_{L^2}\, ds \\\nonumber
		&\qquad + \frac{C_e^4}{8}\sum_{n = 0}^{\ell}\int_{t_n}^{t_{n+1}}\left(\|\p_xu(s)\|^4_{L^2} + \|\p_x u(t_{n+1})\|^4_{L^2}\right)\|\p_xu(t_{n+1}) - u(s)\|^4_{L^2}\, ds\\\nonumber
		&\qquad+ 16\sum_{n = 0}^{\ell}\left\|\int_{t_n}^{t_{n+1}} (B(u(s)) - B(u(t_{n})))\, dW(s)\right\|^4_{L^2} \\\nonumber
		&\qquad+\sum_{n = 0}^{\ell} \left[16\left\|(B(u(t_n)) - B(u_h^n))\Delta W_n\right\|^4_{L^2} - 48k^2\|B(u(t_n)) - B(u^n_h)\|^4_{L^2}\right]\\\nonumber
		&\qquad+  \sum_{n = 0}^{\ell}\left[4\left\|\int_{t_n}^{t_{n+1}} (B(u(s)) - B(u_h^n))\, dW(s)\right\|^2_{L^2} \|\varepsilon^n\|^2_{L^2} \right.\\\nonumber
		&\qquad\qquad\left.- 4\int_{t_{n}}^{t_{n+1}}\|B(u(s)) - B(u^n_h)\|^2_{L^2}\|\varepsilon^n\|^2_{L^2}\, ds\right]\\\nonumber
		&\qquad +\sum_{n = 0}^{\ell}\left[2\left\|(B(u(t_n)) - B(u_h^n))\Delta W_n\right\|^2_{L^2}\|\varepsilon^{n}\|^2_{L^2} - 2k\|B(u(t_n)) - B(u_h^n)\|^2_{L^2}\|\varepsilon^n\|^2_{L^2}\right]\\\nonumber
		&\qquad+ \sum_{n = 0}^{\ell} \left(\int_{t_n}^{t_{n+1}} (B(u(s)) - B(u_h^n))\, dW(s),  \varepsilon^n\right)\|\varepsilon^{n}\|^2_{L^2}  \\\nonumber
		&\qquad + 16 C_B^2\sum_{n = 0}^{\ell}\int_{t_{n}}^{t_{n+1}}\|u(s) - u(t_n)\|^4_{L^2}\, ds + 384 C_B^4C_p^4 k^2 h^{4r}\sum_{n = 0}^{\ell}\|u(t_{n})\|^4_{H^{r}}\\\nonumber
		&\qquad + C_B^2\left(20 + C_B^2+ 384C_B^2k\right)k\sum_{n = 0}^{\ell}\|\varepsilon^n\|^4_{L^2} \\\nonumber
		&\qquad+ k\sum_{n = 0}^{\ell}\left( 9 + \frac{2}{\nu} + \frac{18C_e^2}{\nu}\|u(t_{n+1})\|^2_{L^2} \right)\|\varepsilon^{n+1}\|^4_{L^2}.
	\end{align}
	Denoting $Q_{n+1}:= 9 + \frac{2}{\nu} + \frac{18C_e^2}{\nu}\|u(t_{n+1})\|^2_{L^2}$, the last term on the right-hand side of \eqref{eq3.22} can be further analyzed as follows:
	\begin{align}\label{eq_3.23}
		k\sum_{n = 0}^{\ell}Q_{n+1}\|\varepsilon^{n+1}\|^4_{L^2}	&=k\sum_{n = 0}^{\ell-1}Q_{n+1}\|\varepsilon^{n+1}\|^4_{L^2} + kQ_{\ell+1}\|\varepsilon^{\ell+1}\|^4_{L^2}\\\nonumber
		&= k\sum_{n = 0}^{\ell-1}Q_{n+1}\|\varepsilon^{n+1}\|^4_{L^2} + kQ_{\ell+1}\|\varepsilon^{\ell}\|^4_{L^2} \\\nonumber
		&\qquad+ kQ_{\ell+1}\left(\|\varepsilon^{\ell+1}\|^2_{L^2} + \|\varepsilon^{\ell}\|^2_{L^2}\right)\left(\|\varepsilon^{\ell+1}\|^2_{L^2} - \|\varepsilon^{\ell}\|^2_{L^2}\right)\\\nonumber
		&\leq k\sum_{n = 0}^{\ell}\left(Q_{n} + Q_{n+1}\right)\|\varepsilon^{n}\|^4_{L^2} \\\nonumber
		&\qquad+ 8k^2Q^2_{\ell+1}\left(\|\varepsilon^{\ell+1}\|^2_{L^2} + \|\varepsilon^{\ell}\|^2_{L^2}\right)^2 + \frac{1}{32}\left(\|\varepsilon^{\ell+1}\|^2_{L^2} - \|\varepsilon^{\ell}\|^2_{L^2}\right)^2\\\nonumber
		&\leq k\sum_{n = 0}^{\ell}\left(18 + \frac{4}{\nu} + \frac{18C_e^2}{\nu}\|u(t_{n+1})\|^2_{L^2} + \frac{18C_e^2}{\nu}\|u(t_{n})\|^2_{L^2}\right)\|\varepsilon^{n}\|^4_{L^2} \\\nonumber
		&\qquad+ 8k^2Q^2_{\ell+1}\left(\|\varepsilon^{\ell+1}\|^2_{L^2} + \|\varepsilon^{\ell}\|^2_{L^2}\right)^2 + \frac{1}{32}\left(\|\varepsilon^{\ell+1}\|^2_{L^2} - \|\varepsilon^{\ell}\|^2_{L^2}\right)^2,
	\end{align}
	where the second equality of \eqref{eq_3.23} is obtained by adding $\pm kQ_{\ell+1}\|\varepsilon^{\ell}\|^4_{L^2}$, while the third inequality of \eqref{eq_3.23} is established by using the Young inequality. Additionally, the last term on the right-hand side of \eqref{eq_3.23} will be subsumed into the second term on the left-hand side of \eqref{eq3.22}.

	Now, combining \eqref{eq_3.23} into \eqref{eq3.22}, we obtain the following form for applying the stochastic Gronwall inequality \eqref{Stochastic_Gronwall}:
	\begin{align}\label{eq3.23}
		&	\|\varepsilon^{\ell+1}\|^4_{L^2} + \frac{1}{4}\sum_{n = 0}^{\ell-1} \left(\|\varepsilon^{n+1}\|^2_{L^2} - \|\varepsilon^n\|^2_{L^2}\right)^2 +   \frac{3\nu k}{2}\sum_{n = 0}^{\ell}\|\p_x^2\varepsilon^{n+1}\|^2_{L^2}\|\varepsilon^{n+1}\|^2_{L^2} \\\nonumber
		&\leq \mathcal{O}_{\ell} + Y_{\ell} + \sum_{n = 0}^{\ell} \mathcal{D}_n \|\varepsilon^n\|^4_{L^2},
	\end{align}
	where,
	\begin{align*}
		\mathcal{O}_{\ell}&:= 8kC_p^4 h^{4(r-2)}\sum_{n = 0}^{\ell}\left[\left(\nu + \frac{1}{\nu}\right)^2\|u(t_{n+1})\|^4_{H^{r}} + C_B^4\left(1+ 192k h^8\right) \|u(t_n)\|^4_{H^{r}}\right]\\\nonumber
		&\qquad+ \frac{C_e^4 C_p^4}{2}h^{4(r-2)} k\sum_{n = 0}^{\ell}\left[\|\p_x u(t_{n+1})\|^4_{L^2} + \|u_h^{n+1}\|^4_{L^2}\right]\|u(t_{n+1})\|^4_{H^{r}}\\\nonumber
		&\qquad + \left(2\nu + \frac14\right)^2\sum_{n = 0}^{\ell}\int_{t_n}^{t_{n+1}}\|\p_x^2(u(t_{n+1}) - u(s))\|^4_{L^2}\, ds \\\nonumber
		&\qquad + \frac{C_e^4}{2}\sum_{n = 0}^{\ell}\int_{t_n}^{t_{n+1}}\left(\|\p_xu(s)\|^4_{L^2} + \|\p_x u(t_{n+1})\|^4_{L^2}\right)\|\p_xu(t_{n+1}) - u(s)\|^4_{L^2}\, ds\\\nonumber
		&\qquad+ 64\sum_{n = 0}^{\ell} \left[\left\|\int_{t_n}^{t_{n+1}} (B(u(s)) - B(u(t_{n})))\, dW(s)\right\|^4_{L^2} +  C_B^2\int_{t_{n}}^{t_{n+1}}\|u(s) - u(t_n)\|^4_{L^2}\, ds\right]\\\nonumber
		&\qquad+ 32k^2\left(9 + \frac{2}{\nu} + \frac{18C_e^2}{\nu}\|u(t_{\ell+1})\|^2_{L^2}\right)^2\left(\|\varepsilon^{\ell+1}\|^2_{L^2} + \|\varepsilon^{\ell}\|^2_{L^2}\right)^2 \\\nonumber
		&\qquad+ 8\sum_{n = 0}^{\ell}\int_{t_{n}}^{t_{n+1}}\left( \frac{18C_e^2}{\nu}\left(\|u(t_{n+1}) - u(s)\|^2_{L^2} + \|u(t_{n}) - u(s)\|^2_{L^2}\right)\right)\, ds ,\\\nonumber
		Y_{\ell} &:= \sum_{n = 0}^{\ell} \mathcal{Y}_n,\\\nonumber
		\mathcal{Y}_{n}&:= 4\left[16\left\|(B(u(t_n)) - B(u_h^n))\Delta W_n\right\|^4_{L^2} - 48k^2\|B(u(t_n)) - B(u^n_h)\|^4_{L^2}\right]\\\nonumber
		&\qquad+  4\left[4\left\|\int_{t_n}^{t_{n+1}} (B(u(s)) - B(u_h^n))\, dW(s)\right\|^2_{L^2} \|\varepsilon^n\|^2_{L^2} \right.\\\nonumber
		&\qquad\qquad\left.- 4\int_{t_{n}}^{t_{n+1}}\|B(u(s)) - B(u^n_h)\|^2_{L^2}\|\varepsilon^n\|^2_{L^2}\, ds\right]\\\nonumber
		&\qquad +4\left[2\left\|(B(u(t_n)) - B(u_h^n))\Delta W_n\right\|^2_{L^2}\|\varepsilon^{n}\|^2_{L^2} - 2k\|B(u(t_n)) - B(u_h^n)\|^2_{L^2}\|\varepsilon^n\|^2_{L^2}\right]\\\nonumber
		&\qquad+  4\left(\int_{t_n}^{t_{n+1}} (B(u(s)) - B(u_h^n))\, dW(s),  \varepsilon^n\right)\|\varepsilon^{n}\|^2_{L^2},\\\nonumber
		\mathcal{D}_n &= 8\int_{t_{n}}^{t_{n+1}}\left(9 +10 C_B^2+ 193C_B^4 +\frac{2}{\nu} + \frac{36C_e^2}{\nu}\|u(s)\|^2_{L^2}\right)\, ds.
	\end{align*}
	
	In the next step, we assume that $\{Y_{\ell}: \ell \geq 0\}$ is a Martingale, and using the stochastic Gronwall inequality \eqref{Stochastic_Gronwall} to \eqref{eq3.23} with the following parameters: $q= \frac12, \alpha = \frac{3}{2}, \beta = 3$, we get the following inequality
	\begin{align}\label{eq3.24}
		&\Bigl(\mE\bigl[\sup_{0 \leq \ell \leq M}\|\varepsilon^{\ell}\|^{2}_{L^2}\bigr]\Bigr)^{\frac{1}{2}}  \leq \sqrt[3]{5} \left(\mE\left[\exp\left( \frac{3}{2}\sum_{n = 0}^{M-1}\mathcal{D}_n\right)\right]\right)^{\frac{1}{6}} \,\Bigl(\mE\Bigl[\sup_{0 \leq \ell < M} \mathcal{O}_{\ell}\Bigr]\Bigr)^{\frac14}.
	\end{align}
	To obtain the desired estimate \eqref{eq3.20} with $q = \frac12$, we control the right-hand side of \eqref{eq3.24} as follows: firstly, we control the exponential term by using Lemma \ref{lemma_expo_moment_H} and the inequality $\|u(t)\|_{L^2} \leq \frac{L^2}{4\pi^2}\|\partial_{x}^2 u(t)\|_{L^2}$, 
	\begin{align*}
		\mE\left[\exp\left( \frac{3}{2}\sum_{n = 0}^{M-1}\mathcal{D}_n\right)\right] &= \exp\left(12\left(9 + 10C_B^2 + 193 C_B^4 + \frac{2}{\nu}\right)T \right)\mE\left[ \exp\left(\frac{432 C_e^2}{\nu}\int_{0}^T\|u(s)\|^2_{L^2}\, ds \right)\right]\\\nonumber
		&\leq \exp\left(12\left(9 + 10C_B^2 + 193 C_B^4 + \frac{2}{\nu}\right)T \right)\mE\left[ \exp\left(\frac{27 C_e^2L^4}{\nu\pi^4}\int_{0}^T\|\partial_{x}^2u(s)\|^2_{L^2}\, ds \right)\right]\\\nonumber
		&\leq \exp\left(24\left(9 + 20C_B^2 + 385 C_B^4 + \frac{2}{\nu}\right)T \right) \sqrt{2}T\,\E\left[\exp\left\{ 8\kappa\|u_0\|^2_{L^2}+\kappa c_0T+\kappa \tilde{c}_0 \right\}\right]\\\nonumber
		&:= \tilde{C}_0^{6},
	\end{align*}
	where $\kappa = \frac{54C_e^2L^4}{\nu^2\pi^4} < \frac{1}{32 L_0^2}$. Next, we estimate $\Bigl(\mE\Bigl[\sup_{0 \leq \ell < M} \mathcal{O}_{\ell}\Bigr]\Bigr)^{\frac14}$ as follows:
	\begin{align*}
		\mE\Bigl[\sup_{0 \leq \ell < M} \mathcal{O}_{\ell}\Bigr] &\leq Z_1 + Z_2+Z_3 + Z_4.
	\end{align*}
	
	Using Lemma \ref{lemma_poly_moment_H^s} and Lemma \ref{lemma_highmoment_discrete}, we get
	\begin{align*}
		Z_1 &= 8kC_p^4 h^{4(r-2)}\sum_{n = 0}^{M-1}\mE\left[\left(\nu + \frac{1}{\nu}\right)^2\|u(t_{n+1})\|^4_{H^{r}} + C_B^4\left(1+ 192k h^8\right) \|u(t_n)\|^4_{H^{r}}\right]\\\nonumber
		&\qquad+ \frac{C_e^4 C_p^4}{2}h^{4(r-2)} k\sum_{n = 0}^{M-1}\mE\left[\left(\|\p_x u(t_{n+1})\|^4_{L^2} + \|u_h^{n+1}\|^4_{L^2}\right)\|u(t_{n+1})\|^4_{H^{r}}\right]\\\nonumber
		&\leq CC_p^4\left(C_{r,4} + C_B^4 C_{r,4} + C_e^4(C_{1, 8} + C_{3})C_{r,8} \right)h^{4(r-2)} := C_{X_1} h^{4(r-2)}.
	\end{align*}
	
	Using Lemma \ref{lemma_Holder}, and Lemma \ref{lemma_poly_moment_H^s}, we obtain
	\begin{align*}
		Z_2 &=  \left(2\nu + \frac14\right)^2\sum_{n = 0}^{M-1}\int_{t_n}^{t_{n+1}}\mE\left[\|\p_x^2(u(t_{n+1}) - u(s))\|^4_{L^2}\right]\, ds \\\nonumber
		&\qquad + \frac{C_e^4}{2}\sum_{n = 0}^{M-1}\int_{t_n}^{t_{n+1}}\mE\left[\left(\|\p_xu(s)\|^4_{L^2} + \|\p_x u(t_{n+1})\|^4_{L^2}\right)\|\p_xu(t_{n+1}) - u(s)\|^4_{L^2}\right]\, ds\\\nonumber
		&\leq CT K_{2,2} k^2 + C_e^4 T C_{1,8} K_{1,4} k^2 :=C_{X_2} k^2.
	\end{align*}
	
	Using Burkholder-Davis-Gundy inequality and Lemma \ref{lemma_Holder}, we also get
	\begin{align*}
		Z_3 &= 64\sum_{n = 0}^{M-1} \mE\left[\left\|\int_{t_n}^{t_{n+1}} (B(u(s)) - B(u(t_{n})))\, dW(s)\right\|^4_{L^2} +  C_B^2\int_{t_{n}}^{t_{n+1}}\|u(s) - u(t_n)\|^4_{L^2}\, ds\right],\\\nonumber
		&\leq 64\sum_{n = 0}^{M-1} C_{BDG}\left(\mE\left[\int_{t_n}^{t_{n+1}} \left\|B(u(s)) - B(u(t_{n}))\right\|^2_{L^2}\, ds\right]\right)^{2}+  C_B^2TK_{0,2} k^2,\\\nonumber
		&\leq 64TC_{BDG}K_{0,1}k^3 +  C_B^2TK_{0,2} k^2:= C_{X_3} k^2.
	\end{align*}
	
	Lastly, using Lemma \ref{lemma_Holder}, Lemma \ref{lemma_poly_moment_H^s} and Lemma \ref{lemma_highmoment_discrete}, we have
	\begin{align*}
		Z_4 &\leq 32k^2\mE\left[\sup_{0 \leq \ell \leq M-1}\left(9 + \frac{2}{\nu} + \frac{18C_e^2}{\nu}\|u(t_{\ell+1})\|^2_{L^2}\right)^2\left(\|\varepsilon^{\ell+1}\|^2_{L^2} + \|\varepsilon^{\ell}\|^2_{L^2}\right)^2\right] \\\nonumber
		&\qquad+\mE\left[8\sum_{n = 0}^{\ell}\int_{t_{n}}^{t_{n+1}}\left( \frac{18C_e^2}{\nu}\left(\|u(t_{n+1}) - u(s)\|^2_{L^2} + \|u(t_{n}) - u(s)\|^2_{L^2}\right)\right)\, ds\right]\\\nonumber
		&\leq  C_{X_4} k^2. 
	\end{align*}

	Collecting all the estimates from $X_1, X_2, X_3$, and $X_4$, we obtain the desired estimate \eqref{eq3.20} with $q = \frac12$.
	
	Finally, it is left to verify that $\{Y_{\ell}\}$ is a Martingale. First, using  the Burkholder-Davis-Gundy inequality, the assumption \eqref{Assump_Lipschitz}, Lemma \ref{lemma_poly_moment_H^s}, and Lemma \ref{lemma_highmoment_discrete}, we can easily obtain
	\begin{align}\label{eq_3.26}
		\mE[|Y_{\ell}|]  <\infty.
	\end{align}
	
	In addition, for any $0\leq n \leq M-1$, using the fact that $\mE[(\Delta W_n)^4] = 3k^2$, the independence of the increments $\Delta W_n$, It\^o isometry, and the Martingale property of It\^o integrals,  we also have
	\begin{align*}
		\mE[\mathcal{Y}_n] &=  4\mE\left[16\left\|(B(u(t_n)) - B(u_h^n))\Delta W_n\right\|^4_{L^2} - 48k^2\|B(u(t_n)) - B(u^n_h)\|^4_{L^2}\right]\\\nonumber
		&\qquad+  4\mE\left[4\left\|\int_{t_n}^{t_{n+1}} (B(u(s)) - B(u_h^n))\, dW(s)\right\|^2_{L^2} \|\varepsilon^n\|^2_{L^2} \right.\\\nonumber
		&\qquad\qquad\left.- 4\int_{t_{n}}^{t_{n+1}}\|B(u(s)) - B(u^n_h)\|^2_{L^2}\|\varepsilon^n\|^2_{L^2}\, ds\right]\\\nonumber
		&\qquad +4\mE\left[2\left\|(B(u(t_n)) - B(u_h^n))\Delta W_n\right\|^2_{L^2}\|\varepsilon^{n}\|^2_{L^2} - 2k\|B(u(t_n)) - B(u_h^n)\|^2_{L^2}\|\varepsilon^n\|^2_{L^2}\right]\\\nonumber
		&\qquad+  4\mE\left[\left(\int_{t_n}^{t_{n+1}} (B(u(s)) - B(u_h^n))\, dW(s),  \varepsilon^n\right)\|\varepsilon^{n}\|^2_{L^2}\right]\\\nonumber
		&= 0 + 0 + 0 + 0 = 0.
	\end{align*}
	Then,  the conditional expectation of $Y_{\ell}$ given $\{Y_n\}_{n=0}^{\ell-1}$ is
	\begin{align}\label{eq3.27}
		\mE\bigl[Y_{\ell} | Y_0, Y_1, \cdots, Y_n \bigr]  
		&= \mE\bigl[\mathcal{Y}_0 + \mathcal{Y}_1 + ... + \mathcal{Y}_n| Y_0, Y_1, ... , Y_n \bigr] \\\nonumber
		&\qquad+ \mE \bigl[\mathcal{Y}_{n+1} + ... + \mathcal{Y}_{\ell}| Y_0, Y_1, ... , Y_{n} \bigr]\\\nonumber
		&= Y_{n} + \mE\bigl[\mathcal{Y}_{n+1} + ... + \mathcal{Y}_{\ell} \bigr] = Y_n + 0 = Y_n.
	\end{align}
	
	Thus, we conclude that $\bigl\{Y_{\ell}; \ell \geq 0\bigr\}$ is a Martingale using \eqref{eq_3.26} and \eqref{eq3.27}. 
	
	The proof is completed by combining the triangle inequality, \eqref{eq3.20} and \eqref{projection_ineq}.
	
\end{proof}

\bigskip

Next, we also derive the second moment error estimate in $H^2$ norm.

\begin{theorem}\label{Theorem_second_moment_H^2} 
	Let $u$ be the variational solution to \eqref{Weak_formulation} and $\{u^{n}\}_{n=1}^M$ be generated by \eqref{Scheme_Standard}. Let {$u_0 \in \cap_{i=0}^{r} L^{2^{r+4-i}}(\Omega; H^i(D))$ for any integers $r \geq 4$}. Under the assumptions of Theorem \ref{Theorem_higher_moment}, there holds
	\begin{align}\label{eq3.28}
		\left(\mE\left[\nu k\sum_{n = 1}^M\|\p_x^2(u(t_n) - u_h^n)\|^2_{L^2}\right]\right)^{\frac{1}{2}} 
		\leq \widehat{C}_2\left(k^{\frac12} + h^{r-2}\right),
	\end{align}
	for some constant  $\widehat{C}_2 = C(q,u_0,T,C_B)>0$.
\end{theorem}
\begin{proof}
	The proof is obtained directly from the error inequality \eqref{eq3.16}, and Theorem \ref{Theorem_higher_moment}.
\end{proof}

\subsection{Partial expectation error estimates in the case of general multiplicative noise}\label{sub-sec3.3}
In this subsection, we derive error estimates for the scheme \eqref{Scheme_Standard} without assuming that $B$ is bounded, {i.e., when condition \eqref{Assump_Lineargrowth} is not satisfied}. Since $B$ may be unbounded, the exponential estimate in Lemma \ref{lemma_expo_moment_H} may not be available. So, we cannot reuse the techniques of the error estimates in the previous sections to control the nonlinear term with general multiplicative noise. To overcome such difficulty, we will draw upon the approach from \cite{bessaih2022space,breit2021convergence, carelli2012rates} dealing with the same issue for the Navier-Stokes equation. More specifically, letting $\rho>0$ be given and be chosen later, we introduce the following sequence of subsets of the sample space
\begin{align}
	{A_{\rho, m}} := \left\{\omega \in \Omega; \, \sup_{t \leq t_m} \|u(t)\|^2_{L^2} \leq \rho \right\},
\end{align}
where $u$ is the strong solution from \eqref{Weak_formulation}. We observe that {${A}_{\rho,0} \supset {A}_{\rho,1} \supset ... \supset {A}_{\rho,\ell}$}. 

It should be noted that the probability convergence of $\{u_h^n\}$, produced by the following theorem, is weaker than the results in Section \ref{sub-sec3.3}.

\bigskip

\begin{theorem}\label{Theorem_partial_expectation_error} 
	Let $u$ be the variational solution to \eqref{Weak_formulation} and $\{u_h^{n}\}_{n=1}^M$ be generated by \eqref{Scheme_Standard}. Let $u_0 \in \cap_{i=0}^{r} L^{2^{r+2-i}}(\Omega; H^i(D))$ for any integers $r \geq 4$. Assume that $B$ satisfies the condition \eqref{Assump_Lipschitz}. Then, there holds {for all $\beta>0$ small enough}
	\begin{align}\label{eq3.30}
		&\left(\max_{1\leq n\leq M}\mE\left[{\mathbf{1}_{{A}_{\rho, M}}}\|u(t_n) - u_h^n\|^{2}_{L^2}\right]\right)^{\frac{1}{2}} \\\nonumber
		&\qquad\qquad\qquad+ \left(\mE\left[{\mathbf{1}_{{A}_{\rho, M}}}\nu k\sum_{n = 1}^M\|\p_x^2(u(t_n) - u_h^n)\|^2_{L^2}\right]\right)^{\frac{1}{2}}
		\leq \widehat{C}_3\,\left(k^{\frac12} + h^{r-2}\right)h^{-\frac{\beta}{2}},
	\end{align}
	where $\widehat{C}_3 = C(u_0, T)$ is a positive constant. 
	
	{Moreover, for all $\beta$ sufficiently small, we may pick $\rho = \frac{\nu}{36 C_e^2T} \ln (h^{-\beta})$ so that} $\mP({{{A}_{\rho, M}}}) \rightarrow 1$ as $h \rightarrow 0$.
\end{theorem}
\begin{proof}
	First, multiplying \eqref{eq_3.14} with the indicator function $\mathbf{1}_{A_{\rho, n}}$ and using the fact that $\mathbf{1}_{A_{\rho, n}} \leq 1$, we obtain
	\begin{align}\label{eq_3.31}
		&\frac{1}{2}\left[\mathbf{1}_{A_{\rho, n}}\|\varepsilon^{n+1}\|^2_{L^2} - \mathbf{1}_{A_{\rho, n}}\|\varepsilon^n\|^2_{L^2}\right]  +\frac14\mathbf{1}_{A_{\rho, n}}\|\varepsilon^{n+1} - \varepsilon^n\|^2_{L^2} + \frac{3\nu k}{8}\mathbf{1}_{A_{\rho, n}}\|\p_x^2\varepsilon^{n+1}\|^2_{L^2} \\\nonumber
		&\leq 2kC_p^2h^{2(r-2)}\left(\nu + \frac{1}{\nu}\right)\|u(t_{n+1})\|^2_{H^{r}} + 2kC_B^2C_p^2 h^{2r}\|u(t_n)\|^2_{H^{r}}\\\nonumber&\qquad+ \frac{C_e^2k}{2}C_p^2\left[h^{2(r-1)}\|u(t_{n+1})\|^2_{H^{r}}\|\p_x u(t_{n+1})\|^2_{L^2} + h^{2(r-2)}\|u_h^{n+1}\|^2_{L^2}\|u(t_{n+1})\|^2_{H^{r}}\right]\\\nonumber
		&\qquad+ 2\nu\int_{t_n}^{t_{n+1}}\|\p_x^2(u(t_{n+1}) - u(s))\|^2_{L^2}\, ds + \frac{1}{4} \int_{t_n}^{t_{n+1}}\|\p_x^2(u(t_{n+1}) - u(s))\|^2_{L^2}\, ds\\\nonumber
		&\qquad + \frac{C_e^2}{2}\int_{t_n}^{t_{n+1}}\left(\|\p_xu(s)\|^2_{L^2} + \|\p_x u(t_{n+1})\|^2_{L^2}\right)\|\p_xu(t_{n+1}) - u(s)\|^2_{L^2}\, ds \\\nonumber
		&\qquad+ 2\left\|\int_{t_n}^{t_{n+1}} (B(u(s)) - B(u(t_{n})))\, dW(s)\right\|^2_{L^2} \\\nonumber
		&\qquad + k\mathbf{1}_{A_{\rho, n}}\left( 3 + \frac{2}{\nu} + \frac{18C_e^2}{\nu}\|u(t_{n+1})\|^2_{L^2} \right)\|\varepsilon^{n+1}\|^2_{L^2} + 2kC_B^2\mathbf{1}_{A_{\rho, n}}\|\varepsilon^{n}\|^2_{L^2}\\\nonumber
		&\qquad +2\mathbf{1}_{A_{\rho, n}}\left\|(B(u(t_n)) - B(u_h^n))\Delta W_n\right\|^2_{L^2}  \\\nonumber
		&\qquad+  \mathbf{1}_{A_{\rho, n}}\left(\int_{t_n}^{t_{n+1}} (B(u(s)) - B(u^n))\, dW(s),  \varepsilon^n\right).
	\end{align}
	
	Using the fact that $\mathbf{1}_{A_{\rho, n+1}} \leq \mathbf{1}_{A_{\rho, n}}$ for all $n =0, ..., M-1$, taking the expectation, and then applying the summation $\sum_{n=0}^{\ell}$, we get
	\begin{align}\label{eq3.32}
		&\frac{1}{2}\mE\left[\mathbf{1}_{A_{\rho, \ell+1}}\|\varepsilon^{\ell+1}\|^2_{L^2}   +\frac14\sum_{n = 0}^{\ell}\mathbf{1}_{A_{\rho, n}}\|\varepsilon^{n+1} - \varepsilon^n\|^2_{L^2} + \frac{3\nu k}{8}\sum_{n = 0}^{\ell}\mathbf{1}_{A_{\rho, n+1}}\|\p_x^2\varepsilon^{n+1}\|^2_{L^2}\right] \\\nonumber
		&\leq 2kC_p^2h^{2(r-2)}\sum_{n = 0}^{\ell}\mE\left[\left(\nu + \frac{1}{\nu}\right)\|u(t_{n+1})\|^2_{H^{r}} + C_B^2 \|u(t_n)\|^2_{H^{r}}\right]\\\nonumber&\qquad+ \frac{C_e^2k}{2}C_p^2h^{2(r-2)}\sum_{n=0}^{\ell}\mE\left[\|u(t_{n+1})\|^2_{H^{r}}\|\p_x u(t_{n+1})\|^2_{L^2} + \|u_h^{n+1}\|^2_{L^2}\|u(t_{n+1})\|^2_{H^{r}}\right]\\\nonumber
		&\qquad+ \left(2\nu + \frac14\right)\sum_{n = 0}^{\ell}\int_{t_n}^{t_{n+1}}\mE\left[\|\p_x^2(u(t_{n+1}) - u(s))\|^2_{L^2}\right]\, ds\\\nonumber
		&\qquad + \frac{C_e^2}{2}\sum_{n=0}^{\ell}\int_{t_n}^{t_{n+1}}\mE\left[\left(\|\p_xu(s)\|^2_{L^2} + \|\p_x u(t_{n+1})\|^2_{L^2}\right)\|\p_xu(t_{n+1}) - u(s)\|^2_{L^2}\right]\, ds \\\nonumber
		&\qquad+ 2\sum_{n = 0}^{\ell}\mE\left[\left\|\int_{t_n}^{t_{n+1}} (B(u(s)) - B(u(t_{n})))\, dW(s)\right\|^2_{L^2}\right] + 2kC_B^2\sum_{n = 0}^{\ell}\mE\left[\mathbf{1}_{A_{\rho, n}}\|\varepsilon^{n}\|^2_{L^2}\right]\\\nonumber
		&\qquad + k\sum_{n = 0}^{\ell}\mE\left[\mathbf{1}_{A_{\rho, n}}\left( 3 + \frac{2}{\nu} + \frac{18C_e^2}{\nu}\|u(t_{n+1})\|^2_{L^2} \right)\|\varepsilon^{n+1}\|^2_{L^2}\right] \\\nonumber
		&\qquad +2\sum_{n = 0}^{\ell}\mE\left[\mathbf{1}_{A_{\rho, n}}\left\|(B(u(t_n)) - B(u_h^n))\Delta W_n\right\|^2_{L^2} \right]\\\nonumber
		&\qquad + \sum_{n = 0}^{\ell}\mE\left[ \mathbf{1}_{A_{\rho, n}}\left(\int_{t_n}^{t_{n+1}} (B(u(s)) - B(u^n))\, dW(s),  \varepsilon^n\right)\right]\\\nonumber
		&:= Q_1 + ... + Q_8.
	\end{align}
	
	First of all, $Q_8 =0$ thanks to the Martingale property of It\^o integrals. Next, using Lemma \ref{lemma_poly_moment_H^s} and Lemma \ref{lemma_highmoment_discrete}, we estimate $Q_1$ and $Q_2$ as below
	\begin{align*}
		Q_1 + Q_2 &\leq C_p^2T\left(C_{r,2} + C_B^2 C_{r,2}\right)  h^{2(r-2)} + C_e^2 C_p^2 C_{r,4}\left(C_{1,4} + C_2\right) h^{2(r-2)}:= C_{Q_1,Q_2} h^{2(r-2)}.
	\end{align*}
	
	Next, using Lemma \ref{lemma_Holder}, and Lemma \ref{lemma_poly_moment_H^s} we obtain
	\begin{align*}
		Q_3 + Q_4 &\leq \left(2\nu + \frac14\right)T K_{2,1} k + C_e^2 TC_{1,4} K_{1,2} k :=C_{Q_3,Q_4} k.
	\end{align*}
	
	Using the It\^o isometry and condition \eqref{Assump_Lipschitz}, we get
	\begin{align*}
		Q_5 &= 2\sum_{n = 0}^{\ell}\mE\left[\int_{t_n}^{t_{n+1}} \left\|B(u(s)) - B(u(t_{n}))\right\|^2_{L^2}\, ds\right] + 2kC_B^2\sum_{n = 0}^{\ell}\mE\left[\mathbf{1}_{A_{\rho, n}}\|\varepsilon^{n}\|^2_{L^2}\right]\\\nonumber
		&\leq 2C_B^2\sum_{n = 0}^{\ell}\mE\left[\int_{t_n}^{t_{n+1}} \left\|u(s) - u(t_{n})\right\|^2_{L^2}\, ds\right] + 2kC_B^2\sum_{n = 0}^{\ell}\mE\left[\mathbf{1}_{A_{\rho, n}}\|\varepsilon^{n}\|^2_{L^2}\right]\\\nonumber
		&\leq 2C_B^2T K_{0,1} k + 2kC_B^2\sum_{n = 0}^{\ell}\mE\left[\mathbf{1}_{A_{\rho, n}}\|\varepsilon^{n}\|^2_{L^2}\right]:= C_{Q_5} k + 2kC_B^2\sum_{n = 0}^{\ell}\mE\left[\mathbf{1}_{A_{\rho, n}}\|\varepsilon^{n}\|^2_{L^2}\right].
	\end{align*}
	
	Similarly, 
	\begin{align*}
		Q_7 &\leq 2C_B^2 k\sum_{n = 0}^{\ell} \mE\left[\mathbf{1}_{A_{\rho, n}}\|\varepsilon^n\|^2_{L^2}\right] + 2C_B^2 TC_p^2 C_{r,2} h^{2(r-2)}\\\nonumber
		&:=2C_B^2  k\sum_{n = 0}^{\ell} \mE\left[\mathbf{1}_{A_{\rho, n}}\|\varepsilon^n\|^2_{L^2}\right] + C_{Q_7}\, h^{2(r-2)}.
	\end{align*}
	
	Adding $\pm \|\varepsilon^{n}\|^2_{L^2}$, we estimate $Q_6$ as follows:
	\begin{align*}
		Q_6&=k\sum_{n = 0}^{\ell}\mE\left[\mathbf{1}_{A_{\rho, n}}\left( 3 + \frac{2}{\nu} + \frac{18C_e^2}{\nu}\|u(t_{n+1})\|^2_{L^2} \right)\left(\|\varepsilon^{n+1}\|^2_{L^2} - \|\varepsilon^{n}\|^2_{L^2}\right)\right] \\\nonumber
		&\qquad+  k\sum_{n = 0}^{\ell}\mE\left[\mathbf{1}_{A_{\rho, n}}\left( 3 + \frac{2}{\nu} + \frac{18C_e^2}{\nu}\|u(t_{n+1})\|^2_{L^2} \right) \|\varepsilon^{n}\|^2_{L^2}\right]\\\nonumber
		&=k\sum_{n = 0}^{\ell}\mE\left[\mathbf{1}_{A_{\rho, n}}\left( 3 + \frac{2}{\nu} + \frac{18C_e^2}{\nu}\|u(t_{n+1})\|^2_{L^2} \right)\left(\|\varepsilon^{n+1}\|_{L^2} + \|\varepsilon^{n}\|_{L^2}\right)\left(\|\varepsilon^{n+1}\|_{L^2} - \|\varepsilon^{n}\|_{L^2}\right)\right] \\\nonumber
		&\qquad+  k\sum_{n = 0}^{\ell}\mE\left[\mathbf{1}_{A_{\rho, n}}\left( 3 + \frac{2}{\nu} + \frac{18C_e^2}{\nu}\|u(t_{n+1})\|^2_{L^2} \right) \|\varepsilon^{n}\|^2_{L^2}\right].
	\end{align*}
	
	Next, using the triangle inequality, Young inequality, then Lemma \ref{lemma_poly_moment_H^s}, Lemma \ref{lemma_highmoment_discrete}, and Lemma \ref{lemma_Holder}, we obtain
	\begin{align*}
		Q_6 &\leq k\sum_{n = 0}^{\ell}\mE\left[\mathbf{1}_{A_{\rho, n}}\left( 3 + \frac{2}{\nu} + \frac{18C_e^2}{\nu}\|u(t_{n+1})\|^2_{L^2} \right)\left(\|\varepsilon^{n+1}\|_{L^2} + \|\varepsilon^{n}\|_{L^2}\right)\|\varepsilon^{n+1} - \varepsilon^{n}\|_{L^2}\right] \\\nonumber
		&\qquad+  k\sum_{n = 0}^{\ell}\mE\left[\mathbf{1}_{A_{\rho, n}}\left( 3 + \frac{2}{\nu} + \frac{36C_e^2}{\nu}\|u(t_{n+1}) - u(t_{n})\|^2_{L^2}  + \frac{36C_e^2}{\nu}\|u(t_n)\|^2_{L^2}\right) \|\varepsilon^{n}\|^2_{L^2}\right]\\\nonumber
		&\leq \frac{1}{8} \sum_{n = 0}^{\ell}\mE\left[\mathbf{1}_{A_{\rho, n}}\|\varepsilon^{n+1} - \varepsilon^n\|^2_{L^2}\right] \\\nonumber
		&\qquad+ 2k^2\sum_{n = 0}^{\ell}\mE\left[\mathbf{1}_{A_{\rho, n}}\left( 3 + \frac{2}{\nu} + \frac{18C_e^2}{\nu}\|u(t_{n+1})\|^2_{L^2} \right)^2\left(\|\varepsilon^{n+1}\|_{L^2} + \|\varepsilon^{n}\|_{L^2}\right)^2\right] \\\nonumber
		&\qquad+  k\sum_{n = 0}^{\ell}\mE\left[ \frac{36C_e^2}{\nu}\|u(t_{n+1}) - u(t_{n})\|^2_{L^2}  \|\varepsilon^{n}\|^2_{L^2}\right]\\\nonumber
		&\qquad+  k\sum_{n = 0}^{\ell}\mE\left[\mathbf{1}_{A_{\rho, n}}\left( 3 + \frac{2}{\nu} + \frac{36C_e^2}{\nu}\|u(t_n)\|^2_{L^2}\right) \|\varepsilon^{n}\|^2_{L^2}\right]\\\nonumber
		&\leq \frac{1}{8} \sum_{n = 0}^{\ell}\mE\left[\mathbf{1}_{A_{\rho, n}}\|\varepsilon^{n+1} - \varepsilon^n\|^2_{L^2}\right] +  T C_e^2C_{0,8}(C_{0,8} + C_{8}) k + C_e^2 T K_{0,4}(C_{0,4} + C_4) k\\\nonumber
		&\qquad+  \left( 3 + \frac{2}{\nu} + \frac{36C_e^2}{\nu}\rho\right) k\sum_{n = 0}^{\ell}\mE\left[\mathbf{1}_{A_{\rho, n}} \|\varepsilon^{n}\|^2_{L^2}\right]\\\nonumber
		&:= \frac{1}{8} \sum_{n = 0}^{\ell}\mE\left[\mathbf{1}_{A_{\rho, n}}\|\varepsilon^{n+1} - \varepsilon^n\|^2_{L^2}\right] +  C_{Q_6}k +  \left( 3 + \frac{2}{\nu} + \frac{36C_e^2}{\nu}\rho\right) k\sum_{n = 0}^{\ell}\mE\left[\mathbf{1}_{A_{\rho, n}} \|\varepsilon^{n}\|^2_{L^2}\right],
	\end{align*}
	where the first term $\frac{1}{8} \sum_{n = 0}^{\ell}\mE\left[\mathbf{1}_{A_{\rho, n}}\|\varepsilon^{n+1} - \varepsilon^n\|^2_{L^2}\right] $ will be absorbed to the left-hand side of \eqref{eq3.32}.
	
	Now, collecting all the estimates from $Q_1, ..., Q_8$ and substituting these into \eqref{eq3.32} we arrive at
	\begin{align}\label{eq3.33}
		&\frac{1}{2}\mE\left[\mathbf{1}_{A_{\rho, \ell+1}}\|\varepsilon^{\ell+1}\|^2_{L^2}   +\frac18\sum_{n = 0}^{\ell}\mathbf{1}_{A_{\rho, n}}\|\varepsilon^{n+1} - \varepsilon^n\|^2_{L^2} + \frac{3\nu k}{8}\sum_{n = 0}^{\ell}\mathbf{1}_{A_{\rho, n+1}}\|\p_x^2\varepsilon^{n+1}\|^2_{L^2}\right] \\\nonumber
		&\leq \left(C_{Q_1,Q_2} + C_{Q_7}\right)h^{2(r-2)} + \left(C_{Q_3,Q_4} + C_{Q_5} + C_{Q_6}\right) k\\\nonumber
		&\qquad +\left( 3 + 4C_B^2+ \frac{2}{\nu} + \frac{36C_e^2}{\nu}\rho\right) k\sum_{n = 0}^{\ell}\mE\left[\mathbf{1}_{A_{\rho, n}} \|\varepsilon^{n}\|^2_{L^2}\right].
	\end{align}
	
	Applying the discrete deterministic Gronwall inequality to \eqref{eq3.33}, we obtain
	\begin{align}\label{eq3.34}
		&\frac{1}{2}\mE\left[\mathbf{1}_{A_{\rho, \ell+1}}\|\varepsilon^{\ell+1}\|^2_{L^2}   +\frac18\sum_{n = 0}^{\ell}\mathbf{1}_{A_{\rho, n}}\|\varepsilon^{n+1} - \varepsilon^n\|^2_{L^2} + \frac{3\nu k}{8}\sum_{n = 0}^{\ell}\mathbf{1}_{A_{\rho, n+1}}\|\p_x^2\varepsilon^{n+1}\|^2_{L^2}\right] \\\nonumber
		&\leq \left\{\left(C_{Q_1,Q_2} + C_{Q_7}\right)h^{2(r-2)} + \left(C_{Q_3,Q_4} + C_{Q_5} + C_{Q_6}\right) k\right\}\exp\left(T\left( 3 + 4C_B^2+ \frac{2}{\nu} + \frac{36C_e^2}{\nu}\rho\right)\right)\\\nonumber
		&=  \left\{\left(C_{Q_1,Q_2} + C_{Q_7}\right)h^{2(r-2)} + \left(C_{Q_3,Q_4} + C_{Q_5} + C_{Q_6}\right) k\right\}h^{-\beta}\exp\left(T\left( 3 + 4C_B^2+ \frac{2}{\nu}\right)\right)\\\nonumber
		&:= \widehat{C}_3 \left(k + h^{2(r-2)}\right) h^{-\beta}.
	\end{align}
	
	Then, the desired estimate \eqref{eq3.30} is obtained by using the triangle inequality and combining with \eqref{eq3.34} and \eqref{projection_ineq}.
	
	Lastly, using the Markov inequality and Lemma \ref{lemma_poly_moment_H^s}, we also obtain
	\begin{align*}
		\mP\left(A_{\rho, n}\right) &= 1 - 	\mP\left(A^c_{\rho, n}\right) \\\nonumber
		&\geq 1 - \frac{\mE\left[\sup_{t \leq t_n}\|u(t)\|^2_{L^2}\right]}{C\ln (h^{-\beta})} \rightarrow 1 \qquad\text{as } h \rightarrow 0.
	\end{align*}
	
	This finishes the proof.
	
\end{proof}

{	\section{Numerical experiments} \label{sec-4}
	
	In this section, we present numerical experiments to validate the theoretical results established in the preceding sections. The computational domain is $D = (0,1) \subset \mathbb{R}$, with final time $T = 0.01$ and diffusion coefficient $\nu = 0.01$. The Wiener process $W(t)$ in \eqref{eq1.1} is assumed to be real-valued and is simulated using a minimal time step $k_0 = 2^{-16}$.
	
	We consider two representative choices for the diffusion coefficient $G$, namely $G(u) = \alpha_1 \sin(u)$ and $G(u) = \alpha_2 u$, with $\alpha_1, \alpha_2 > 0$. These correspond to bounded and unbounded multiplicative noise, respectively, and are consistent with the assumptions used in the error estimates of Theorem~\ref{Theorem_second_moment_H^2} and Theorem~\ref{Theorem_partial_expectation_error}.
	
	All expectations are approximated by a standard Monte Carlo method with $J = 400$ independent sample paths. The spatial discretization is carried out using $P_3$ finite elements. In all experiments, we aim to verify the global convergence rate of order $O(k^{\frac{1}{2}} + h^2)$ predicted by Theorem~\ref{Theorem_second_moment_H^2} and Theorem~\ref{Theorem_partial_expectation_error}, under the time step scaling $k = h^4$.
	
	We implement the Main Algorithm and measure the error of the numerical solution in the norm
	\begin{align*}
		L^2_{\omega}L^{2}_tH^2_x(u)
		:= \left(\mathbb{E}\left[k\sum_{n=1}^{M}\|u(t_n) - u_h^n\|^2_{H^2}\right]\right)^{1/2}.
	\end{align*}
	Since the exact solution is not available, the error is approximated by
	\begin{align*}
		\left(\frac{1}{J}\sum_{j=1}^{J}\left(k\sum_{n=1}^{M}\|u_{h_{\mathrm{ref}}}^n(\omega_j) - u_h^n(\omega_j)\|^2_{H^2}\right)\right)^{1/2},
	\end{align*}
	where $u_{h_{\mathrm{ref}}}^n(\omega_j)$ is a reference solution computed on the same sample path using a finer mesh with $h_{\mathrm{ref}} = h/2$.
	
	To solve the nonlinear system arising from the fully implicit scheme, we employ a fixed-point iteration (see, e.g., \cite[Algorithm~3]{vo2024high} and \cite[Section~4.1]{feng2017finite}), which is more robust and efficient than Newton-type methods in the stochastic setting.
	
	\medskip
	\noindent
	\textbf{Test 1.} We first examine the convergence behavior under periodic boundary conditions, which are assumed throughout the analysis. The initial condition is
	\begin{align*}
		u_0(x) = 0.2 \sin(2\pi x).
	\end{align*}
	The results in Tables~\ref{tab:1} and \ref{tab:2} confirm the theoretically predicted second-order convergence rate in space and half-order convergence rate in time.
	
	\begin{table}[htb]
		\renewcommand{\arraystretch}{1.2}
		\centering
		\begin{tabular}{|c|c|c|c|c|}
			\hline
			\hline
			$h$ & $k$ & $L_{\omega}^2L_t^2H^2_x(u)$ error & spatial rate & time rate \\
			\hline
			\hline
			$2^{-1}$ & $2^{-4}$  & 2.92176869e-02 & ---    & ---    \\
			\hline
			$2^{-2}$ & $2^{-8}$  & 5.76473492e-03 & 2.3415 & 0.5854 \\
			\hline
			$2^{-3}$ & $2^{-12}$ & 1.41012220e-03 & 2.0314 & 0.5078 \\
			\hline
			$2^{-4}$ & $2^{-16}$ & 3.52583309e-04 & 1.9998 & 0.4999 \\
			\hline
		\end{tabular}
		\caption{Errors, spatial convergence rates, and time convergence rates for $k=h^4$ with $B(u)=u$ under periodic boundary conditions.}
		\label{tab:1}
	\end{table}

	\begin{table}[htb]
		\renewcommand{\arraystretch}{1.2}
		\centering
		\begin{tabular}{|c|c|c|c|c|}
			\hline
			\hline
			$h$ & $k$ & $L_{\omega}^2L_t^2H^2_x(u)$ error & spatial rate & time rate \\
			\hline
			\hline
			$2^{-1}$ & $2^{-4}$  & 3.20175727e-02 & ---    & ---    \\
			\hline
			$2^{-2}$ & $2^{-8}$  & 8.61849229e-03 & 1.8934 & 0.4734 \\
			\hline
			$2^{-3}$ & $2^{-12}$ & 2.10420639e-03 & 2.0342 & 0.5086 \\
			\hline
			$2^{-4}$ & $2^{-16}$ & 5.01052324e-04 & 2.0702 & 0.5176 \\
			\hline
		\end{tabular}
		\caption{Errors, spatial convergence rates, and time convergence rates for $k=h^4$ with $B(u)=5\sin(u)$ under periodic boundary conditions.}
		\label{tab:2}
	\end{table}
	
	\medskip
	\noindent
	\textbf{Test 2.} We next consider Dirichlet boundary conditions
	\begin{align*}
		u(0,t) = u(1,t) = 0, \qquad \partial_x u(0,t) = \partial_x u(1,t) = 0,
	\end{align*}
	with initial condition $u_0(x) = x^2(1-x)^2$. Although this case is not covered by the theoretical analysis, the numerical results in Tables~\ref{tab:3} and \ref{tab:4} indicate that the method remains stable and achieves a second-order convergence rate in space and a half-order convergence rate in time.
	
	\begin{table}[htb]
		\renewcommand{\arraystretch}{1.2}
		\centering
		\begin{tabular}{|c|c|c|c|c|}
			\hline
			\hline
			$h$ & $k$ & $L_{\omega}^2L_t^2H^2_x(u)$ error & spatial rate & time rate \\
			\hline
			\hline
			$2^{-1}$ & $2^{-4}$  & 1.71777696e-03 & ---    & ---    \\
			\hline
			$2^{-2}$ & $2^{-8}$  & 1.49418015e-03 & 0.2012 & 0.0503 \\
			\hline
			$2^{-3}$ & $2^{-12}$ & 3.71488610e-04 & 2.0080 & 0.5020 \\
			\hline
			$2^{-4}$ & $2^{-16}$ & 9.23212504e-05 & 2.0086 & 0.5021 \\
			\hline
		\end{tabular}
		\caption{Errors, spatial convergence rates, and time convergence rates for $k=h^4$ with $B(u)=u$ under Dirichlet boundary conditions.}
		\label{tab:3}
	\end{table}
	
	\begin{table}[htb]
		\renewcommand{\arraystretch}{1.2}
		\centering
		\begin{tabular}{|c|c|c|c|c|}
			\hline
			\hline
			$h$ & $k$ & $L_{\omega}^2L_t^2H^2_x(u)$ error & spatial rate & time rate \\
			\hline
			\hline
			$2^{-1}$ & $2^{-4}$  & 4.73065413e-03 & ---    & ---    \\
			\hline
			$2^{-2}$ & $2^{-8}$  & 4.72798046e-03 & 0.0008 & 0.0002 \\
			\hline
			$2^{-3}$ & $2^{-12}$ & 1.24395806e-03 & 1.9263 & 0.4816 \\
			\hline
			$2^{-4}$ & $2^{-16}$ & 2.58472530e-04 & 2.2669 & 0.5667 \\
			\hline
		\end{tabular}
		\caption{Errors, spatial convergence rates, and time convergence rates for $k=h^4$ with $B(u)=10\sin(u)$ under Dirichlet boundary conditions.}
		\label{tab:4}
	\end{table}

}

\section{Conclusion}\label{sec-5}

We derived and analyzed a fully discrete finite element approximation for the stochastic Kuramoto--Sivashinsky (SKS) equation with multiplicative noise. Our method combines a standard continuous finite element scheme in space with the implicit Euler--Maruyama method in time. For bounded multiplicative noise, we established optimal strong convergence rates in full expectation by exploiting a discrete stochastic Gronwall inequality, exponential stability estimates, and higher-moment bootstrapping arguments. For the more general case of unbounded multiplicative noise, we developed a localization strategy that yields sub-optimal convergence rates in probability, thereby extending the applicability of the method to a broader noise regime. 

The results presented here provide the first rigorous error analysis for fully discrete finite element schemes applied to the SKS equation. Not only does it offer a theoretical foundation for the numerical study of this challenging stochastic PDE, but it also highlights techniques, such as the interplay of stability estimates, stochastic Gronwall inequalities, and localization, that are likely to prove useful in the numerical study of other nonlinear SPDEs with non-Lipschitz drift and multiplicative noise. Future research directions include extending the analysis to higher spatial dimensions of the SKS equation and other nonlinear SPDEs.

{\section*{Acknowledgment} 
	
	The authors would like to thank anonymous referees for their providing a thorough review of this work. We appreciate their careful reading and insightful comments, which have improved the manuscript.}

\section*{Declarations of Funding}\,

The author Liet Vo was supported by the National Science Foundation (NSF) under Grant No. DMS-2530211.

\section*{Declarations of Author Contribution} \,

All authors contributed equally to this work. HDN was responsible for Section 2, including the proofs of the PDE results, such as the regularity of solutions and Holder continuity estimates. LV was responsible for Section 3, including the development of the main numerical scheme and the derivation of its stability and error estimates. He also prepared the abstract and introduction. All authors contributed to the reading, editing, and final approval of the manuscript.

\appendix

\section{Auxiliary results} \label{auxiliary}

First, we state the following discrete stochastic Gronwall inequality from \cite[Theorem 1]{kruse2018discrete}, which plays a vital role in our error analysis in Section \ref{sec3}.

\begin{lemma}{\cite[Theorem 1]{kruse2018discrete}}\label{Stochastic_Gronwall}
	Let $\left\{M_n\right\}_{n \in \mathbb{N}}$ be an $\left\{\mathcal{F}_n\right\}_{n \in \mathbb{N}}$-Martingale satisfying $M_0=0$ on a filtered probability space $\left(\Omega, \mathcal{F},\left\{\mathcal{F}_n\right\}_{n \in \mathbb{N}}, \mathbb{P}\right)$. Let $\left\{X_n\right\}_{n \in \mathbb{N}},\left\{F_n\right\}_{n \in \mathbb{N}}$, and $\left\{G_n\right\}_{n \in \mathbb{N}}$ be sequences of nonnegative and adapted random variables with $\mathbb{E}\left[X_0\right]<\infty$ such that
	\begin{align}\label{ineq2.3}
		X_n \leq F_n+M_n+\sum_{k=0}^{n-1} G_{\ell} X_{\ell} \quad \text { for all } n \in \mathbb{N}.
	\end{align}
	Then, for any $q \in(0,1)$ and a pair of conjugate numbers $\alpha, \beta \in[1, \infty]$, i.e.,  $\frac{1}{\alpha}+\frac{1}{\beta}=1$, satisfying $q \alpha<1$,  there holds 
	\begin{align}\label{ineq2.4}
		\mathbb{E}\left[\sup _{0 \leq \ell \leq n} X_{\ell}^q\right] \leq\left(1+\frac{1}{1-\alpha q}\right)^{\frac{1}{\alpha}}\left\|\prod_{\ell=0}^{n-1}\left(1+G_{\ell}\right)^q\right\|_{L^{\beta}(\Omega)}\left(\mathbb{E}\left[\sup _{0 \leq \ell \leq n} F_{\ell}\right]\right)^q.
	\end{align}
\end{lemma}

In the remainder of this section, we proceed to collect auxiliary results that are employed to establish the energy estimates in Section \ref{sec2}.	As mentioned in Section \ref{sec2}, the operator $\nu\partial_x^4+\partial_x^2$ is not guaranteed to be positive for arbitrary values of $\nu$. To establish Lemma \ref{lemma_expo_moment_H}, it is crucial to exploit the nonlinear nature of $u\partial_x u$ to navigate the sign-definite issue from the linear term. To this end, we will closely follow the framework of \cite{collet1993global,gao2025exponential,nicolaenko1985some}, tailored to our setting of multiplicative noise. More specifically, we introduce ${H}^{k}([0,2L])$ as the Sobolev spaces of periodic functions on $[0,2L]$ with vanishing integrals. We record the following inequality that was previously proven in \cite{gao2025exponential}.
\begin{lemma}{\cite[Lemma A.3]{gao2025exponential}}
	For all $\nu>0$, there exists a $2L$-periodic function $\varphi\in C^2(D)$ such that for all $u\in H^2(D)$ and $b\in{\mathbb R}$, it holds that 
	\begin{align} \label{ineq:phi}
		&\frac{1}{2}\nu \|\partial_x^2 u\|^2_{L^2([0,2L])}-\|\partial_x u\|^2_{L^2([0,2L])}+\frac{1}{2}( u^2,\partial_x\varphi(\cdot+b))_{L^2([0,2L])}  \notag  \\
		&\ge \frac{1}{8}\nu \|\partial_x^2u\|^2_{L^2([0,2L])} +\frac{1}{2}\|u\|^2_{L^2([0,2L])}-\frac{1}{4L}\big|( u,\partial_x\varphi(\cdot+b))_{L^2([0,2L])}  \big|^2.
	\end{align}
	
\end{lemma}
In the above, $( \cdot,\cdot)_{L^2([0,2L])}$ denotes the inner product in ${L^2([0,2L])}$
\begin{align*}
	( u,v)_{L^2([0,2L])} = \int_{0}^{2L} u(x)v(x)d x.
\end{align*}
Also, for $b(\cdot)\in C^1([0,\infty);\mathbb{R})$, we denote by $\varphi_{b(\cdot)}$ the translation of $\varphi$ defined as
\begin{align} \label{form:phi_b(t)}
	\varphi_{b(t)}(x) = \varphi(x+b(t)) ,\quad x\in\mathbb{R}.
\end{align}
In particular, given $u(t)$ the solution of \eqref{eq1.1}, we introduce the process $b(t)$ satisfying the equation
\begin{align}\label{eqn:b(t)}
	\frac{d}{d t}b(t) = \frac{1}{4L}\left( u,\partial_x\varphi_{b(t)}\right)_{L^2([0,2L])},\quad b(0)=0.
\end{align}
We note that since $u\in C([0,\infty);L^2(D))$, the solution $b(t)$ of \eqref{eqn:b(t)} is guaranteed to exist. See \cite[Appendix A]{collet1993global}. In particular, this auxiliary function is used to establish Lemma \ref{lemma_expo_moment_H}.

Next, we turn to the study of the process $z(t)$ defined in \eqref{eqn:z} that is employed to establish Lipschitz bounds in Lemma \ref{lemma_Holder}. We note that thanks to the choice of $a$ in \eqref{cond:a}, the linear operator $ \nu \partial_x^4 +\partial_x^2 +a$ is strictly positive. We will exploit this crucial property to derive useful estimates on $z(t)$ through Lemma \ref{lem_sup_|z|_Hn} and Lemma \ref{lemma_Holder_z} below. These results appear in the analysis of the solution $u(t)$ in Section \ref{sec2}. 

We start with Lemma \ref{lem_sup_|z|_Hn} giving a moment bound in $H^m$ norm on $z(t)$.

\begin{lemma}\label{lem_sup_|z|_Hn}
	Given $m\ge 0$ and $q\ge 2$, suppose that $u_0\in \cap_{i=0}^m L^{2^{m-i}q}(\Omega;H^i(D))$ and that $B$ satisfies condition \eqref{Assump_Lipschitz}. Then, the following holds
	\begin{align} \label{ineq:sup_|z|_Hn}
		\sup_{t\in[0,T]}  \E\left[\|z(t)\|_{H^m}^q\right] \le C,
	\end{align}
	for some positive constant $C=C(m,q,u_0,T)$. In the above, $z(t)$ is the process defined in \eqref{eqn:z}.
\end{lemma}
\begin{proof} We apply It\^o's formula to \eqref{eqn:z} and obtain
	\begin{align*}
		&d \|\partial_x^mz\|^q_{L^2} - q\|\partial_x^m z\|^{q-2}_{L^2}\|\partial_x^{m+1}z\|^2_{L^2}\, dt + q\nu\|\partial_x^m z\|^{q-2}_{L^2}\|\partial_x^{m+2}z\|^2_{L^2}\,dt\\
		&\quad=  q\|\partial_x^m z\|^{q-2}_{L^2}(\partial_x^m z,\partial_x^m B(u)dW) + \frac{1}{2}q\|\partial_x^m z\|^{m-2}_{L^2}\|\partial_x^m B(u)\|^{2}_{L^2}\,dt\\
		&\qquad+  \frac{1}{2}q(q-2)|\|\partial_x^m z\|^{q-4}_{L^2}(\partial_x^m z,\partial_x^m B(u))|^2\,dt.
	\end{align*}
	We employ integration by parts and H\"older's inequality to infer
	\begin{align*}
		-q\|\partial_x^m z\|^{q-2}_{L^2}\|\partial_x^{m+1}z\|^2_{L^2} = q\|\partial_x^m z\|^{q-2}_{L^2} (\partial_x^m z,\partial_x^{n+2} z)  \le \frac{\nu}{100}\|\partial_x^m z\|^{q-2}_{L^2}\|\partial_x^{m+2}z\|^2_{L^2} + C\|\partial_x^m z\|^q_{L^2}.
	\end{align*}
	Also, the Cauchy-Schwarz inequality produces
	\begin{align*}
		\frac{1}{2}q\|\partial_x^m z\|^{q-2}_{L^2}\|\partial_x^m B(u)\|^2_{L^2}&+  \frac{1}{2}q(q-2)|\|\partial_x^m z\|^{q-4}_{L^2}|(\partial_x^m z,\partial_x^m B(u))|^2\\
		&\le C\|\partial_x^m z\|^q_{L^2}+C\|\partial_x^m B(u)\|^q_{L^2}.
	\end{align*}
	As a consequence, since $B(u)$ is dominated by $u$ in $H^m$ norm, we obtain (recalling $z(0)=0$)
	\begin{align*}
		\E\left[\|\partial_x^m z(t)\|^q_{L^2}\right] & \le C\int_0^t \E\left[\|\partial_x^mz(s)\|^q_{L^2}\right]\, ds +C\int_0^t \E\left[\|\partial_x^m B(u(s))\|^q_{L^2}\right] \,ds \\
		&\le C\int_0^t \E\left[\|\partial_x^mz(s)\|^q_{L^2}\right]\,ds +C\int_0^t \E\left[\|\partial_x^m u(s)\|^q_{L^2}\right]\, ds+Ct.
	\end{align*}
	In light of Lemma \ref{lemma_poly_moment_H^s}, we have
	\begin{align*}
		\sup_{s\in[0,T]}\E\left[\|\partial_x^m u(s)\|^q_{L^2}\right]  \le C\sum_{i=0}^{m}\E\left[\|u_0\|^{2^{m-i}q}_{H^i}\right].
	\end{align*}
	In turn, Gronwall's inequality implies that 
	\begin{align*}
		\sup_{s\in[0,T]}\E\left[\|\partial_x^mz(s)\|^q_{L^2} \right] \le C\E\left[\|u_0\|^{2^{m-i}q}_{H^i}\right],
	\end{align*}
	which establishes \eqref{ineq:sup_|z|_Hn}, as claimed.
	
\end{proof}

Having obtained moment bounds on $z(t)$, we state and prove the following H\"older continuity property: 

\begin{lemma} \label{lemma_Holder_z}
	(a) Given an integer $m\ge 0$, suppose that $B$ satisfies condition \eqref{Assump_Lipschitz} and that 
	\begin{align*}
		u_0\in \cap_{i=0}^{m+1} L^{2^{m+2-i}}(\Omega;H^i(D)).
	\end{align*}
	Then, for any $0\le s\le t\le T$, the following holds
	\begin{align} \label{ineq:|z(t)-z(s)|^2_Hn}
		\E\left[\|\partial_x^m (z(t) -z(s)) \|_{L^2}^{2} \right]\le C(t-s).
	\end{align}
	for a positive constant $C=C(m,T,u_0)$.
	
	(b) Suppose further that given an integer $q\ge 2$,
	\begin{align*}
		u_0\in \cap_{i=0}^{(m+1)q} L^{2^{(m+1)q-i}2q}(\Omega;H^i(D)).    
	\end{align*}
	Then, for any $0\le s\le t\le T$, the following holds
	\begin{align} \label{ineq:|z(t)-z(s)|^q_Hn}
		\E\left[\|\partial_x^m (z(t) -z(s)) \|_{L^2}^{2q}\right] \le C(t-s)^{q}.
	\end{align}
	for a positive constant $C=C(m,q,T,u_0)$.
\end{lemma}

\begin{proof} 
	For notational convenience, we set $\rho=\frac{2\pi}{L}$ {and denote
		\begin{align*}
			e_{\ell}(x) = \frac{1}{\sqrt{L}}e^{\textup{i}\ell \rho x},\quad \ell \in\Zbb\setminus \{0\}.
		\end{align*}
		We note that $\{e_\ell\}_{\ell \in\Zbb\setminus \{0\}}$ forms an orthonormal basis in $L^2(D)$. So, given $z\in L^2(D)$, we may recast $z$ through the Fourier representation
		\begin{align*}
			z = \sum_{\ell \in\Zbb\setminus \{0\}} (z,e_\ell)e_\ell=: \sum_{\ell \in\Zbb\setminus \{0\}} z_\ell e_\ell.
		\end{align*}
		In particular, since $z(t)$ solves \eqref{eqn:z}, 
		$z_\ell(t)$ satisfies the equation
		\begin{align*}
			d z_\ell = -(\nu \rho^4\ell^4- \rho^2 \ell^2+a)z_\ell dt + (B(u),e_\ell)dW(t),
		\end{align*}
		with initial condition $z_\ell(0)=0$. It follows that
		\begin{align*}
			z_\ell(t) =  \int_0^t e^{-(\nu \rho^4\ell^4- \rho^2 \ell^2+a)(t-s)}(B(u(s)),e_{\ell})\,dW(s).
	\end{align*}}
	In turn, $z(t)$ admits the representation
	\begin{align*}
		z(t) =\sum_{\ell\in\Zbb \setminus\{0\}}  \int_0^t e^{-(\nu \rho^4\ell^4- \rho^2 \ell^2+a)(t-s)}(B(u(s)),e_{\ell})\,dW(s) \,e_{\ell}.
	\end{align*}
	{Since $\partial_x^m e_\ell = (\textup{i}\ell \rho)^m e_\ell$, we get
		\begin{align*}
			\partial_x^m z(t) &= \sum_{\ell\in\Zbb \setminus\{0\}} \int_0^t e^{-(\nu \rho^4\ell^4- \rho^2 \ell^2+a)(t-s)}(B(u(s)),e_{\ell})\,dW(s) \,\partial_x^m e_{\ell} \\
			&=  \sum_{\ell\in\Zbb \setminus\{0\}} (\textup{i}\ell \rho)^m\int_0^t e^{-(\nu \rho^4\ell^4- \rho^2 \ell^2+a)(t-s)}(B(u(s)),e_{\ell})\,dW(s) \,e_{\ell}  .
		\end{align*}
		Thus, for $0\le s\le t$, we have
		\begin{align*}
			&\|\partial_x^m z(t)-\partial_x^m z(s)\|^2_{L^2}\\
			& \quad =    \rho^{2m}\sum_{{\ell\in\Zbb \setminus\{0\}}}\ell^{2m}\Big|\int_0^t e^{-(\nu  \rho^4 \ell^4- \rho^2 \ell^2+a)(t-\xi)}(B(u(\xi)),e_{\ell})dW(\xi) \\
			&\hspace{4cm}- \int_0^s e^{-(\nu  \rho^4 \ell^4- \rho^2 \ell^2+a)(s-\xi)}(B(u(\xi)),e_{\ell})dW(\xi)  \Big|^2   \\
			&\quad = \rho^{2m}\sum_{{\ell\in\Zbb \setminus\{0\}}}\ell^{2m}\Big|\int_s^t e^{-(\nu  \rho^4 \ell^4- \rho^2 \ell^2+a)(t-\xi)}(B(u(\xi)),e_{\ell})dW(\xi) \\
			&\qquad+\big[e^{-(\nu   \rho^4 \ell^4 - \rho^2 \ell^2+a)(t-s)}-1\big]\int_0^s e^{-(\nu   \rho^4 \ell^4 - \rho^2 \ell^2+a)(s-\xi)}(B(u(\xi)),e_{\ell})dW(\xi) \Big|^2.
		\end{align*}
		Using the inequality $(a+b)^2 \le 2(a^2+b^2)$,} it follows immediately that
	\begin{align} \label{ineq:|partial_x^n_z(t)-z(s)|^2}
		&\|\partial_x^m z(t)-\partial_x^m z(s)\|^2_{L^2} \notag \\
		&\le 2\rho^{2m} \sum_{{\ell\in\Zbb \setminus\{0\}}}\ell^{2m}\Big|\int_s^t e^{-(\nu   \rho^4 \ell^4 - \rho^2 \ell^2+a)(t-\xi)}(B(u(\xi)),e_{\ell})dW(\xi)\Big|^2  \\
		&\quad+2\rho^{2m}\sum_{{\ell\in\Zbb \setminus\{0\}}}\ell^{2m}\Big|\big[e^{-(\nu   \rho^4 \ell^4 - \rho^2 \ell^2+a)(t-s)}-1\big]\notag \\
		&\qquad\qquad\qquad\qquad\times\int_0^s e^{-(\nu   \rho^4 \ell^4 - \rho^2 \ell^2+a)(s-\xi)}(B(u(\xi)),e_{\ell})dW(\xi) \Big|^2\notag \\
		& =: 2(I_1+I_2).\notag 
	\end{align}
	Concerning $I_1$, we employ It\^o's isometry and the fact that $\nu   \rho^4 \ell^4 - \rho^2 \ell^2+a>0$ to deduce
	\begin{align*}
		\E [I_1] & = \sum_{{\ell\in\Zbb \setminus\{0\}}}\rho^{2m}\ell^{2m}\int_s^t e^{-2(\nu   \rho^4 \ell^4 - \rho^2 \ell^2+a)(t-\xi)}\E|(B(u(\xi)),e_{\ell})|^2\,d\xi\\
		&\le  \int_s^t \E\left[\|B(u(\xi))\|^2_{H^m}\right]\,d\xi\\
		&\le (t-s)\sup_{\xi\in[0,T]}\E\Big[\|B(u(\xi))\|^2_{H^m}\Big].
	\end{align*}
	Likewise, we have
	\begin{align*}
		\E[I_2] & = \sum_{{\ell\in\Zbb \setminus\{0\}}}\rho^{2m}\ell^{2m}\big[e^{-(\nu   \rho^4 \ell^4 - \rho^2 \ell^2+a)(t-s)}-1\big]^2\\
		&\qquad\qquad\qquad\times\int_0^s e^{-2(\nu   \rho^4 \ell^4 - \rho^2 \ell^2+a)(s-\xi)}\E\left[|(B(u(\xi)),e_{\ell})|^2\right]\,d \xi.
	\end{align*}
	We employ the elementary inequality $1-e^{-x}\le \sqrt{x}$, $x\ge 0$ to infer
	\begin{align*}
		\big[e^{-(\nu  \rho^4 \ell^4- \rho^2 \ell^2+a)(t-s)} -1\big]^2 \le (\nu  \rho^4 \ell^4- \rho^2 \ell^2+a)(t-s),
	\end{align*}
	whence
	\begin{align*}
		& \big[e^{-(\nu  \rho^4 \ell^4- \rho^2 \ell^2+a)(t-s)}-1\big]^2\int_0^s e^{-2(\nu  \rho^4 \ell^4- \rho^2 \ell^2+a)(s-\xi)}\E|(B(u(\xi)),e_{\ell})|^2\,d \xi\\
		&\le (t-s) (\nu  \rho^4 \ell^4- \rho^2 \ell^2+a)\int_0^s e^{-2(\nu  \rho^4 \ell^4- \rho^2 \ell^2+a)(s-\xi)}\E|(B(u(\xi)),e_{\ell})|^2\,d \xi.
	\end{align*}
	We invoke the fact that $e^{-x}\le 1/\sqrt{x}$, $x>0$, to further estimate
	\begin{align*}
		&\int_0^s e^{-2(\nu  \rho^4 \ell^4- \rho^2 \ell^2+a)(s-\xi)}\E|(B(u(\xi)),e_{\ell})|^2\,d \xi \\
		&\le \frac{1}{\sqrt{2(\nu  \rho^4 \ell^4- \rho^2 \ell^2+a)}}\int_0^s \frac{1}{\sqrt{ s-\xi}}\E|(B(u(\xi)),e_{\ell})|^2\,d \xi.
	\end{align*}
	So, 
	\begin{align*}
		\E[I_2] & \le C(t-s)\int_0^s \frac{1}{\sqrt{s-\xi}} \sum_{{\ell\in\Zbb \setminus\{0\}}}\rho^{2m}\ell^{2m}\sqrt{\nu  \rho^4 \ell^4- \rho^2 \ell^2+a}\, \E\left[|(B(u(\xi)),e_{\ell})|^2\right]\,d \xi\\
		&\le C(t-s)\int_0^s \frac{1}{\sqrt{s-\xi}} \sum_{{\ell\in\Zbb \setminus\{0\}}}\rho^{2(m+1)}\ell^{2(m+1)} \E\left[|(B(u(\xi)),e_{\ell})|^2\right]\,d \xi\\
		&\le C (t-s) \int_0^s\frac{1}{\sqrt{s-\xi}}  \E\left[\|B(u(\xi))\|^2_{H^{m+1}}\right] \,d\xi\\
		&\le C(t-s)\sup_{\xi\in[0,T]}\E\Big[\|B(u(\xi))\|^2_{H^{m+1}}\Big].
	\end{align*}
	Altogether with \eqref{ineq:|partial_x^n_z(t)-z(s)|^2}, we get
	\begin{align*}
		\E \|\partial_x^m z(t)-\partial_x^m z(s)\|^2_{L^2}  \le C(t-s)\sup_{\xi\in[0,T]}\E\Big[\|B(u(\xi))\|^2_{H^{m+1}}\Big].
	\end{align*}
	Since $B(u)$ is dominated by $u$ in $H^{m+1}$ norm and that $u_0\in \cap_{i=0}^{m+1} L^{2^{m+2-i}}(\Omega;H^i(D))$, we invoke Lemma \ref{lemma_poly_moment_H^s}, part 1, cf. \eqref{ineq:poly_moment_H^s:int_0^T} with $q=2$, to infer
	\begin{align*}
		\sup_{\xi\in[0,T]}\E\Big[\|(B(u(\xi))\|^2_{H^{m+1}}\Big] \le C\left(1+\sup_{\xi\in[0,T]} \E\Big[\|\partial_x^{m+1} u(\xi)\|^2_{L^2}\Big]\right) \le C.
	\end{align*}
	As a consequence, we obtain
	\begin{align*}
		\E \|\partial_x^m z(t)-\partial_x^m z(s)\|^2_{L^2} \le C(t-s).
	\end{align*}
	This establishes \eqref{ineq:|z(t)-z(s)|^2_Hn}, thereby completing the proof of part (a).

	(b) With regard to \eqref{ineq:|z(t)-z(s)|^q_Hn}, we employ \eqref{ineq:|partial_x^n_z(t)-z(s)|^2} again to infer for $q> 1$
	\begin{align*}
		&\|\partial_x^m z(t)-\partial_x^m z(s)\|^{2q}_{L^2} \\
		&\le C\Big( \sum_{{\ell\in\Zbb \setminus\{0\}}}\ell^{2m}\Big|\int_s^t e^{-(\nu  \rho^4 \ell^4- \rho^2 \ell^2+a)(t-\xi)}(B(u(\xi)),e_{\ell})dW(\xi)\Big|^2\Big)^{q} \\
		&\qquad+C\Big(\sum_{{\ell\in\Zbb \setminus\{0\}}} \ell^{2m}\Big|\big[e^{-(\nu  \rho^4 \ell^4- \rho^2 \ell^2+a)(t-s)}-1\big]\\
		&\qquad\qquad\qquad\qquad\times\int_0^s e^{-(\nu  \rho^4 \ell^4- \rho^2 
			\ell^2+a)(s-\xi)}(B(u(\xi)),e_{\ell})dW(\xi) \Big|^2\Big)^{q}\\
		&=:C(J_1+J_2).
	\end{align*}
	With regard to $J_1$, since $q>1$, we use H\"older's inequality to infer
	\begin{align*}
		J_1&\le  \Big(\sum_{{\ell\in\Zbb \setminus\{0\}}}\ell^{(2m +1)q}\Big|\int_s^t e^{-(\nu  \rho^4 \ell^4- \rho^2 \ell^2+a)(t-\xi)}(B(u(\xi)),e_{\ell})dW(\xi)\Big|^{2q}\Big)\\
		&\qquad\qquad\times\Big(  \sum_{{\ell\in\Zbb \setminus\{0\}}} \frac{1}{\ell^{q/(q-1)}}\Big)^{q-1}.
	\end{align*}
	By virtue of Burkholder's inequality and the fact that $\nu  \rho^4 \ell^4- \rho^2 \ell^2+a>0$, we obtain the bound in expectation
	\begin{align*}
		\E[J_1] &\le C \sum_{{\ell\in\Zbb \setminus\{0\}}}\ell^{(2m+1)q}\E\left[\Big|\int_s^t e^{-2(\nu  \rho^4 \ell^4- \rho^2 \ell^2+a)(t-\xi)}|(B(u(\xi)),e_{\ell})|^2d\xi\Big|^{q}\right]\\
		&\le C \sum_{{\ell\in\Zbb \setminus\{0\}}}\ell^{(2m+1)q}\Big(\int_s^t e^{-2\frac{q}{q-1}(\nu  \rho^4 \ell^4- \rho^2 \ell^2+a)(t-\xi)}d \xi\Big)^{q-1}\E\left[\int_s^t |(B(u(\xi)),e_{\ell})|^{2q}\,d\xi\right]\\
		&\le C (t-s)^{q-1} \int_s^t \E\Big[|B(u(\xi))|^{2q-2}_{L^2}\|B(u(\xi))\|^2_{H^{(m+1)q}}\Big]d \xi\\
		&\le C(t-s)^q \sup_{\xi\in[0,T]} \E\left[\|B(u(\xi))\|^{2q}_{H^{(m+1)q}}\right].
	\end{align*}
	Turning to $J_2$, a similar argument to that of $J_1$ produces
	\begin{align*}
		\E [J_2] & \le C \sum_{{\ell\in\Zbb \setminus\{0\}}} \ell^{(2m+1)q}\big[e^{-(\nu  \rho^4 \ell^4- \rho^2 \ell^2+a)(t-s)}-1\big]^{2q}\\
		&\qquad\qquad\times\E\left[\Big|\int_0^s e^{-(\nu  \rho^4 \ell^4- \rho^2 \ell^2+a)(s-\xi)}(B(u(\xi)),e_{\ell}) dW(\xi) \Big|^{2q}\right].
	\end{align*}
	Similar to the argument of $I_2$ from part (a), we have
	\begin{align*}
		\big[e^{-(\nu  \rho^4 \ell^4- \rho^2 \ell^2+a)(t-s)}-1\big]^{2q} \le (\nu  \rho^4 \ell^4- \rho^2 \ell^2+a)^q(t-s)^q.
	\end{align*}
	On the other hand, we invoke Burkholder's and H\"older's inequalities  again to deduce
	\begin{align*}
		&\E\left[\Big|\int_0^s e^{-(\nu  \rho^4 \ell^4- \rho^2 \ell^2+a)(s-\xi)}(B(u(\xi)),e_{\ell})dW(\xi) \Big|^{2q}\right] \\
		&\quad \le  C\,\E\left[\Big|\int_0^s e^{-2(\nu  \rho^4 \ell^4- \rho^2 \ell^2+a)(s-\xi)}|(B(u(\xi)),e_{\ell})|^2\,d\xi \Big|^q\right]\\
		&\quad \le  C\Big(\int_0^s e^{-\frac{q}{q-1}(\nu  \rho^4 \ell^4- \rho^2 \ell^2+a)(s-\xi)}\,d\xi\Big)^{q-1}\int_0^s e^{-q(\nu  \rho^4 \ell^4- \rho^2 \ell^2+a)(s-\xi)}\E\left[|(B(u(\xi)),e_{\ell})|^{2q}\right]\,d\xi \\
		&  \quad \le  C \frac{1}{(\nu  \rho^4 \ell^4 - \rho^2 \ell^2+a)^{q-1}}\int_0^s e^{-q(\nu  \rho^4 \ell^4- \rho^2 \ell^2+a)(s-\xi)}\E\left[|(B(u(\xi)),e_\ell)|^{2q}\right]\,d\xi.
	\end{align*}
	It follows that
	\begin{align*}
		\E[J_2] & \le C (t-s)^q 
		\sum_{{\ell\in\Zbb \setminus\{0\}}} \ell^{(2m+1)q}(\nu  \rho^4 \ell^4- \rho^2 \ell^2+a)\\
		&\qquad\qquad\qquad\qquad\times\int_0^s e^{-q(\nu  \rho^4 \ell^4- \rho^2 \ell^2+a)(s-\xi)}\E\left[|(B(u(\xi)),e_{\ell})|^{2q}\right]\,d \xi\\
		&\le C(t-s)^q \int_0^s  \sum_{{\ell\in\Zbb \setminus\{0\}}} \ell^{(2m+1)q}\E\left[|(B(u(\xi)),e_{\ell})|^{2q}\right]\, d \xi\\
		&\le C(t-s)^{q} \sup_{\xi\in[0,T]} \E\left[\|B(u(\xi))\|^{2q}_{H^{(m+1)q}}\right].
	\end{align*}
	Altogether, we obtain
	\begin{align*}
		\E\left[\|\partial^m_x z(t)-\partial^m_x z(s)\|^{2q}_{L^2}\right] \le C(t-s)^{q} \sup_{\xi\in[0,T]} \E\left[\|B(u(\xi))\|^{2q}_{H^{(m+1)q}}\right].
	\end{align*}
	Since $u_0\in \cap_{i=0}^{(m+1)q} L^{2^{(m+1)q-i}2q}(\Omega;H^i(D))$, in light of Lemma \ref{lemma_poly_moment_H^s}, part 1, it holds that 
	\begin{align*}
		\sup_{\xi\in[0,T]} \E\left[\|B(u(\xi))\|^{2q}_{H^{(m+1)q}}\right] \le C\left(1+\sup_{\xi\in[0,T]} \E\left[\|u(\xi)\|^{2q}_{H^{(m+1)q}}\right] \right)\le C(m,q,T,u_0).
	\end{align*}
	In turn, we deduce
	\begin{align*}
		\E\left[\|\partial^m_x z(t)-\partial_x^m z(s)\|^{2q}_{L^2}\right] \le C(t-s)^q,
	\end{align*}
	which establishes \eqref{ineq:|z(t)-z(s)|^q_Hn}. The proof is thus complete.
\end{proof}

Lastly, we provide the proof of Lemma \ref{lemma_Holder_v}, while making use of auxiliary estimates on the process $z$ established in Lemma \ref{lem_sup_|z|_Hn}.

\begin{proof}[Proof of Lemma \ref{lemma_Holder_v}]
	First of all, a routine calculation produces the identity
	\begin{align} \label{eqn:d/dt|v(t)-v(s)|^2_Hn}
		0&=\frac{1}{2}\frac{d}{dt}\|\partial_x^m (v-v(s))\|^2_{L^2}- (\partial_x^{m+1}(v-v(s)),\partial_x^{m+1}v)+\nu (\partial_x^{m+2}(v-v(s)),\partial_x^{m+2}v)\notag \\
		&\qquad-a(\partial_x^{m}(v-v(s)),\partial_x^m z)+ (\partial_x^{m}(v-v(s)),\partial_x^m[(v+z)\partial_x(v+z)])\\
		&  =: \frac{1}{2}\frac{d}{dt}\|\partial_x^m(v-v(s))\|^2_{L^2}+I_1+\dots I_4.\notag 
	\end{align}
	Concerning $I_1+I_2+I_3$, we employ the Cauchy-Schwarz inequality with the Sobolev embedding to infer
	\begin{align*}
		|I_1+I_2+I_3| \le C\big(\|\partial_x^{m+2} v\|^2_{L^2}+ \|\partial_x^{m+2} v(s)\|^2_{L^2}+ \|\partial_x^{m} z\|^2_{L^2}\big).
	\end{align*}
	Turning to $I_4$, we compute 
	\begin{align*}
		I_4& = \big(\partial_x^{m}(v-v(s)),(v+z)\partial_x^{m+1}(v+z)\big)+\sum_{i=1}^m c_i\big(\partial_x^{m}(v-v(s)),\partial_x^i(v+z)\partial_x^{m+1-i}(v+z)\big)\\
		&= -\big(\partial_x^{m}(v-v(s)),\partial_x(v+z)\partial_x^{m}(v+z)\big)-\big(\partial_x^{m+1}(v-v(s)),(v+z)\partial_x^{m}(v+z)\big)\\
		&\qquad \qquad+\sum_{i=1}^m c_i\big(\partial_x^{m}(v-v(s)),\partial_x^i(v+z)\partial_x^{m+1-i}(v+z)\big).
	\end{align*}
	For $i=1,\dots,m$, it holds that
	\begin{align*}
		\big(\partial_x^{m} &(v-v(s)),\partial_x^i(v+z)\partial_x^{m+1-i}(v+z)\big) \\
		&\le \|\partial_x^{m}(v-v(s))\|_{L^\infty}\|\partial_x^i(v+z)\|_{L^2}\|\partial_x^{m+1-i}(v+z) \|_{L^2}\\
		&\le c\|\partial_x^{m+1}(v-v(s))\|_{L^2}\|\partial_x^m(v+z)\|^2_{L^2}\\
		&\le c\big(\|\partial_x^{m+1}v\|^2_{L^2}+ \|\partial_x^{m+1}v(s)\|^2_{L^2} + \|\partial_x^{m}v\|^4_{L^2}+\|\partial_x^{m+1}z\|^4_{L^2}\big).
	\end{align*}
	Likewise
	\begin{align*}
		\big(\partial_x^{m +1}&(v-v(s)),(v+z)\partial_x^{m}(v+z)\big) \\
		&\le \|\partial_x^{m+1}(v-v(s))\|_{L^\infty}\|(v+z)\|_{L^2}\|\partial_x^{m}(v+z) \|_{L^2}\\
		&\le c\big(\|\partial_x^{m+2}v\|^2_{L^2}+ \|\partial_x^{m+2}v(s)\|^2_{L^2} + \|\partial_x^{m}v\|^4_{L^2}+\|\partial_x^{m+1}z\|^4_{L^2}\big).  
	\end{align*}
	It follows that
	\begin{align*}
		I_4\le c\big(\|\partial_x^{m+2}v\|^2_{L^2}+ \|\partial_x^{m+2}v(s)\|^2_{L^2} + \|\partial_x^{m}v\|^4_{L^2}+\|\partial_x^{m+1}z\|^4_{L^2}\big).
	\end{align*}
	We collect the estimates on $I_1,\dots, I_4$ with expression \eqref{eqn:d/dt|v(t)-v(s)|^2_Hn} to arrive at the a.s. bound
	\begin{align} \label{ineq:|v(t)-v(s)|^2_Hn<int_s^t|u|+|z|}
		&\|\partial_x^m(v(t)-v(s))\|^2_{L^2} \notag \\
		&\le  c\int_s^t \big(\|\partial_x^{m+2}v(\xi)\|^2_{L^2}+ \|\partial_x^{m+2}v(s)\|^2_{L^2} + \|\partial_x^{m}v(\xi)\|^4_{L^2}+\|\partial_x^{m+1}z(\xi)\|^4_{L^2}+1\big) d\xi,
	\end{align}
	whence for $q\ge 1$,
	\begin{align*}
		&\|\partial_x^m(v(t)-v(s))\|^{2q}_{L^2} \\
		&\le  c\Big(\int_s^t \big(\|\partial_x^{m+2}v(\xi)\|^2_{L^2}+ \|\partial_x^{m+2}v(s)\|^2_{L^2} + \|\partial_x^{m}v(\xi)\|^4_{L^2}+\|\partial_x^{m+1}z(\xi)\|^4_{L^2}+1\big) d\xi\Big)^{q}\\
		&\le c(t-s)^{q-1}\Big(\int_s^t \big(\|\partial_x^{m+2}v(\xi)\|^{2q}_{L^2}+ \|\partial_x^{m+2}v(s)\|^{2q}_{L^2} + \|\partial_x^{m}v(\xi)\|^{4q}_{L^2}+\|\partial_x^{m+1}z(\xi)\|^{4q}_{L^2}+1\big) d\xi.
	\end{align*}
	Recalling $v=u-z$, we get
	\begin{align*}
		& \E \|\partial_x^m(v(t)-v(s))\|^{2q}_{L^2} \\
		&\le c(t-s)^q\sup_{\xi\in[0,T]}\big(\E\|\partial_x^{m+2}v(\xi)\|^{2q}_{L^2}+\E\|\partial_x^{m}v(\xi)\|^{4q}_{L^2}+\|\partial_x^{m+1}z(\xi)\|^{4q}_{L^2}+1\big)\\
		&\le c(t-s)\sup_{\xi\in[0,T]}\big(\E\|\partial_x^{m+2}u(\xi)\|^{2q}_{L^2}+\E\|\partial_x^{m+2}z(\xi)\|^{2q}_{L^2}+\E\|\partial_x^{m}u(\xi)\|^{4q}_{L^2}+\|\partial_x^{m+1}z(\xi)\|^{4q}_{L^2}+1\big).
	\end{align*}
	To determine the condition so that the above supremum is finite, we invoke Lemma \ref{lemma_poly_moment_H^s} part 1, cf. \eqref{ineq:poly_moment_H^s:int_0^T}, and observe that if
	\begin{align*}
		u_0 \in \bigcap_{i=0}^{m+2} L^{2^{m+2-i}2q}(\Ome;H^i(D))=\bigcap_{i=0}^{m+2} L^{2^{m+3-i}q}(\Ome;H^i(D)),
	\end{align*}
	it holds that
	\begin{align*}
		\sup_{\xi\in[0,T]}\E\|\partial_x^{m+2}u(\xi)\|^{2q}_{L^2} \le C(T,u_0,m,q).
	\end{align*}
	Also, if 
	\begin{align*}
		u_0 \in \bigcap_{i=0}^{m} L^{2^{m-i}4q}(\Ome;H^i(D))=\bigcap_{i=0}^{m} L^{2^{m+2-i}q}(\Ome;H^i(D)),
	\end{align*}
	we have
	\begin{align*}
		\sup_{\xi\in[0,T]}\E\|\partial_x^{m}u(\xi)\|^{4q}_{L^2} \le C(T,u_0,m,q).
	\end{align*}
	So, provided
	\begin{align*}
		u_0 \in \bigcap_{i=0}^{m+2} L^{2^{m+3-i}q}(\Ome;H^i(D)),
	\end{align*}
	we obtain
	\begin{align*}
		\sup_{\xi\in[0,T]}\big(\E\|\partial_x^{m+2}u(\xi)\|^{2q}_{L^2}+\E\|\partial_x^{m}u(\xi)\|^{4q}_{L^2}\big)\le C(T,u_0,m,q).
	\end{align*}
	Likewise, in view of Lemma \ref{lem_sup_|z|_Hn}, under the same condition on $u_0$, we deduce
	\begin{align*}
		\sup_{\xi\in[0,T]}\big(\E\|\partial_x^{m+2}z(\xi)\|^{2q}_{L^2}+\E\|\partial_x^{m}z(\xi)\|^{4q}_{L^2}\big)\le C(T,u_0,m,q).
	\end{align*}
	Altogether, we arrive at the bound in expectation
	\begin{align*}
		& \E \|\partial_x^m(v(t)-v(s))\|^{2q}_{L^2} \le C(t-s)^q,
	\end{align*}
	where $C=C(T,u_0,m,q)$ is independent of the difference $(t-s)$. This establishes \eqref{ineq:|v(t)-v(s)|^(2q)_Hn}, thereby completing the proof.

\end{proof}
	
	
	\bibliographystyle{abbrv}
	\bibliography{references}

\end{document}